\documentclass[11pt]{thesis}
\usepackage{upgreek}
\usepackage{graphicx}
\usepackage{amssymb}
\usepackage{epstopdf}
\usepackage{textcomp}
\usepackage{amsmath}
\usepackage{amsthm}
\usepackage{mathdots}
\usepackage[all,cmtip]{xy}
\usepackage[usenames,dvipsnames]{color}
\input diagxy
\usepackage{cancel}
\usepackage{url}

\title{A Geometric Construction of Cyclic Cocycles on Twisted Convolution Algebras}
\author{Eitan}{Angel}
\advisor{Associate Professor}{Alexander Gorokhovsky}
\degree{Doctor of Philosophy}{Ph.D., Mathematics}
\dept{Department of}{Mathematics}
\otherdegrees{
B.A., University of California at Berkeley, 2003}
\reader{Prof.~Arlan Ramsay}		

\abstract{  \OnePageChapter	
In this thesis we give a construction of cyclic cocycles on convolution algebras twisted by gerbes over discrete translation groupoids.  In his seminal book, Connes constructs a map from the equivariant cohomology of a manifold carrying the action of a discrete group into the periodic cyclic cohomology of the associated convolution algebra. Furthermore, for proper \'etale groupoids, J.-L.~Tu and P.~Xu provide a map between the periodic cyclic cohomology of a gerbe twisted convolution algebra and twisted cohomology groups.

Our focus will be the convolution algebra with a product defined by a gerbe over a discrete translation groupoid.  When the action is not proper, we cannot construct an invariant connection on the gerbe; therefore to study this algebra, we instead develop simplicial notions related to ideas of J.~Dupont to construct a simplicial form representing the Dixmier-Douady class of the gerbe.  Then by using a JLO formula we define a morphism from a simplicial complex twisted by this simplicial Dixmier-Douady form to the mixed bicomplex of certain matrix algebras. Finally, we define a morphism from this complex to the mixed bicomplex computing the periodic cyclic cohomology of the twisted convolution algebras. 
}

\dedication[Dedication]{	
	To my father who sowed the seeds of my mathematical journey and to Beethoven for the joy he provided during this recent leg.
	}

\acknowledgements{	\OnePageChapter	
	I would first like to thank my advisor, Professor Alexander Gorokhovsky, for his patience in guiding me through this project.
	I would also like to thank Professor Judy Packer for her support and terrific instruction, Professor Siye Wu from whom learning mathematics was inspirational, Professor Arlan Ramsay for his helpful suggestions in producing this thesis, as well as Professors Carla Farsi, Markus Pflaum, Marty Walter, and my friends in the Mathematics Department at CU for creating a wonderful community I have been proud to be a part of.
	Finally, I would like to thank my family for their wholehearted support.
	}

\ToCisShort	

	\emptyLoF	

	\emptyLoT	

\begin{document}

\newcommand{\A}{\mathcal{A}}
\newcommand{\U}{\mathcal{U}}
\newcommand{\G}{\mathcal{G}}
\newcommand{\Ca}{\mathcal{C}}
\newcommand{\La}{\mathcal{L}}
\newcommand{\Gl}[1][{}]{\mathcal{G}_{(#1)}}
\newcommand{\Cl}[1][{}]{\mathcal{C}_{(#1)}}
\newcommand{\N}{\mathbb{N}}
\newcommand{\Z}{\mathbb{Z}}
\newcommand{\R}{\mathbb{R}}
\newcommand{\C}{\mathbb{C}}
\newcommand{\id}{\mathrm{id}}
\newcommand{\E}{\mathcal{E}}
\newcommand{\End}{\mathrm{End}\,}
\newcommand{\tr}{\mathrm{tr}}
\newcommand{\Hom}{\mathrm{Hom}\,}
\newcommand{\Log}{\mathrm{Log}\,}
\newcommand{\sgn}{\mathrm{sgn}\,}

\newcommand{\set}[2]{ \left\{ #1 \,  : #2 \right\} }

\theoremstyle{plain}
\newtheorem{thm}{Theorem}[chapter]
\newtheorem{cor}[thm]{Corollary}
\newtheorem{prop}[thm]{Proposition}
\newtheorem{lem}[thm]{Lemma}

\theoremstyle{definition}
\newtheorem{defn}[thm]{Definition}
\newtheorem{ex}[thm]{Example}

\theoremstyle{remark}
\newtheorem{rmk}[thm]{Remark}

\chapter{Introduction}


In the early 1980's, Alain Connes described diverse motivations for studying various non-commutative algebras and presented a framework for the subject now known as noncommutative geometry in \cite{connes94}. One way to study such non-commutative algebras geometrically is by means of their cyclic cohomology.  A smooth manifold $M$ may be studied through its algebra of smooth functions $C^\infty_c (M)$.  In \cite{connes85}, Connes identifies the periodic cyclic cohomology of $C^\infty_c (M)$ with the de Rham cohomology of $M$.  This shows that cyclic cohomology can be seen as the appropriate generalization of de Rham cohomology to non-commutative algebras.

If, in addition, $M$ carries the action of a discrete group $\Gamma$, a natural question to ask is what the orbit space $M / \Gamma$ looks like.  In general the orbit space is not a manifold or even a Hausdorff space. 
The prescription of Connes' noncommutative geometry to deal with such ``bad'' quotients is to instead study the convolution algebra of smooth functions with compact support on the translation groupoid $M \rtimes \Gamma$.  This is the algebra, $C^\infty_c (M \rtimes \Gamma)$, of smooth functions with compact support on $M \times \Gamma$ with the discrete convolution product
$$(f_1 * f_2) (x,g) = \sum_{g_1 g_2 = g} f_1 (x, g_1) f_2 (x g_1 , g_2)$$
which in general defines a noncommutative algebra.  Furthermore, in \cite{connes86} and \cite{connes94} Chapter III.2, Connes describes the cyclic cohomology of the convolution algebra $C^\infty_c (M \rtimes \Gamma)$ via a morphism $\Phi$ from a complex computing the de Rham homology of the homotopy quotient $M_\Gamma = M \times_\Gamma E \Gamma$, i.e.\ the Borel model of equivariant cohomology $H^*_\Gamma (M;\C) = H^* (M_\Gamma;\C)$, to the $(b,B)$-bicomplex computing the periodic cyclic cohomology of $C^\infty_c (M \rtimes \Gamma)$, $HP^*(C^\infty_c (M \rtimes \Gamma))$.

The main object of study in this thesis is the twisted convolution algebra $C^\infty_c (M \rtimes \Gamma, L)$ and the periodic cyclic cohomology $HP^* ( C^\infty_c (M \rtimes \Gamma, L))$ defined as follows.  Given a line bundle $L \rightarrow M \rtimes \Gamma$ as well as a collection of line bundle isomorphisms $\mu_{g_1 , g_2} : L_{g_1} \otimes (L_{g_2})^{g_1} \xrightarrow{\sim} L_{g_1 g_2}$ for every $g_1 , g_2 \in \Gamma$, where $L_g$ denotes the restriction of $L$ to $M_g = M \times \{g\}$, the sections $C^\infty (M \rtimes \Gamma, L)$ define a twisted convolution algebra with a convolution product defined by
\begin{equation*}
(f_1 * f_2) (x,g) = \sum_{g_1 g_2 = g}  \mu_{g_1 , g_2} (f_1|_{M_{g_1}} \otimes (f_2|_{M_{g_2}})^{g_1} ) (x , g).
\end{equation*}
The data $(L , \mu)$ is a gerbe over the translation groupoid $M \rtimes \Gamma$.

Gerbes were introduced by J. Giraud \cite{giraud71} and largely developed in a differential geometric framework by J.-L. Brylinski \cite{brylinski93} to study central extensions of loop groups, line bundles on loop spaces, and the Dirac monopole among other applications.  Just as line bundles may be seen as a geometric realization of classes in $H^2 (M, \Z)$, gerbes on manifolds may be seen as a geometric realization of classes in $H^3(M,\Z)$.  In \cite{block_daenzer}, the authors associate a certain type of differential graded algebra, called a curved dga, to a gerbe on an \'etale groupoid.  Furthermore, considerable attention has been paid to gerbes recently due to their appearance in string theory and quantum field theory.  
For instance, in string theory, the $B$-field may be viewed as a gerbe 
as made explicit by D. Freed and E. Witten in \cite{freed_witten99}.  The appearance of gerbes in string theory and quantum field theory has been studied by numerous mathematicians and physicists, e.g., in \cite{bcmms02}, \cite{sss09}, \cite{gawedzki_reis02} and \cite{cmm00}.

The main result of this thesis is the construction of a morphism analogous to Connes' $\Phi$ map into the $(b,B)$-bicomplex computing the periodic cyclic cohomology of a twisted convolution algebra, $HP^* ( C^\infty_c (M \rtimes \Gamma, L))$, where the twisting is determined by a gerbe $(L, \mu)$ on the translation groupoid $M \rtimes \Gamma$.  Such a morphism, given by a JLO-type formula, was constructed by J.-L. Tu and P. Xu in \cite{tu_xu06} in the case of a gerbe over any proper \'etale groupoid.  The construction of this thesis removes the properness requirement in the case of a discrete translation groupoid.  In our situation one cannot construct a $\Gamma$-invariant connection on $L \rightarrow M \rtimes \Gamma$ as is done in \cite{tu_xu06}.  In particular, this means that we cannot in general choose a connection $\nabla$ on $L \rightarrow M \rtimes \Gamma$ that satisfies the cocycle condition
\begin{equation*}
\mu_{g,h}^* (\nabla_{gh}) = \nabla_g \otimes 1 + 1 \otimes (\nabla_h)^g
\end{equation*}
where $\nabla_g$ denotes $\nabla |_{M_g}$.  Instead there is some discrepancy $\alpha \in \Omega^1 ( (M \rtimes \Gamma_{(2)} )$
\begin{equation*}
\nabla_g \otimes 1 + 1 \otimes (\nabla_h)^g - \mu_{g,h}^* (\nabla_{g.h}) = \alpha (g,h).
\end{equation*}
The curvature forms of $\nabla_g$, $\theta_g \in \Omega^2 (M)$, define a form $\theta \in \Omega^2 ( (M \rtimes \Gamma)_{(1)} )$.  The $\theta_g$ satisfy
\begin{equation*}
\theta_g -\theta_{gh}+ \theta_h^g = d\alpha(g,h)
\end{equation*}
so that $(-\alpha, \theta)$ is a cocycle that represents the Dixmier-Douady class in $H^3 ( (M \rtimes \Gamma)_\bullet )$, the degree 3 de Rham cohomology of the groupoid $M \rtimes \Gamma$.

To overcome this difficulty, the geometric data which forms the source of our morphism is constructed using simplicial notions similar to those of J. Dupont in \cite{dupont}.  Given a simplicial manifold $M_\bullet$, Dupont defines a complex of compatible forms 
$(\Omega^* (M_\bullet) , d_{dR} + d_\Delta)$.  Elements of $\Omega^k (M_\bullet)$ are sequences of $k$-forms $\{\omega_{(n)} \}$, where $\omega_{(n)} \in \Omega^k (M_n \times \Delta^n)$, and which satisfy the compatibility condition
\begin{equation*}
(\id \times \partial_i)^* \omega_{(n)} = (\delta_i \times \id)^* \omega_{(n-1)}
\end{equation*}
on $\Omega^k (M_n \times \Delta^{n-1})$ for $0 \leq i \leq n$ and $n \geq 1$.  Here $\delta_i$ denotes the $i^{\text{th}}$ face map on $M_n$ and $\partial_i$ denotes the $i^{\text{th}}$ face map on the geometric simplex $\Delta^{n-1}$. Also $d_{dR}$ denotes the exterior differential on the simplicial manifold $M_\bullet$ and $d_\Delta$ denotes the exterior differential on the geometric simplex.

The geometric data we use is a complex termed the \textit{twisted simplicial complex} 
\begin{equation*}
(\Omega^*_- ( (M \rtimes \Gamma)_\bullet )[u] , u \tilde{d}_{dR} + d_\Delta - \Theta_u \wedge \cdot)
\end{equation*}
where $\Omega^*_- (M_\bullet)$ is a rescaled version of the compatible forms $\Omega^* (M_\bullet)$ and $u$ is a formal variable.  The twisting is by a simplicial version of the Dixmier-Douady form $\Theta_u$ in $\Omega^3_- ((M \rtimes \Gamma)_\bullet)$.  The 3-form $\Theta_u$ is constructed using the data of the gerbe $(L, \mu)$ and additionally a choice of connection $\nabla$ on $L \rightarrow M \rtimes \Gamma$.  From $L \rightarrow M \rtimes \Gamma$ we define an infinite dimensional vector bundle $\E \rightarrow M$ as the direct sum bundle $\E = \bigoplus_{g \in \Gamma} L_g$ with the direct sum connection $\nabla^\E$.  Then $\End \E \rightarrow M$ is a vector bundle of finite rank endomorphisms over $M$.


From the data $(L, \mu)$ we define connections $\nabla^k$ on the vector bundles $\E_k \times \Delta^k \rightarrow (M \rtimes \Gamma)_{(k)} \times \Delta^k$, where $\E_k = p_k^* \E$ and $p_k : M \times \Gamma^k \rightarrow M$ is the map $p_k : (m, g_1 , \ldots , g_k) \mapsto m \in M$.  While $\End \E_k \rightarrow (M \rtimes \Gamma)_{(k)}$ forms a simplicial vector bundle, the sequence of $\End \E_k$-valued 2-forms defined by $(\nabla^k)^2$ do not define a simplicial form.  Instead we are able to define a simplicial 2-form by the formula
\begin{equation*}
\vartheta_{(k)} (g_1 , \ldots , g_k) = (\nabla^k)^2 - \sum_{i=1}^k t_i \theta_{g_1 \cdots g_i}  + \sum_{1 \leq i < j \leq k} \alpha(g_1 \cdots g_i , g_{i+1} \cdots g_j) (t_i dt_j - t_j dt_i).
\end{equation*}
Here the difference between $\vartheta_{(k)}$ and $(\nabla^k)^2$ is a scalar valued form.  The simplicial Dixmier-Douady form $\Theta_u$ is defined by $(\Theta_u)_{(k)} = \nabla^k_u (\vartheta_u)_{(k)}$, where $\nabla^k_u$ and $\vartheta_u$ are rescaled versions of $\nabla^k$ and $\vartheta$.

A simplicial version of the character formula of A. Jaffe, A. Lesniewski, and K. Osterwalder \cite{jlo88}, similar to \cite{mathai_stevenson06} and \cite{gorokhovsky99}, given by 
\begin{multline*}
\tau_\nabla (  \omega ) (\tilde{a}_0 , a_1 , \ldots , a_n) = \sum_k \int_M \int_{\Delta^k}  \\ \omega_{(k)} \wedge \left( \int_{\Delta^n} \tr \left(\tilde{a}_0 e^{-\sigma_0 \left( \vartheta_u \right)_{(k)}} \nabla_u^k (a_1) e^{-\sigma_1  \left( \vartheta_u \right)_{(k)}} \cdots \nabla_u^k (a_n) e^{-\sigma_n  \left( \vartheta_u \right)_{(k)}} \right) d \sigma_1 \cdots d \sigma_n \right) 
\end{multline*}
for $\omega \in \Omega^*_- ( (M \rtimes \Gamma)_\bullet)$ and $a_i \in  C^\infty_c (M , \End \E) $, defines a morphism into group cochains valued in the $(b,B)$ bicomplex of sections of $\End \E$,
\begin{multline*}
\tau_\nabla : (\Omega^*_- (M_\bullet) [u] , u\tilde{d}_{dR} + d_\Delta - \Theta_u \wedge \cdot) \\ \to (C^\bullet (\Gamma , \overline{C}^\bullet ( C^\infty_c (M , \End \E) )  [u^{-1} , u] ) , b+uB + {\delta_\Gamma}').
\end{multline*}
Following this we define algebraic morphisms
\begin{multline*}
(C^\bullet (\Gamma , \overline{C}^\bullet ( C^\infty_c (M , \End \E) )  [u^{-1} , u] ) , b+uB + {\delta_\Gamma}') \\ \to \overline{C}^\bullet (C^\infty_c (M \rtimes \Gamma , L)  [u^{-1} , u], b+uB)
\end{multline*}
into the periodic cyclic complex of the twisted convolution algebra $C^\infty_c (M \rtimes \Gamma , L)$.

The rest of this thesis is organized as follows.

The content of Chapter 2 is standard material on simplicial cohomology, cyclic cohomology, and gerbes.

In Chapter 3 we discuss gerbes on the translation groupoid and we define the twisted simplicial complex which is the source of the morphism.

Chapter 4 describes the construction of the simplicial Dixmier-Douady form $\Theta_u$.

Finally, Chapter 5 contains the main result, a morphism from the twisted simplicial complex of a gerbe $(L, \mu)$ on $M \rtimes \Gamma$ to the periodic cyclic complex of the twisted convolution algebra $C^\infty_c (M \rtimes \Gamma, L)$.


\chapter{Preliminaries}
\label{chap:1}


\section{Categorical notions}
\label{sec:cat}

First we will introduce some standard notions from category theory.  This material is largely taken from \cite{connes94} and \cite{loday98}.


\subsection{The simplicial category}
\begin{defn}
\label{def:simp_cat}
The \textit{simplicial category} $\Delta$ is the small category of objects $[n]$, where $[n]$ is the ordered set of $n+1$ points, $[n] = \{0 < 1 < \cdots < n \}$, $n \in \N$, and arrows the nondecreasing maps $f: [n] \rightarrow [m]$.  For $n,m \in \N$, a map $f : [n] \rightarrow [m]$ is called \textit{nondecreasing} if $f(i) \geq f(j)$ whenever $i > j$.
\end{defn}

There are certain morphisms in $\Delta$ called face maps and degeneracy maps that we focus on, as any other morphism may be written in terms of face maps and degeneracy maps.

\begin{defn}
\label{def:face_deg}
In $\Delta$ the \textit{face maps} are the injections $\updelta^n_i : [n-1] \rightarrow [n]$, $0 \leq i \leq n$, which skip over the $i^\mathrm{th}$ point, i.e.\ such that $\updelta^n_i (i-1) = i-1$, $\updelta^n_i (i) = i+1$.  The \textit{degeneracy maps} are the surjections $\upsigma^n_j : [n+1] \rightarrow [n]$, $0\leq j \leq n$, that sends both $j$ and $j+1$ to $j$, i.e.\ such that $\upsigma^n_j (j) = \upsigma^n_j (j+1) = j$. When the context is clear, we will omit the superscript $n$.
\end{defn}

\begin{prop}
\label{prop:face_deg_comm}
The category $\Delta$ is generated by the morphisms $\updelta^n_i$ and $\upsigma^n_j$ and has the presentation
\begin{align}
\label{eq:face_deg_comm}
\updelta_j \updelta_i &= \updelta_i \updelta_{j-1} \qquad \textrm{for } i<j \\
\upsigma_j \upsigma_i &= \upsigma_i \upsigma_{j+1} \qquad \textrm{for } i \leq j \\
\upsigma_j \updelta_i &= \begin{cases} \updelta_i \upsigma_{j-1} &\textrm{for } i<j \\ \id_{[n]} &\textrm{if } i=j \textrm{ or } i=j+1 \\ \updelta_{i-1} \upsigma_j &\textrm{for } i>j+1 \end{cases} 
\end{align}
\begin{proof}
We can check that the above relations are obeyed by direct computation. We also need to show that this generates all of the morphisms in $\Delta$.
\end{proof}
\end{prop}

\begin{defn}
A \textit{simplicial object} (respectively \textit{cosimplicial object}) in a category $\Ca$ is a contravariant functor $X_\bullet : \Delta \rightarrow \Ca$ (respectively a covariant functor $X^\bullet : \Delta \rightarrow \Ca$).  We will often denote the objects $X_n=X([n])$, $n \in \N$, and the morphisms $\delta_i^n= X_\bullet (\updelta^n_i)$ and $\sigma^n_j=X_\bullet (\upsigma^n_j)$, $0 \leq i , j \leq n$. In view of Proposition \ref{prop:face_deg_comm} , such a functor is determined entirely by the objects $X_n$ and the morphisms $d^n_i$ and $s^n_j$. When the morphisms $s^n_j$ are not specified, the functor is instead called a \textit{pre-simplicial object}.
\end{defn}

\begin{ex} Some examples of simplicial objects are given by the following.
\begin{enumerate}
\item A contravariant functor $X_\bullet :\Delta \rightarrow \mathbf{Sets}$ to the category of sets is called a simplicial set.
\item A contravariant functor $X_\bullet : \Delta \rightarrow \mathbf{Top}$ to the category of topological spaces is called a simplicial space.
\item A contravariant functor $X_\bullet : \Delta \rightarrow \mathbf{Man}$ to the category of smooth manifolds is called a simplicial manifold. This example will be discussed further in Section \ref{sec:Simplicial_manifolds}.
\end{enumerate}
\end{ex}

\begin{defn}
The \textit{geometric $n$-simplex} is the cosimplicial space $\mathbf{\Delta}^\bullet : \Delta \rightarrow \mathbf{Top}$ defined by $\mathbf{\Delta}^\bullet ([n]) = \Delta^n$, where $\Delta^n$ is the standard $n$-simplex
\begin{align}
\Delta^n &= \set{ (t_0 , \ldots , t_n ) \in \R^{n+1} }{ 0 \leq t_i \leq 1, \ 0 \leq i \leq n, \sum_{i =0}^n t_i = 1} & \textrm{(barycentric)} \\
&= \set{(t_1 , \ldots , t_n) \in \R^{n} }{ t_1 , \ldots , t_n \geq 0, \ \sum_{i=1}^n t_i \leq 1 } & \textrm{(Cartesian).}
\end{align}
The Cartesian coordinates may be obtained from the barycentric coordinates eliminating $t_0 = 1 - t_1 - \cdots - t_n$. For barycentric coordinates
\begin{align}
\mathbf{\Delta}^\bullet (\updelta^n_i) (t_0 , \ldots , t_{n-1}) &= (t_0 , \ldots , t_{i-1} , 0 , t_i , \ldots , t_{n-1}) \\
\mathbf{\Delta}^\bullet (\upsigma^n_j) (t_0 , \ldots ,  t_{n+1}) &=  (t_0 , \ldots , t_{j-1} , t_j + t_{j+1}, \ldots , t_{n+1}),
\end{align}
and for Cartesian coordinates
\begin{align}
\mathbf{\Delta}^\bullet (\updelta^n_i) (t_1 , \ldots , t_{n-1}) &= 
\begin{cases}
(1-t_1 - \cdots - t_{n-1}, t_1 , \ldots , t_{n-1}) & i = 0 \\ 
(t_1 , \ldots , t_{i-1} , 0 , t_i , \ldots , t_{n-1}) & 1 \leq i \leq n
\end{cases} \\
\mathbf{\Delta}^\bullet (\upsigma^n_j) (t_1 , \ldots ,  t_{n+1}) &= 
\begin{cases}
( t_2, \ldots , t_{n+1}) & i=0 \\
(t_1 , \ldots , t_j + t_{j+1} , \ldots , t_{n+1}) & 1 \leq i \leq n
\end{cases}
\end{align}
We will denote $\partial^n_i= \mathbf{\Delta}^\bullet (\updelta^n_i)$ and $\varsigma^n_j=\mathbf{\Delta}^\bullet (\upsigma^n_j)$.
\end{defn}

\begin{defn}
Given a simplicial set (or space or manifold) $X_\bullet$, the \textit{fat realization} of $X_\bullet$, denoted $\| X_\bullet \|$, is the space
\begin{equation}
\| X_\bullet \| = \bigcup_{n \geq 0} X_n \times \Delta^n / \sim
\end{equation}
where $\sim$ is the equivalence relation generated by $(\delta^n_i (x) , t) \sim ( x , \partial^n_i (t))$, for $(x,t) \in X_n \times \Delta^{n-1}$. The \textit{geometric realization} of $X$, denoted by $| X |$ is the space
\begin{equation}
| X_\bullet | = \bigcup_{n \geq 0} X_n \times \Delta^n / \approx
\end{equation}
where $\approx$ is the equivalence relation generated by $(f^* (x) , t) \approx (x , f_* (t))$ for $(x,t) \in X_n \times \Delta^{n-1}$ and any $f \in \hom (\Delta)$.
\end{defn}

\begin{defn}
\label{def:nerve}
Let $\Cl[n]$ consist of $n$-tuples of composable arrows in a (small) category $\Ca$, with $\Cl[0]$ the objects of $\Ca$.  The \textit{nerve} of $\Ca$ is the simplicial set $\Ca_\bullet : \Delta \rightarrow \mathbf{Sets}$ given on objects by $\Ca_\bullet ([n]) = \Cl[n]$ and on morphisms by the faces $\Ca_\bullet (\updelta^n_i) = \delta^n_i$ which compose adjacent morphisms in the $i^\mathrm{th}$ place and by the degeneracies $\Ca_\bullet (\upsigma^n_j) = \sigma^n_j$ which insert the identity morphism in the $j^\mathrm{th}$ place, i.e.\
\begin{align}
\label{eq:nerve}
\delta^n_i (f_1 , \ldots , f_n) &= \begin{cases} (f_2 , \ldots , f_n) & \textrm{if } i=0 \\ (f_1 , \ldots , f_{i} f_{i+1} , \ldots , f_n) & \textrm{for } 1 \leq i \leq n-1 \\ (f_1 , \ldots , f_{n-1}) & \textrm{if } i=n \end{cases} \\
\sigma^n_j (f_1 , \ldots , f_n) &= (f_1 , \ldots , f_j , \id , f_{j+1} , \ldots , f_n)
\end{align}
with $\delta^1_0 (f)$ the terminal object of $f$ and $\delta^1_1 (f)$ the initial object of $f$ and where $f_i f_{i+1} = f_{i+1} \circ f_i$.  The \textit{classifying space} $B \Ca$ of the small category $\Ca$ is the geometric realization of the nerve of $\Ca$.
\end{defn}

\begin{ex}
\begin{enumerate}
\item A discrete group $\Gamma$ may be considered as a category with only one object $*$, with morphisms the elements of $\Gamma$ and composition equal to the group product. The nerve of $\Gamma$ is denoted by $\Gamma_\bullet$ and is given by $\Gamma_{(n)} = \Gamma^n$ and
\begin{align*}
\delta^n_i (g_1 , \ldots , g_n) &=
\begin{cases}
(g_2 , \ldots , g_n) & \text{for $i=0$}, \\
(g_1 , \ldots , g_i g_{i+1} , \ldots , g_n) & \text{for $1\leq i \leq n-1$} , \\
(g_1 , \ldots , g_{n-1}) & \text{for $i=n$}
\end{cases} \\
s^n_j (g_1 , \ldots , g_n) &= (g_1 , \ldots , g_j , 1_\Gamma , g_{j+1} , \ldots , g_n).
\end{align*}
\item A topological $G$-space $X$ where $G$ is a topological group (possibly with the discrete topology), i.e.\ a space $X$ with an action of $G$, $X \times G \to X$, can be considered as a category. The objects are elements of $X$ and every $g \in G$ such that $xg = x'$ determines a morphism from $x \in X$ to $x' \in X$. The nerve of this category is a simplicial space with objects $X \times G^n$ and its geometric realization is the Borel space $X \times_G EG$.
\end{enumerate}
\end{ex}

\section{\'Etale groupoids}
\label{sec:gpd}


We will recall some standard material about groupoids.  See \cite{moerdijk02}, \cite{crainic99}, or \cite{moerdijk_mrcun} for example.

\begin{defn}
A \textit{groupoid} is a (small) category $\G$ in which every arrow is invertible. The set of objects is denoted by $\Gl[0]$ and the set of arrows is denoted by $\Gl[1]$.  The set of arrows will often be denoted simply by $\G$. For each arrow in $\Gl[1]$ there is a source object and a range object given by the \textit{range and source maps}, $r,s : \Gl[1] \rightrightarrows \Gl[0]$.  To denote that $g \in \Gl[1]$ is an arrow with source $s(g) = x$ and range $r(g) = y$ we write either $g : x \rightarrow y$ or $x \xrightarrow{g} y$.

As $\G$ is a category, there is a rule of composition.  Given two arrows $y \xrightarrow{g_1} z, x \xrightarrow{g_2} y \in \Gl[1]$ such that $s(g_1) = r(g_2)$, their composition is an arrow $x \xrightarrow{g_1 g_2} z$ called their \textit{product}.  If we define pairs of composable arrows as
\begin{equation}
\Gl[2] = \Gl[1] \times_{\Gl[0]} \Gl[1] = \{ (g_1 , g_2) \in \Gl[1] \times \Gl[1] : s(g_1) = r(g_2) \}
\end{equation}
then the product defines the \textit{multiplication map}
\begin{equation}
m : \Gl[2] \rightarrow \Gl[1], \qquad m(g_1 , g_2) = g_1g_2.
\end{equation}
As the composition in a category is required to be associative, the groupoid product is associative.

For any object $x \in \Gl[0]$ there is an identity morphism $1_x : x \rightarrow x$ that satisfies $1_x g = g 1_y = g$ for any arrow $x \xrightarrow{g} y \in \Gl[1]$. The \textit{unit map} is then $$u : \Gl[0] \rightarrow \Gl[1], \qquad u(x) = 1_x.$$

Every arrow $g : x \rightarrow y$ in $\Gl[1]$ is invertible so we will denote the inverse of $g$ by $g^{-1} : y \rightarrow x$.  With this we can define the \textit{inverse map} $$i : \Gl[1] \rightarrow \Gl[1], \qquad i(g) = g^{-1}.$$ Then $g^{-1} g = 1_x$ and $g g^{-1} = 1_y$.
\end{defn}

\begin{rmk}
As a groupoid $\G$ is a category, definition \ref{def:nerve} applies and there is a simplicial set $\G_\bullet$ called the nerve of $\G$.
\end{rmk}

\begin{defn}
As in \cite{bgnt07}, Section 2.3, we will define the following projection maps. For $0 \leq i \leq n$ and $\G_\bullet$ the nerve of $\G$ let $\varpi^n_i : \Gl[n] \rightarrow \Gl[0]$ be the initial object of the first morphism when $i=0$ and the final object of the $i^\mathrm{th}$ morphism (which is the initial object of the $(i-1)^\mathrm{st}$ morphism) when $i \neq 0$. For $0 \leq j \leq n$, $0 \leq i_j \leq m$, the map $\varpi^n_{i_0} \times \cdots \times \varpi^n_{i_m} : \Gl[n] \rightarrow (\Gl[0])^m$ factors into maps $\Gl[n] \rightarrow \Gl[m]$ and the canonical projection $\Gl[m] \rightarrow (\Gl[0])^m$. The map $\Gl[n] \rightarrow \Gl[m]$ will be denoted $\varpi^n_{i_0 \cdots i_m}$. We may write $\varpi^n_{i_0 \cdots i_m}$ explicitly as
\begin{equation*}
\varpi^n_{i_0 \cdots i_m} : (g_1 , \ldots, g_n) \mapsto ( g_1 ' , \ldots , g_m ')
\end{equation*}
where, for $1 \leq k \leq m$
\begin{equation*}
g_k ' = \begin{cases}
g_{i_{k-1} +1} \cdots g_{i_{k}} & \text{if $i_{k-1} < i_{k}$} \\
\id_{s(g_{i_k})} & \text{if $i_{k-1} = i_{k}$} \\
(g_{i_{k}} \cdots g_{i_{k-1} - 1})^{-1} & \text{if $i_{k-1} > i_{k}$}
\end{cases}
\end{equation*}
\end{defn}

\begin{rmk}
The maps $\delta^n_i$ of Definition \ref{eq:nerve} may be written $\varpi^n_{0\cdots \hat{i} \cdots n} : \Gl[n] \rightarrow \Gl[n-1]$ where $\hat{i}$ denotes the omission of $i$.
\end{rmk}

\begin{defn}
A \textit{Lie groupoid} is a groupoid in which the sets $\Gl[0]$ and $\Gl[1]$ are smooth manifolds and the structure maps $r,s,u,i,m$ are smooth.  Furthermore, the range and source maps $r,s : \Gl[1] \rightrightarrows \Gl[0]$ are required to be submersions so that the domain of the multiplication map, $\Gl[2]$, is a manifold.  An \textit{\'etale groupoid} is a Lie groupoid in which the source map is \'etale, i.e.\ a local diffeomorphism.  In this case the other structure maps are \'etale as well.
\end{defn}

\begin{ex} We will now give some examples of groupoids.
\begin{enumerate}
\item Any group $G$ may be viewed as a groupoid by taking the object set to be $\Gl[0] = pt$ and the arrows to be the group elements, $\Gl[1] = G$. The composition in the groupoid is the product in the group. If $G$ is a Lie group then $\G$ is a Lie groupoid.
\item Given an manifold $M$ with an open cover $\mathcal{U} = \{U_\alpha\}$ the \textit{\v{C}ech groupoid} associated to $\mathcal{U}$ has objects $\{U_\alpha\}$ and morphisms $U_{\alpha \beta} = U_\alpha \cap U_\beta$.  The range and source maps are given by $s(U_{\alpha \beta}) = U_\alpha$ and $r(U_{\alpha \beta}) = U_\beta$.
\item Given a manifold $M$ carrying a right action of a group $G$, there is a groupoid called the \textit{translation groupoid} or \textit{action groupoid} $M \rtimes G$ with objects $(M \rtimes G)_{(0)} = M$ and arrows $(M \rtimes G)_{(1)} = M \times G$. The source and range maps are given by $s(x,g) = xg$ and $r(x,g) = x$. Given arrows $xg_1 \xleftarrow{g_1} x$ and $xg_1 g_2 \xleftarrow{g_2} xg_1$ their composition is $m((x,g_1) ,(xg_1, g_2)) = (x, g_1 g_2) = xg_1 g_2 \xleftarrow{g_1 g_2} x$.
\item For a connected manifold $M$, the \textit{fundamental groupoid} $\Pi (M)$ of $M$ is the groupoid with objects $\Pi (M)_{(0)} = M$ and arrows $x \rightarrow y$ in $\Pi (M)_{(1)}$ are homotopy classes classes of paths from $x$ to $y$ for any points $x,y \in M$.  
\item An orbifold may be considered as a proper \'etale groupoid in the sense of \cite{moerdijk02}.  To any orbifold one may associate a proper \'etale groupoid $\G$ called a \textit{presentation} of the orbifold.  While a presentation of an orbifold is not unique, the Morita equivalence classes of proper \'etale groupoids are in one-to-one correspondence with orbifolds.

%

\item Given a foliation $\mathcal{F}$ of a manifold $M$, there is an associated \'etale groupoid called the \textit{holonomy groupoid}, $\mathrm{Hol} (M, \mathcal{F})$, which is well defined up to Morita equivalence.  See \cite{haefliger} for a discussion.


\end{enumerate}
\end{ex}

\begin{defn}
Let $\Omega^*(M)$ be the graded algebra of differential forms on a manifold $M$ and let $d$ be the de Rham differential. Denote by $\delta^n$ the alternating sum of the face maps of \eqref{eq:nerve}.  Then we can consider $\Omega^* (\G_\bullet)$ as a double complex with the de Rham differential $d : \Omega^k (\Gl[n]) \rightarrow \Omega^{k+1} (\Gl[n])$ and $\delta^* : \Omega^k (\Gl[n]) \rightarrow \Omega^k (\Gl[n+1])$, i.e.,
\begin{equation}
\label{eq:gpd_de_rham_cplx}
\xymatrix{
\vdots & \vdots & \vdots &  \\
\Omega^2 (\Gl[0]) \ar[r]^{\delta^*} \ar[u]_{d} & \Omega^2 (\Gl[1]) \ar[r]^{\delta^*} \ar[u]_{d} & \Omega^2 (\Gl[2]) \ar[r]^{\delta^*} \ar[u]_{d} & \cdots \qquad \\
\Omega^1 (\Gl[0]) \ar[r]^{\delta^*} \ar[u]_{d} & \Omega^1 (\Gl[1]) \ar[r]^{\delta^*} \ar[u]_{d} & \Omega^1 (\Gl[2]) \ar[r]^{\delta^*} \ar[u]_{d} & \cdots \qquad  \\
\Omega^0 (\Gl[0]) \ar[r]^{\delta^*} \ar[u]_{d} & \Omega^0 (\Gl[1]) \ar[r]^{\delta^*} \ar[u]_{d} & \Omega^0 (\Gl[2]) \ar[r]^{\delta^*} \ar[u]_{d} & \cdots \qquad
}
\end{equation}
The cohomology of the double complex $(\Omega^* (\G_\bullet) , d, \delta^*)$ is the \textit{de Rham cohomology} of $\G$, denoted $H^*_{dR} (\G)$.
\end{defn}


\section{Simplicial manifolds}
\label{sec:Simplicial_manifolds}

We will recount the ideas of Dupont in \cite{dupont} (and \cite{felisatti_neumann}) to explicitly construct the cohomology of a simplicial manifold.

\subsection{Simplicial forms}

\begin{defn}
\label{def:compatible_form}
Given a simplicial manifold $M_\bullet : \Delta \rightarrow \mathbf{Man}$, a \textit{simplicial differential $k$-form} on $M_\bullet$ is a sequence of $k$-forms $\{\omega_{(n)} \}$, where $\omega_{(n)} \in \Omega^k (M_n \times \Delta^n)$, and which satisfy the compatibility condition
\begin{equation}
\label{eq:compatibility_cond}
(\id \times \partial_i)^* \omega_{(n)} = (\delta_i \times \id)^* \omega_{(n-1)}
\end{equation}
on $\Omega^k (M_n \times \Delta^{n-1})$ for $0 \leq i \leq n$ and $n \geq 1$.
\begin{itemize}
\item The collection of simplicial differential forms on $M_\bullet$ will be denoted $\Omega^* (M_\bullet)$.
\item The exterior differential $d : \Omega^k (M_n \times \Delta^n) \rightarrow \Omega^{k+1} (M_n \times \Delta^n)$ induces an exterior differential $d : \Omega^k (M_\bullet) \rightarrow \Omega^{k+1} (M_\bullet)$.
\item The complex $(\Omega^* (M_\bullet) , d)$ is called the complex of \textit{compatible forms}. 
\item The wedge product $\wedge : \Omega^k (M_n \times \Delta^n) \otimes \Omega^l (M_n \times \Delta^n) \rightarrow \Omega^{k+l} (M_n \times \Delta^n)$ also induces a wedge product $\wedge : \Omega^k (M_\bullet) \otimes \Omega^l (M_\bullet) \rightarrow \Omega^{k+l} (M_\bullet)$.
\end{itemize}
\end{defn}

\begin{rmk}
Notice that $\{\omega_{(n)} \}$ defines a $k$-form on $\displaystyle \coprod_{n \geq 0} M_n \times \Delta^n$.
In view of this, the compatibility condition \eqref{eq:compatibility_cond} is precisely the condition required for $\{ \omega_{(n)} \}$ to define a form on $\|M_\bullet\|$.

The complex $(\Omega^* (M_\bullet), d)$ can actually be thought of as the total complex of a certain double complex $(\Omega^{r,s} (M_\bullet) , d_{dR} , d_\Delta)$ called the \textit{bicomplex of compatible forms}.  To define this double complex, let
\begin{equation}
\label{eq:bicompatible_form}
 \Omega^{r,s} (M_\bullet) = \left( \coprod_{n \geq 0} \Omega^r (M_n) \otimes \Omega^s (\Delta^n)  \right) / \sim
 \end{equation}
Here $\sim$ is similar to the compatibility condition \eqref{eq:compatibility_cond}.  If $\omega_{(n)}^{r,s} \in \Omega^r (M_n) \otimes \Omega^s (\Delta^n)$ we may write $\omega_{(n)}^{r,s}$ as a linear combination of forms $\alpha_{(n)}^r \otimes \beta_{(n)}^s$ where $\alpha_{(n)}^r \in \Omega^r (M_n)$ and $\beta_{(n)}^s \in \Omega^s (\Delta^n)$.  So $\omega_{(n)}^{r,s} = \alpha_{(n)}^r \otimes \beta_{(n)}^s \sim \alpha_{(n-1)}^r \otimes \beta_{(n-1)}^s = \omega_{(n-1)}^{r,s}$ if and only if 
\begin{equation}
\label{eq:bicompatibility_cond}
\alpha_{(n)}^r \otimes \partial_i^* (\beta_{(n)}^s) = \delta_i^* (\alpha_{(n-1)}^r ) \otimes \beta_{(n-1)}^s.
\end{equation}
We define the differentials by
\begin{itemize}
\item $d_{dR} = d : \Omega^{r,s} (M_\bullet) \rightarrow \Omega^{r+1,s} (M_\bullet)$ is the exterior differential on each $M_n$ induced by $d \otimes \id : \Omega^r (M_n) \otimes \Omega^s (\Delta^n) \rightarrow \Omega^{r+1} (M_n) \otimes \Omega^s (\Delta^n)$.
\item $d_\Delta = (-1)^s d : \Omega^{r,s} (M_\bullet) \rightarrow \Omega^{r,s+1} (M_\bullet)$ is the exterior differential on each $\Delta^n$ induced by $\id \otimes d : \Omega^r (M_n) \otimes \Omega^s (\Delta^n) \rightarrow \Omega^r (M_n) \otimes \Omega^{s+1} (\Delta^n)$.
\end{itemize}
The induced differentials on $\Omega^{r,s} (M_\bullet)$ are well defined due to the compatibility condition \eqref{eq:bicompatibility_cond}.  Hence for $r+s=k$ we can alternatively define 
\begin{equation}
\label{eq:bicompatible_sum}
\Omega^k (M_\bullet) = \bigoplus_{r+s=k} \Omega^{r,s} (M_\bullet).
\end{equation}
and $(\Omega^* (M_\bullet) , d_{dR} + d_\Delta)$ to be the complex of compatible forms.  Given an element of the complex of of compatible forms $\omega \in \Omega^k (M_\bullet)$ we will denote the component of $\omega$ in $\Omega^{r,s} (M_\bullet)$ by $\omega^{r,s} = \omega |_{\Omega^{r,s} (M_\bullet)}$.

The space $\Omega^{r,s} (M_\bullet)$ of $(r+s)$-forms are locally of the form
\begin{equation}
\omega^{r,s} |_{M_n \times \Delta^n} = \sum \omega_{i_1 \cdots i_r j_1 \cdots j_s} dx_{i_1} \wedge \cdots \wedge dx_{i_r} \wedge dt_{j_1} \wedge \cdots \wedge dt_{j_s}
\end{equation}
where $\{x_i\}$ are local coordinates of $M_n$ and $(t_0 , \ldots , t_n)$ are barycentric coordinates of $\Delta^n$.  Again, $d_{dR}$
is the exterior differential with respect to the local coordinates of $M_n$ and $d_\Delta$ is $(-1)^s$ times the exterior derivative with respect to the barycentric coordinates of $\Delta^n$.  
\end{rmk}

\begin{defn}
Dupont defines a complex called the \textit{simplicial de Rham complex} $(\A^* (M_\bullet) , \delta+d)$. Define
\begin{equation}
\A^k (M_\bullet) = \bigoplus_{r+s = k} \A^{r,s} (M_\bullet)
\end{equation}
where $\A^{r,s} (M_\bullet) = \Omega^s (M_r)$ is the collection of differential $s$-forms on the manifold $M_r$.  Also, the differential
\begin{equation}
\delta : \A^{r,s} (M_\bullet) \rightarrow \A^{r+1,s} (M_\bullet)
\end{equation}
is the alternating sum of the pullback of the face maps $\delta^{r+1}_i = M_\bullet (\updelta^n_i)$, $0 \leq i \leq r+1$ and 
\begin{equation}
d : \A^{r,s} (M_\bullet) \rightarrow \A^{r,s+1} (M_\bullet)
\end{equation}
is the exterior differential $d : \Omega^s (M_r) \rightarrow \Omega^{s+1} (M_r)$. Then $(\A^* (M_\bullet) , \delta+d)$ is the total complex of the double complex $(\A^{r,s} (M_\bullet), d, \delta)$.
\end{defn}

\begin{rmk}
If a simplicial manifold $\G_\bullet$ is the nerve of a groupoid $\G$ then the simplicial de Rham bicomplex $(\A^* (\G_\bullet) ,d, \delta)$ is the same as the de Rham complex of the groupoid $\G$, $(\Omega^* (\G_\bullet) , d, \delta^*)$, of Definition \ref{eq:gpd_de_rham_cplx}.
\end{rmk}

\begin{defn}
Given a commutative ring $R$ we can associate the \textit{simplicial singular cochain complex} $(C^*(M_\bullet ; R) , \delta + \partial)$ to a simplicial manifold $M_\bullet$.  Define
\begin{equation}
C^k (M_\bullet ; R) = \bigoplus_{r+s=k} C^{r,s} (M_\bullet ; R)
\end{equation}
where $C^{r,s} (M_\bullet ; R) = C^s (M_r ; R)$ is the collection of singular cochains of degree $s$ on $M_r$.  Also, the differential $\partial = \partial_1 + \partial_2$ where
\begin{equation}
\delta : C^{r,s} (M_\bullet; R) \rightarrow C^{r+1,s} (M_\bullet; R)
\end{equation}
is the alternating sum of the pullback of the face maps $\delta^{r+1}_i = M_\bullet (\updelta^n_i)$, $0 \leq i \leq r+1$ and
\begin{equation}
\partial : C^{r,s} (M_\bullet; R) \rightarrow C^{r,s+1} (M_\bullet; R)
\end{equation}
is the usual coboundary on singular cochains.  Then $(C^* (M_\bullet ; R) , \delta + \partial)$ is the total complex of the double complex $(C^{r,s} (M_\bullet ; R) , \delta , \partial)$.
\end{defn}

From \cite{dupont} Proposition 5.15
\begin{thm}
For a simplicial manifold $M_\bullet$, we have the isomorphism
\begin{equation}
H^* (\| M_\bullet \|; R) \cong H(C^* (M_\bullet ; R) , \delta + \partial)
\end{equation}
\end{thm}

From \cite{dupont} Proposition 6.1,
\begin{thm}[Simplicial de Rham theorem]
The integration map
\begin{equation}
\mathcal{I} : \A^{r,s} (M_\bullet) \rightarrow C^{r,s} (M_\bullet)
\end{equation}
defined by $\mathcal{I} (\omega^{r,s}) (c_{r,s}) = \int_{c_{r,s}} \omega^{r,s}$ for $\omega^{r,s} \in \A^{r,s} (M_\bullet)$ and $c_{r,s} \in C_s (M_r)$, the collection of singular $s$-chains on $M_r$, gives a morphism of double complexes. Furthermore, this integration map induces an isomorphism
\begin{equation}
H(\A^* (M_\bullet) , \delta + d) \cong H(C^* (M_\bullet) , \delta + \partial)
\end{equation}
on cohomology.
\end{thm}

Furthermore, there is another morphism of complexes given by Stokes' theorem. From \cite{dupont} Theorem 6.4,
\begin{thm}
\label{thm:simplex_integration}
Let 
\begin{equation}
\mathcal{I}_\Delta : (\Omega^{r,s} (M_\bullet) , d_{dR} , d_\Delta) \rightarrow (\A^{r,s} (M_\bullet) , \delta , d)
\end{equation}
be the map defined by integration over the standard simplex, i.e. the map defined on $\Omega^s (M_r \times \Delta^r)$ by
\begin{equation}
\mathcal{I}_\Delta : \omega_{(r)} \mapsto \int_{\Delta^r} \omega_{(r)}.
\end{equation}
Then $\mathcal{I}_\Delta$ is a morphism that induces an isomorphism
\begin{equation}
H(\Omega^* (M_\bullet) , d) \cong H((\A^* (M_\bullet) , \delta+d)
\end{equation}

\end{thm}

\subsection{Connections and Curvature}

\begin{defn}
Let $G$ be a Lie group.  A \textit{simplicial $G$-bundle} $\pi_\bullet : E_\bullet \rightarrow M_\bullet$ over a simplicial manifold $M_\bullet$ is a sequence of principal $G$-bundles $\{ \pi_n : E_n \rightarrow M_n\}$ where $E_\bullet$ is itself a simplicial manifold and the diagrams
\begin{equation}
\label{eq:simp_bund_diag}
\xymatrix{
E_{n+1}  \ar[rr]^{E_{n+1} (\updelta_i)} \ar[d]_{\pi_{n+1}} && E_{n} \ar[d]^{\pi_{n}} && E_{n-1} \ar[rr]^{E_{n-1} (\upsigma_j)} \ar[d]_{\pi_{n-1}} && E_{n} \ar[d]^{\pi_{n}} \\
M_{n+1} \ar[rr]^{M_{n+1} (\updelta_i)} && M_{n} & \textrm{and} & M_{n-1} \ar[rr]^{M_{n-1} (\upsigma_j)} && M_{n}
}
\end{equation}
commute.  Given a simplicial $G$-bundle $\pi_\bullet : E_\bullet \rightarrow M_\bullet$, the geometric realization $|\pi_\bullet| :| E_\bullet | \rightarrow | M_\bullet|$ is a principal $G$-bundle with $G$-action induced by
\begin{equation}
E_n \times \Delta^n \times G \rightarrow E_n \times \Delta^n, \quad (x,t,g) \mapsto (xg, t).
\end{equation}
If we only require that the first diagram of \eqref{eq:simp_bund_diag} commute then we may still consider the fat realization $\|\pi_\bullet\| :\| E_\bullet \| \rightarrow \| M_\bullet \|$ which is a principal $G$-bundle.
\end{defn}

\begin{defn}
A \textit{connection} in a simplicial $G$-bundle $\pi_\bullet : E_\bullet \rightarrow M_\bullet$ is a $1$-form $\omega \in \Omega^1 (E_\bullet ; \mathfrak{g})$ on $E_\bullet$ (in the sense of definition \ref{def:compatible_form}) with coefficients in $\mathfrak{g}$, the Lie algebra of $G$, such that $\omega_{(n)} = \omega|_{E_n \times \Delta^n}$ is a connection in the usual sense on the bundle $\pi_n \times \id : E_n \times \Delta^n \rightarrow M_n \times \Delta^n$. The curvature $\theta$ of a connection $\omega$ is the differential form
\begin{equation}
\theta = d\omega + \frac{1}{2} [\omega , \omega] \in \Omega^2 (M_\bullet ; \mathfrak{g}).
\end{equation}

\end{defn}


\section{Cyclic cohomology}

\begin{defn}
The \textit{cyclic category} $\Lambda$ has objects $[n]$, $n \in \N$, and morphisms generated by face maps $\updelta^n_i$ and degeneracy maps $\upsigma^n_j$, $0 \leq i,j \leq n$, as in definitions \ref{def:simp_cat} and \ref{def:face_deg} subject to the relations of \eqref{eq:face_deg_comm} as well as cyclic operators $\uptau_n : [n] \rightarrow [n]$ further satisfying the relations
\begin{align}
\updelta_i \uptau_n &= \begin{cases} \uptau_{n-1} \updelta_{i-1} & \textrm{if } 1 \leq i \leq n \\ \updelta_n & \textrm{if } i= 0 \end{cases} \\
\upsigma_j \uptau_n &= \begin{cases} \uptau_{n+1} \upsigma_{j-1} & \textrm{if } 1 \leq i \leq n \\ \uptau_{n+1}^2 \upsigma_n & \textrm{if } i=0 \end{cases} \\
\uptau_n^{n+1} &= \id_{[n]}
\end{align}
A \textit{cyclic object} (respectively, \textit{cocyclic object}) in a category $\Ca$ is a contravariant (respectively, covariant) functor $X_\bullet : \Lambda \rightarrow \Ca$ (respectively $X^\bullet : \Lambda \rightarrow \Ca$).  For a cyclic object we will typically denote $X_\bullet ([n]) = X_n$, $X_\bullet (\updelta_i^n) = d_i^n$, $X_\bullet (\upsigma_j^n) = s_j^n$, and $X_\bullet (\uptau_n) = \tau_n$.  For a cocyclic object we will typically denote $X^\bullet ([n]) = X^n$, $X^\bullet (\updelta^n_i) = d^i_n$, $X^\bullet (\upsigma_j^n) = s^j_n$, and $X^\bullet (\uptau_n) = \tau_n$.
\end{defn}

\begin{ex}
Let $\A$ be a unital algebra over a field $k$. One may associate to $\A$ the cyclic object $\A^\natural_\bullet : \Lambda \rightarrow k-\mathbf{Mod}$ in the category of $k$-modules.  On objects we define $\A^\natural_\bullet ([n]) = \A^\natural_n = \A^{\otimes (n+1)}$ and on morphisms $\A^\natural_\bullet (\updelta^n_i) = d^n_i$, $\A^\natural_\bullet (\upsigma^n_i) = s^n_i$, and $\A^\natural_\bullet (\uptau_n) = \tau_n$ where
\begin{align}
d_i (a_0 \otimes a_1 \otimes \cdots \otimes a_n) &= \begin{cases} a_0 \otimes a_1 \otimes \cdots \otimes a_i a_{i+1} \otimes \cdots \otimes a_n & \textrm{if } 0 \leq i \leq n-1 \\ a_n a_0 \otimes a_1 \otimes \cdots \otimes a_{n-1} & \textrm{if } i=n \end{cases} \\
s_i (a_0 \otimes a_1 \otimes \cdots \otimes a_n) &= a_0 \otimes a_1 \otimes \cdots \otimes a_i \otimes 1 \otimes a_{i+1} \otimes \cdots \otimes a_n \\
\tau (a_0 \otimes a_1 \otimes \cdots \otimes a_n) &= a_n \otimes a_0 \otimes a_1 \otimes \cdots \otimes a_{n-1}
\end{align}
Dually, one may associate a cocyclic object $\A_\sharp^\bullet : \Lambda \rightarrow k-\mathbf{Mod}$ to $\A$.  On objects we define $\A_\sharp^n = \Hom (\A^{\otimes (n+1)} , k)$ and on morphisms $\A_\sharp^\bullet(\updelta^n_i) = d_n^i$, $\A_\sharp^\bullet (\upsigma^n_i) = s_n^i$, and $\A_\sharp^\bullet (\uptau_n) = \tau^n$ where
\begin{align}
d^i f (a_0 , a_1 , \ldots , a_n) &= \begin{cases} f(a_0 , a_1 , \ldots , a_i a_{i+1} , \ldots , a_n) & \textrm{if } 0 \leq i \leq n-1 \\ f(a_n a_0 , a_1 , \ldots , a_{n-1} ) & \textrm{if } i=n \end{cases} \\
s^i f (a_0 , a_1 , \ldots , a_n) &= f(a_0 , a_1, \ldots , a_i , 1 , a_{i+1}, \ldots , a_n) \\
\tau f(a_0, a_1, \ldots , a_n) &= f(a_n , a_0, a_1 , \ldots , a_{n-1})
\end{align}
for any $f \in \Hom (\A^{\otimes (n+1)} , k)$.  We will denote $\Hom (\A^{\otimes (n+1)})$ by $C^n (\A)$.
\end{ex}

We are now in a position to define the cyclic (co)homology of a (co)cyclic object in an abelian category.  Let $X_\bullet$ be a cyclic object in an abelian category. Let $N_n : X_n \rightarrow X_n$ be defined by $N_n = 1 + \tau_n + \tau_n^2 + \cdots + \tau_n^n$.  There are operators $b_n: X_n \rightarrow X_{n-1}$ and $B_n : X_n \rightarrow X_{n+1}$ defined by
\begin{align}
b_n &= \sum_{i=0}^n d_i^n \\
B_n &= (1-\tau_n) s_n N_n
\end{align}
that satisfy the relations $b_n^2 = B_n^2 = b_n B_n + B_n b_n = 0$.  The subscripts will typically be omitted.
Define
\begin{equation}
C_{pq} (X_\bullet) = \begin{cases} X_{q-p} & \textrm{if } q \geq p \\ 0 & \textrm{otherwise} \end{cases}
\end{equation}
where $p,q \in \N$.  Then $b : C_{pq} \rightarrow C_{p (q-1)}$ and $B: C_{pq} \rightarrow C_{(p-1)q}$.  The homology of the bicomplex $\mathcal{B} (X_\bullet)= (C_{pq} , b, B)$ is the \textit{cyclic homology} of $X_\bullet$ which is denoted $HC_* (X_\bullet)$.  [If we instead allow $p,q \in \Z$ then the homology of the bicomplex $\mathcal{B}^\textrm{per} (X_\bullet) = (C_{pq} , b, B)$ is the \textit{periodic cyclic homology} of $X_\bullet$ which is denoted $HP_* (X_\bullet)$.]  The cyclic cohomology of a cocyclic object $X^\bullet$ is similarly defined. Now let $N^n = 1+ \tau^n + (\tau^n)^2 + \cdots + (\tau^n)^n$ and $b^n : X^n \rightarrow X^{n+1}$ and $B^n : X^n \rightarrow X^{n-1}$ be defined by
\begin{align}
b^n &= \sum_{i=0}^n d^i_n \label{eq:Hochschild_coboundary}\\
B^n &= (1-\tau^n) s^n N^n \label{eq:Connes_coboundary}
\end{align}
and define
\begin{equation}
C^{pq} (X^\bullet) = \begin{cases} X^{q-p} & \textrm{if } q \geq p \\ 0 & \textrm{otherwise} \end{cases}
\end{equation}
for $p,q \in \N$.  Then the cohomology of the \textit{mixed bicomplex} $\mathcal{B} (X^\bullet) = (C^{pq} , b , B)$ is the \textit{cyclic cohomology} of $X^\bullet$ which is denoted $HC^* (X^\bullet)$.  If we instead allow $p,q \in \Z$ then the cohomology of the bicomplex $\mathcal{B}^\textrm{per} (X^\bullet) = (C^{pq} , b, B)$ is the \textit{periodic cyclic cohomology} of $X^\bullet$, which is denoted $HP^* (X^\bullet)$.

In particular, for a unital algebra $\A$ over a field $k$, the bicomplex $\mathcal{B} (\A) = \mathcal{B} (\A^\bullet_\sharp)$ is written
\begin{equation}
\xymatrix{
&& & \\
C^2 (\A) \ar[r]^B \ar[u] & C^1 (\A) \ar[r]^B \ar[u] & C^0 (\A) \ar[u] \\
C^1 (\A) \ar[r]^B \ar[u]^b & C^0 (\A) \ar[u]^b & \\
C^0 (\A) \ar[u]^b & &
}
\end{equation}
and the bicomplex $\mathcal{B}^\textrm{per} (\A) = \mathcal{B}^\textrm{per} (\A^\bullet_\sharp)$ is written
\begin{equation}
\xymatrix{
& & & \vdots &\vdots & \iddots \\
& & \ddots \ar[r] & C^2 (\A) \ar[r]^B \ar[u] & C^1 (\A) \ar[r]^B \ar[u] & C^0 (\A) \ar[u] \\
& \ddots \ar[r] & C^2 (\A)  \ar[r]^B \ar[u] & C^1 (\A) \ar[r]^B \ar[u]^b & C^0 (\A) \ar[u]^b & \\
\cdots \ar[r] & C^2 (\A) \ar[r]^B \ar[u] & C^1 (\A) \ar[r]^B \ar[u]^b & C^0 (\A) \ar[u]^b & &\\
\cdots \ar[r] & C^1 (\A) \ar[r]^B \ar[u]^b & C^0 (\A) \ar[u] &&& \\
& \vdots \ar[u] &  &&&
}
\end{equation}
We denote $ \mathcal{B} (\A^\bullet_\sharp)$ as $\mathcal{B} (\A)$ as we will be dealing with cyclic cohomology.

In this case we can also write $b^n$ and $B^n$ explicitly as
\begin{align}
b^n f (a_0 , \ldots , a_n) &= \sum_{i=0}^{n-1} (-1)^i f(a_0 , \ldots , a_i a_{i+1} , \ldots , a_n) + (-1)^n f(a_n a_0 , \ldots , a_{n-1}) \\
B^n f (a_0 , \ldots , a_n) &= \sum_{i=0}^n (-1)^{ni} f(1 , a_i , \ldots , a_n , a_0 , \ldots , a_{i-1}) \nonumber \\
&- (-1)^{ni} f(a_i , 1 , a_{i+1} , \ldots , a_n , a_0 , \ldots , a_{i-1}) 
\end{align}
for any $f \in C^n(\A)$.

There is another related complex that we will make use of denoted by $\overline{\mathcal{B}} (\A)$. Let 
\begin{equation}
\overline{C}^{pq} (\A) = \Hom (\A \otimes \overline{\A}^{\otimes (q-p)}, k)
\end{equation}
where $\overline{\A} = \A / (k\cdot 1)$.  Also denote $\overline{C}^n (\A) = \Hom(\A \otimes \overline{\A}^n, k)$.  Then the operators $b$ \eqref{eq:Hochschild_coboundary} and $B$ \eqref{eq:Connes_coboundary} are still well defined on this bicomplex.  The operator $B : \overline{C}^n (\A) \rightarrow \overline{C}^{n+1} (\A)$ is given explicitly by 
\begin{equation}
B f (a_0 , \ldots , a_n) = \sum_{i=0}^n (-1)^{ni} f(1 , a_i , \ldots , a_n , a_0 , \ldots , a_{i-1})
\end{equation}
The complex $\overline{\mathcal{B}} (\A) = (\overline{C}^{pq} (\A) , b ,B)$ is quasi-isomorphic to $\mathcal{B} (\A) $. We will denote 
\begin{equation}
\bar{s}^i_n f (a_0 , \ldots , a_n) = (-1)^{ni} f(1 , a_i , \ldots, a_n , a_0 , \ldots , a_{i-1})
\end{equation}
so that $B = \sum_{i=0}^n \bar{s}^i_n$.

For an algebra $\A$ that is not necessarily unital we can use the same sort of construction to compute the cyclic cohomology of $\A$. Let $\widetilde{\A} = \A \oplus k$ be the unitization of $\A$. Then
\begin{equation}
\overline{C}^{pq} (\widetilde{\A}) = \Hom (\widetilde{\A} \otimes \A^{\otimes (q-p)}, k) = \Hom (\A^{\otimes (q-p+1)} \oplus \A^{\otimes (q-p)},k).
\end{equation}
This bicomplex is quasi-isomorphic to $\mathcal{B} (\A)$.

Furthermore, we will make use of the following description of $\mathcal{B} (\A)$ due to \cite{getzler_jones} (see also \cite{loday98}).  Let $u$ be a formal variable of degree $+2$. Then the complex $(C^\bullet (\A) \otimes k[u], b+uB)$, where $k[u]$ denotes polynomials in $u$ over $k$, is quasi-isomorphic to $\mathcal{B} (\A)$ and the complex $(\overline{C}^\bullet (\A) \otimes k[u], b+uB)$ is quasi-isomorphic to $\overline{\mathcal{B}} (\A)$ and all may be referred to as the \textit{cyclic complex} of $\A$. The complex $(C^\bullet (\A) \otimes k[u^{-1}, u], b+ uB)$ is quasi-isomorphic to $\mathcal{B}^\textrm{per} (\A)$ and the complex $(\overline{C}^\bullet (\A) \otimes k[u^{-1}, u], b+ uB)$ is quasi-isomorphic to $\overline{\mathcal{B}}^\textrm{per} (\A)$ and all may be referred to as the \textit{periodic cyclic complex} of $\A$.

\section{Gerbes}


In \cite{hitchin01} N. Hitchin describes a gerbe on a manifold.  We will first discuss an approach similar to Hitchin's. Then we will define a gerbe on a smooth \'etale groupoid.  The latter notion agrees with the former in the case of the \v{C}ech groupoid.

\subsection{Gerbes on Manifolds}
In describing a gerbe on a manifold we will follow \cite{hitchin01} and especially \cite{benameur_gorokhovsky}. See also \cite{brylinski93} and \cite{ruffino09}. 

Given a smooth manifold $M$, a \textit{descent datum on $M$} is a collection $\{U_\alpha , \La_{\alpha \beta} , \mu_{\alpha \beta \gamma} \}$ denoted by $\La$ where $\U = \{U_\alpha\}_{\alpha \in I}$ is an open cover of $M$, $\{\La_{\alpha \beta}\}$ is a collection of line bundles over $U_{\alpha \beta}$, and $\mu_{\alpha \beta \gamma} : \La_{\alpha \beta} \otimes \La_{\beta \gamma} \xrightarrow{\sim}\La_{\alpha \gamma}$ are isomorphisms of line bundles over each triple intersection $U_{\alpha \beta \gamma}$ that satisfy the associativity condition that the diagram
\begin{equation}
\label{eq:gerbe_manifold_associativity}
\bfig
\Square|aarb|[\La_{\alpha \beta} \otimes \La_{\beta \gamma} \otimes \La_{\gamma \delta} ` \La_{\alpha \gamma} \otimes \La_{\gamma \delta} ` \La_{\alpha \beta} \otimes \La_{\beta \delta} ` \La_{\alpha \delta}; \mu_{\alpha \beta \gamma} \otimes \id ` \id \otimes \mu_{\beta \gamma \delta} ` \mu_{\alpha \gamma \delta} ` \mu_{\alpha \beta \delta}]
\efig
\end{equation}
commutes over $U_{\alpha \beta \gamma \delta}$.

An \textit{isomorphism} of descent datum $\psi: \{U_\alpha , \La_{\alpha \beta} , \mu_{\alpha \beta \gamma} \} \rightarrow \{U_\alpha , \La'_{\alpha \beta} , \mu'_{\alpha \beta \gamma} \}$ on the same open cover $\mathcal{U}$ is defined by a collection of line bundles $K_\alpha$ on $U_\alpha$ and a collection of line bundle isomorphisms $\psi_{\alpha \beta} : {K_\alpha}^{-1} \otimes \La_{\alpha \beta} \otimes K_\beta \rightarrow \La'_{\alpha \beta}$ on $U_{\alpha \beta}$ such that the diagram
\begin{equation}
\label{eq:gerbe_manifold_isomorphism}
\bfig
\Square|aarb|[{K_\alpha}^{-1} \otimes \La_{\alpha \beta} \otimes K_\beta \otimes {K_\beta}^{-1} \otimes \La_{\beta \gamma} \otimes K_\gamma ` {K_{\alpha}}^{-1} \otimes \La_{\alpha \gamma} \otimes K_\gamma ` \La_{\alpha \beta}' \otimes \La_{\beta \gamma}' ` \La_{\alpha \gamma}'; \id \otimes \mu_{\alpha \beta \gamma} \otimes \id ` \psi_{\alpha \beta} \otimes \psi_{\beta \gamma} ` \psi_{\alpha \gamma} ` \mu_{\alpha \beta \gamma}']
\efig
\end{equation}
commutes on $U_{\alpha \beta \gamma}$.

In order to address the dependence of a descent datum on the open cover $\mathcal{U}$, recall that the cover $\mathcal{V} = \{V_\beta\}_{\beta \in J}$ is a \textit{refinement} of the cover $\mathcal{U} = \{U_\alpha\}_{\alpha \in I}$ if there is a map $\phi: J \rightarrow I$ such that $V_\beta \subset U_{\phi (\beta)}$.  Such a refinement $\phi$ defines a gerbe on $\mathcal{V}$
\begin{equation}
\{V_\alpha , \La_{\phi(\alpha) \phi(\beta)}|_{V_{\alpha \beta}} , \mu_{\phi(\alpha) \phi(\beta) \phi(\gamma)}|_{V_{\alpha \beta \gamma}} \}_{\alpha, \beta, \gamma \in J}.
\end{equation}
Also recall by standard theory (e.g. \cite{bott_tu}) that any two open covers of $M$, $\mathcal{U}$ and $\mathcal{U}'$ have a common refinement $\mathcal{V}$ and that any open cover of $M$ has a refinement which is a good cover.  Note this means that the collection of covers on $M$ is a directed set and good covers are cofinal in the set of all covers.  Therefore we may define an isomorphism between descent datum $\{U_\alpha , \La_{\alpha \beta} , \mu_{\alpha \beta \gamma} \}$ and $\{U'_\alpha , \La'_{\alpha \beta} , \mu'_{\alpha \beta \gamma} \}$ on different open covers $\mathcal{U} = \{U_\alpha\}_{\alpha \in I}$ and $\mathcal{U}' = \{U'_\alpha\}_{\alpha \in I'}$ as an isomorphism
\begin{multline}
\{V_\alpha , \La_{\phi(\alpha) \phi(\beta)}|_{V_{\alpha \beta}} , \mu_{\phi(\alpha) \phi(\beta) \phi(\gamma)}|_{V_{\alpha \beta \gamma}} \}_{\alpha, \beta, \gamma \in J} \\ \to \{V_\alpha , \La_{\phi'(\alpha) \phi'(\beta)}|_{V_{\alpha \beta}} , \mu_{\phi'(\alpha) \phi'(\beta) \phi'(\gamma)}|_{V_{\alpha \beta \gamma}} \}_{\alpha, \beta, \gamma \in J}
\end{multline}
where $\mathcal{V} = \{V_\alpha\}_{\alpha \in J}$ is a common refinement with $\phi : J \rightarrow I$ and $\phi' : J \rightarrow I'$.  A \textit{Dixmier-Douady gerbe} on $M$ is an equivalence class of descent data on $M$ under isomorphism and will be denoted by $[\La]$ where $\La$ is some descent datum in the equivalence class.

If $\U$ is a good cover, then there is a trivialization of every line bundle $\La_{\alpha \beta}$ over $U_{\alpha \beta}$.  A trivialization of $\La_{\alpha \beta}$ may be seen as a non-vanishing global section $\lambda_{\alpha \beta} : U_{\alpha \beta} \rightarrow \La_{\alpha \beta}$.  Upon choosing a global section $\lambda_{\alpha \beta}$ of each $\La_{\alpha \beta}$, we may view the isomorphisms $\{ \mu_{\alpha \beta \gamma} \}$ as a \v{C}ech $2$-cochain with values in the sheaf $\underline{\C}_M^\times$ of multiplicative complex numbers, $\check{C}^2 ( \mathcal{U} , \underline{\C}_M^\times)$, by defining
\begin{equation}
\mu_{\alpha \beta \gamma} (\lambda_{\alpha \beta} \otimes \lambda_{\beta \gamma}) = \check{\mu}_{\alpha \beta \gamma} \lambda_{\alpha \gamma} \textrm{ on } U_{\alpha \beta \gamma}
\end{equation}
where $\check{\mu}_{\alpha \beta \gamma} : U_{\alpha \beta \gamma} \rightarrow \C^\times$ is a smooth invertible function. This makes sense as $\mu_{\alpha \beta \gamma} (\lambda_{\alpha \beta} \otimes \lambda_{\beta \gamma})$ is a non-vanishing section of $\La_{\alpha \gamma} |_{U_{\alpha \beta \gamma}}$ and we may obtain any non-vanishing section of $\La_{\alpha \gamma}|_{U_{\alpha \beta \gamma}}$ through multiplication of $\lambda_{\alpha \gamma}$ by a non-zero complex number at each point of $U_{\alpha \beta \gamma}$.  Abusing notation, we will write $\mu_{\alpha \beta \gamma}$ for $\check{\mu}_{\alpha \beta \gamma}$.  By the commutativity of \eqref{eq:gerbe_manifold_associativity}, $\{ \mu_{\alpha \beta \gamma} \}$ defines a \v{C}ech 2-cocycle $\underline{\mu} \in \check{Z}^2 (\mathcal{U} , \underline{\C}_M^\times)$, hence a class $[ \mu ] \in \check{H}^2 (\mathcal{U} , \underline{\C}_M^\times)$. Since the \v{C}ech cohomology is the same for all good coverings of $M$ and the \v{C}ech cohomology of $M$ is defined to be the direct limit $\check{H}^* (M , \underline{\C}_M^\times) = \lim_{\mathcal{U}} \check{H}^* (\mathcal{U} , \underline{\C}_M^\times)$, such a gerbe in fact defines a class $[\underline{\mu}] \in \check{H}^2 (M , \underline{\C}_M^\times )$ called the \textit{Dixmier-Douady class} of $\mathcal{L}$.  We may denote this class by $c_1( \mathcal{L})$.

\begin{prop}
Two descent datum $\mathcal{L} = \{U_\alpha , \La_{\alpha \beta} , \mu_{\alpha \beta \gamma} \}$ and $\mathcal{L}' = \{U_\alpha , \La'_{\alpha \beta} , \mu'_{\alpha \beta \gamma} \}$ are isomorphic if and only if $[\mathcal{L}] = [\mathcal{L}']$. In other words a Dixmier-Douady gerbe has a well defined Dixmier-Douady class and a Dixmier-Douady class gives rise to a well defined Dixmier-Douady gerbe.
\begin{proof}
It is sufficient to assume that $\{U_\alpha\}$ is a good cover. Let $\psi : \mathcal{L} \rightarrow \mathcal{L}'$ be an isomorphism given by $\psi_{\alpha \beta} : {K_\alpha}^{-1} \otimes \La_{\alpha \beta} \otimes K_\beta \rightarrow \La'_{\alpha \beta}$ on $U_{\alpha \beta}$.  As every $U_\alpha$ is contractible, we can pick non-vanishing global sections $\kappa_{\alpha}$ for each $K_{\alpha}$ and as before we can pick non-vanishing global sections $\lambda_{\alpha \beta}$ for each $\La_{\alpha \beta}$. Furthermore, over each $U_{\alpha \beta}$, we may write the isomorphism $\psi_{\alpha \beta} ({\kappa_{\alpha}}^{-1} \otimes \lambda_{\alpha \beta} \otimes \kappa_{\beta}) = s_{\alpha \beta} \lambda_{\alpha \beta}'$ for some $s_{\alpha \beta} : U_{\alpha \beta} \rightarrow \underline{\C}_M^\times$.  By the isomorphism condition \eqref{eq:gerbe_manifold_isomorphism}, $\psi_{\alpha \gamma} \circ (\id \otimes \mu_{\alpha \beta \gamma} \otimes \id ) = \mu_{\alpha \beta \gamma}' \circ (\psi_{\alpha \beta} \otimes \psi_{\beta \gamma})$, so that on $U_{\alpha \beta \gamma}$
\begin{align*}
\psi_{\alpha \gamma} \circ (\id \otimes \mu_{\alpha \beta \gamma} \otimes \id ) ({\kappa_{\alpha}}^{-1} \otimes \lambda_{\alpha \beta} \kappa_{\beta} {\kappa_{\beta}}^{-1} \lambda_{\beta \gamma} \kappa_{\gamma}) &= \psi_{\alpha \gamma} ({\kappa_{\alpha}}^{-1} (\mu_{\alpha \beta \gamma} \lambda_{\alpha \gamma}) \kappa_{\gamma} )\\
&= s_{\alpha \gamma} \mu_{\alpha \beta \gamma} \lambda_{\alpha \gamma}' \\
\intertext{and}
\mu_{\alpha \beta \gamma}' \circ \psi_{\alpha \beta} \otimes \psi_{\beta \gamma} ({\kappa_{\alpha}}^{-1} \otimes \lambda_{\alpha \beta} \kappa_{\beta} {\kappa_{\beta}}^{-1} \lambda_{\beta \gamma} \kappa_{\gamma}) &= \mu_{\alpha \beta \gamma}' (s_{\alpha \beta} \lambda_{\alpha \beta} ' s_{\beta \gamma} \lambda_{\beta \gamma} ' ) \\
&= s_{\alpha \beta} s_{\beta \gamma} \mu_{\alpha \beta \gamma} ' \lambda_{\alpha \gamma} '
\end{align*}
are equal.  Therefore
\begin{equation}
\mu_{\alpha \beta \gamma} {\mu_{\alpha \beta \gamma}'}^{-1} = s_{\alpha \beta} s_{\beta \gamma} {s_{\alpha \gamma} }^{-1}
\end{equation}
showing that the difference between $\underline{\mu}$ and $\underline{\mu'}$ is a \v{C}ech coboundary.
\end{proof}
\end{prop}

The short exact sequence
\begin{equation}
0 \rightarrow 2 \pi i \Z \rightarrow \underline{\C}_M \xrightarrow{e^{2 \pi i \cdot}} \underline{\C}^\times_M \rightarrow 0
\end{equation}
of sheaves defines a long exact sequence on \v{C}ech cohomology. The connecting homomorphism of the long exact sequence provides an isomorphism
\begin{equation}
\check{H}^2 (M , \underline{\C}^\times_M) \cong \check{H}^3 (M , 2 \pi i \Z) \cong  \check{H}^3 (M , \Z)
\end{equation}
since $\underline{\C}_M$ is a soft sheaf which implies $\check{H}^i (M , \underline{\C}_M) = 0$ for $i>0$.  Explicitly, we can write down a \v{C}ech 3-cocycle $\underline{\nu}$ with values in $2\pi i \Z$ that represents the Dixmier-Douady class of a descent datum $\La = \{U_\alpha , \La_{\alpha \beta} , \mu_{\alpha \beta \gamma} \}$.  Given a branch $\Log (\mu_{\alpha \beta \gamma})$ of the logarithm of $\mu_{\alpha \beta \gamma}$ over $U_{\alpha \beta \gamma}$, set
\begin{equation}
\nu_{\alpha \beta \gamma \delta} = \Log (\mu_{\beta \gamma \delta}) - \Log (\mu_{\alpha \gamma \delta}) + \Log (\mu_{\alpha \beta \delta}) - \Log (\mu_{\alpha \beta \gamma}) \textrm{ on } U_{\alpha \beta \gamma \delta}.
\end{equation}
Then $\underline{\nu}$ defines a $2 \pi i \Z$-valued \v{C}ech 3-cocycle $\underline{\nu} \in \check{Z}^3 (\mathcal{U} , 2\pi i \Z)$ which represents $c_1 (\La)$.

Given a smooth map $f: M' \rightarrow M$ between smooth manifolds, the pullback of a descent datum $\La$ on $M$ is a descent datum $f^* \La$ on $M'$. The pullback respects isomorphisms so that the pullback of a Dixmier-Douady gerbe is a Dixmier-Douady gerbe.

A \textit{unitary descent datum} is a descent datum $\La = \{U_\alpha , \La_{\alpha \beta} , \mu_{\alpha \beta \gamma} \}$ with a choice of metric on each $\La_{\alpha \beta}$ such that $\mu_{\alpha \beta \gamma}$ is additionally an isometry.  A unitary equivalence of unitary descent data $\La = \{U_\alpha , \La_{\alpha \beta} , \mu_{\alpha \beta \gamma} \}$ and $\La' = \{U_\alpha , \La_{\alpha \beta}' , \mu_{\alpha \beta \gamma}' \}$ sharing the same open cover $\{U_\alpha\}$ is an isomorphism of descent data $\psi_{\alpha \beta} : {K_\alpha}^{-1} \otimes \La_{\alpha \beta} \otimes K_\beta \rightarrow \La'_{\alpha \beta}$ where the line bundles $\{ K_\alpha \}$ are Hermitian and each $\psi_{\alpha \beta}$ is an isometry. Again, a unitary Dixmier-Douady gerbe is an equivalence class of unitary gerbes.

\begin{prop}
Given a descent datum $\La = \{U_\alpha , \La_{\alpha \beta} , \mu_{\alpha \beta \gamma} \}$, there exists a collection of connections on $\La$ denoted $\nabla = \{ \nabla_{\alpha \beta} \}$ satisfying the condition 
\begin{equation}
\label{eq:gerbe_manifold_connection}
\mu_{\alpha \beta \gamma}^* \nabla_{\alpha \gamma} = \nabla_{\alpha \beta} \otimes 1 + 1 \otimes \nabla_{\beta \gamma} \textrm{ on } U_{\alpha \beta \gamma}.
\end{equation}

\begin{proof}
Following \cite{benameur_gorokhovsky}, for any $\alpha, \beta$ such that $U_{\alpha \beta} \neq \emptyset$, fix a connection $\nabla_{\alpha \beta} '$ on $\La_{\alpha \beta}$.  Then for $U_{\alpha \beta \gamma} \neq \emptyset$ define
\begin{equation}
A_{\alpha \beta \gamma} = \mu_{\alpha \beta \gamma}^* \nabla_{\alpha \gamma} ' - \nabla_{\alpha \beta} \otimes 1 - 1 \otimes \nabla_{\beta \gamma}.
\end{equation}
Using the identification
\begin{equation}
\End (\La_{\alpha \beta} \otimes \La_{\beta \gamma}) \cong U_{\alpha \beta \gamma} \times \C
\end{equation}
we see that $A_{\alpha \beta \gamma}$ is identified with a differential 1-form on the open set $U_{\alpha \beta \gamma}$.  We also have for $U_{\alpha \beta \gamma \delta} \neq \emptyset$
\begin{equation}
A_{\beta \gamma \delta} + A_{\alpha \beta \delta} = A_{\alpha \gamma \delta} + A_{\alpha \beta \gamma}.
\end{equation}
Therefore, there exists $A' = (A_{\alpha \beta}')$ such that
\begin{equation}
A_{\alpha \beta \gamma} = A_{\beta \gamma}' - A_{\alpha \gamma}' + A_{\alpha \beta}'.
\end{equation}
Hence the collection of connections $(\nabla_{\alpha \beta} = \nabla_{\alpha \beta} ' + A_{\alpha \beta}' )$ is a connection on $\La$.
\end{proof}
\end{prop}

\begin{prop}
Given a connection $\nabla$ on $\La$ as above, let $\theta_{\alpha \beta} = \nabla_{\alpha \beta}^2$ be the curvatures of the connections $\nabla_{\alpha \beta}$.  Then there exists a collection of 2-forms $\theta_\alpha \in \Omega^2 (U_\alpha)$ that satisfy
\begin{equation}
\theta_{\alpha \beta} = \theta_\beta - \theta_\alpha
\end{equation}
\begin{proof}
Following \cite{benameur_gorokhovsky}, we have a collection $\theta_{\alpha \beta}$ of 2-forms that satisfy
\begin{equation}
\theta_{\alpha \gamma} = \theta_{\alpha \beta} - \theta_{\beta \gamma}
\end{equation}
which means that the collection of $\theta_{\alpha \beta}$ is a \v{C}ech cocycle.  Since the \v{C}ech complex of differential forms is acyclic, there exists some collection of forms $\theta_\alpha$ on each $U_\alpha$ that satisfy
\begin{equation}
\theta_{\alpha \beta} = \theta_\beta - \theta_\alpha
\end{equation}
\end{proof}
\end{prop}

We call $(\nabla_{\alpha \beta} , \theta_\alpha)$ a connection on the descent datum $\La$.

\begin{prop}
Let $\La = \{U_\alpha , \La_{\alpha \beta} , \mu_{\alpha \beta \gamma} \}$ be a descent datum with a connection $( \{ \nabla_{\alpha \beta} \} , \{\theta_\alpha\})$.  Then there exists a well-defined closed form $\Omega \in \Omega^3 (M)$ called the curvature $3$-form of the connection that is locally $\Omega |_{U_\alpha} = \frac{d\theta_\alpha}{2 \pi i}$. Furthermore, if $\Omega'$ is the curvature 3-form of another connection then $\Omega ' = \Omega + d\eta$ for $\eta \in \Omega^2 (M) / d\Omega^1(M)$.
\begin{proof}
As in \cite{benameur_gorokhovsky}, since $d \theta_{\alpha \beta} = 0$ for any $\alpha , \beta$,
\begin{equation*}
d \theta_\alpha |_{U_{\alpha \beta}} = d \omega_\beta |_{U_{\alpha \beta}}
\end{equation*}
which shows the existence of the closed form $\Omega$.  Restriction to a refinement clearly does not change $\Omega$ nor does replacing the connection by an equivalent one.
\end{proof}
\end{prop}

The curvature $3$-form of a connection $(\nabla_{\alpha \beta} , \theta_\alpha)$ on a descent datum $\La$ is a representative in de Rham cohomology of $c_1 (\La)$.

\subsection{Gerbes on Groupoids}

\begin{defn}
Given an \'etale groupoid $\G$, a \textit{gerbe} $(L,\mu)$ over $\G$ is a line bundle over $\Gl[1]$, $L \rightarrow \G$, together with an isomorphism 
\begin{equation}
\label{eq:gerbe_multiplication}
\mu : (\varpi^2_{01})^* L \otimes (\varpi^2_{12})^* L \xrightarrow{\sim} (\varpi^2_{02})^* L
\end{equation}
of line bundles on $\Gl[2]$.  This isomorphism may be written as $\mu_{(g_1 , g_2)} : L|_{g_1} \otimes L|_{g_2} \xrightarrow{\sim} L|_{g_1 g_2}$ on the fiber over $(g_1 , g_2) \in \Gl[2]$.  Here $L_g$ denotes the fiber over $g$.  In addition, $\mu$ must satisfy the associativity condition that the diagram
\begin{equation}
\bfig
\Square|aarb|[(\varpi^3_{01})^* L \otimes (\varpi^3_{12})^* L \otimes (\varpi^3_{23})^* L ` (\varpi^3_{012})^* L \otimes (\varpi^3_{23})^* L ` (\varpi^3_{01})^* L \otimes (\varpi^3_{123})^* L ` (\varpi^3_{0123})^* L; (\varpi^3_{012})^* (\mu) \otimes \id ` \id \otimes (\varpi^3_{123})^* (\mu) ` (\varpi^3_{023})^* (\mu) ` (\varpi^3_{013})^* (\mu)]
\efig
\end{equation}
commutes. On fibers this condition means that the diagram
\begin{equation}
\bfig
\Square|aarb|[L|_{g_1} \otimes L|_{g_2} \otimes L|_{g_3} ` L|_{g_1 g_2} \otimes L|_{g_3} ` L|_{g_1} \otimes L|_{g_2 g_3} ` L|_{g_1 g_2 g_3} ; \mu_{(g_1 , g_2)} \otimes \id ` \id \otimes \mu_{(g_2, g_3)} ` \mu_{(g_1 g_2 , g_3)} ` \mu_{(g_1 , g_2 g_3)}]
\efig
\end{equation}
commutes. 
\end{defn}


Let $\nabla$ be a connection on $L \rightarrow \G$, i.e.\ a linear operator $\nabla : C^\infty_c (\G , L) \rightarrow \Omega^1_c (\G , L)$ where $C^\infty_c (\G , L)$ denotes compactly supported sections of $L \rightarrow \G$ and $\Omega^1_c (\G, L)$ denotes compactly supported $L$-valued $1$-forms on $\G$. Then for a local section $s_U : U \rightarrow L$, $U \subset \Gl[1]$, we may write the operator $\nabla$ as $\nabla s_U = (d+\omega) s_U$ for some $\omega \in \Omega^1 (\G, L)$.  We call $\omega$ the \textit{connection form} of $\nabla$. A connection $\nabla$ on $L \rightarrow \G$ determines a connection on $(\varpi^2_{01})^* L \otimes (\varpi^2_{12})^* L$ given by $(\varpi^2_{01})^* \nabla \otimes 1 + 1 \otimes (\varpi^2_{12})^* \nabla$ as well as a connection $(\mu)^* (\varpi^2_{02})^* \nabla$ on $(\mu)^* (\varpi^2_{02})^* L$.  Their discrepancy is given by a $1$-form
\begin{equation}
\label{eq:gerbe_connection1}
(\varpi^2_{01})^* \nabla \otimes 1 + 1 \otimes (\varpi^2_{12})^* \nabla - (\mu)^* (\varpi^2_{02})^* \nabla = \alpha \in \Omega^1 (\Gl[2])
\end{equation}
Suppose additionally that there is a two form $\theta \in \Omega^2 (\Gl[1])$ satisfying
\begin{equation}
\label{eq:gerbe_connection2}
(\varpi^2_{01})^* \theta - (\varpi^2_{02})^* \theta +  (\varpi^2_{12})^* \theta  = d \alpha,
\end{equation}
in other words, satisfying $\delta^* \theta = d \alpha$. Then $(-\alpha, \theta) \in \Omega^1(\Gl[2]) \oplus \Omega^2 (\Gl[1])$ determines a class in $H^3_{dR} (\G)$, which we call the \textit{Dixmier-Douady class} of $\nabla$ on $L\rightarrow \G$.  Although the Dixmier-Douady class depends on $\nabla$ throughout, it is expected this dependence can be removed.  We call $(\nabla, \theta)$ a connection on $L$ and $\alpha$ the \textit{discrepancy} of $\nabla$. 


This notion of gerbe and Dixmier-Douady class agrees with the notion of a gerbe on a manifold when $\G$ is the \v{C}ech groupoid.


%

\chapter{Twisted cohomology}
\label{chap:2}

In this chapter we will develop the notion of a twisted simplicial complex for the simplicial manifold $(M \rtimes \Gamma)_\bullet$ determined by the translation groupoid arising from the action of a discrete group on a manifold.


\section{Gerbes on the (discrete) translation groupoid}

First we will make a few remarks about the discrete translation groupoid.  Let $\Gamma$ be a discrete group and let $M$ be a manifold carrying a (right) action of $\Gamma$, $M \times \Gamma \rightarrow M$, denoted $(x,g) \mapsto xg$.  The nerve of $M \rtimes \Gamma$ is a simplicial manifold $(M \rtimes \Gamma)_\bullet$ which we will often denote merely by $M_\bullet$.  On objects, $M_\bullet$ is given by $M_{(k)} = M \times \Gamma^k$ and on morphisms $M_\bullet$ is given by
\begin{align}
\label{eq:simp_mfld_face_maps}
\delta^k_i (x, g_1 , \ldots , g_k) &= \begin{cases} (xg_1 , g_2 , \ldots , g_k) & \textrm{if } i=0 \\ (x, g_1 , \ldots , g_{i} g_{i+1} , \ldots , g_k) & \textrm{for } 1 \leq i \leq k-1 \\ (x, g_1 , \ldots , g_{k-1}) & \textrm{if } i=k \end{cases} \\
\label{eq:simp_mfld_deg_maps}
\sigma^k_j (x, g_1 , \ldots , g_k) &= (x, g_1 , \ldots , g_j , 1_\Gamma , g_{j+1} , \ldots , g_k)
\end{align}
in the notation of \eqref{eq:nerve}.  In particular, $\delta^1_0 (x, g_1) = xg_1$ and $\delta^1_1 (x , g_1) = x$.

Let $L \rightarrow M \rtimes \Gamma$ be a gerbe with an isomorphism $\mu$ as in \eqref{eq:gerbe_multiplication}.  As $(M \rtimes \Gamma)_{(1)} = M \times \Gamma$ and $\Gamma$ is a discrete group, we may view $L$ as a collection of line bundles $L_g$ on $M_g = M \times \{ g \}$ for each $g \in \Gamma$. The isomorphism $\mu$ may be restricted in the same way to a collection of isomorphisms
\begin{equation}
\mu_{g_1 , g_2} : L_{g_1} \otimes (L_{g_2})^{g_1} \xrightarrow{\sim} L_{g_1 g_2}
\end{equation}
where $g_1 , g_2 \in \Gamma$ and $(L_{g_2})^{g_1}$ denotes the line bundle $L_{g_2}$ shifted by the action of $g_1$.  In other words $(L_{g_2})^{g_1}$ has sections $s^{g_1}$ where $s^{g_1} (x) = s(xg_1)$ denotes the (left) action of $\Gamma$ for $s: M_{g_2} \rightarrow L|_{g_2}$ a section of $L|_{g_2}$.

Now let $\nabla$ be a connection on $L$ and let $\nabla_g = \nabla|_{M_g}$ be a connection on $L_g$ for each $g \in \Gamma$.   In this case the condition \eqref{eq:gerbe_connection1} may be written
\begin{equation}
\label{eq:gerbe_connection1_group}
\nabla_g \otimes 1 + 1 \otimes (\nabla_h)^g - \mu_{g_1 , g_2}^* (\nabla_{gh}) = \alpha (g,h)
\end{equation}
for all $g,h \in \Gamma$ and for some $\alpha (g,h) \in \Omega^1(M)$. Hence $\alpha \in \Omega^1 ( (M \rtimes \Gamma)_{(2)}) = \Omega^1 (M \times \Gamma \times \Gamma)$.

Define $\theta_g \in \Omega^2 (M)$ by $[\theta_g , s] = (\nabla_g)^2 s$ for a section $s : M_g \rightarrow L_g$. Then we have a collection $\theta = (\theta_g)_{g \in \Gamma} \in \Omega^2 (M \times \Gamma) = \Omega^2 ( (M \rtimes \Gamma)_{(1)} )$ and the following relation.
\begin{prop}
\label{prop:gerbe_connection2_group}
The 2-forms $\theta_g$ satisfy condition \eqref{eq:gerbe_connection2}, which may be written
\begin{equation}
\label{eq:gerbe_connection2_group}
\theta_g + \theta_h^g - \theta_{gh} = d\alpha(g,h)
\end{equation}
\begin{proof}
This follows from the condition \ref{eq:gerbe_connection1_group}.  The curvature of the sum of connections $\nabla_g \otimes 1 + 1 \otimes (\nabla_h)^g$ is $\theta_g + \theta_h^g$ and the curvature of the connection plus scalar form $\mu_{g_1 , g_2}^* (\nabla_{gh}) + \alpha (g,h)$ is $ \theta_{gh} + d\alpha(g,h)$.
\end{proof}
\end{prop}
In view of the previous proposition, a connection $\nabla$ on $L \rightarrow M \rtimes \Gamma$ induces a connection $(\nabla , \theta)$ with discrepancy $\alpha \in \Omega^2 ( (M \rtimes \Gamma)_{(1)})$ on the gerbe $(L,\mu)$ over $M \rtimes \Gamma$.  We will nonetheless often denote such a connection on $(L,\mu)$ by $\nabla$.  As defined above, $(-\alpha , \theta)$ represents the Dixmier-Douady class of $\nabla$ on $(L , \mu)$ over $M \rtimes \Gamma$.


\section{Twisted bundles}

\begin{defn}
Given a gerbe $(L,\mu)$ on $M \rtimes \Gamma$, an $L$-twisted vector bundle $\E$ is a vector bundle on $M$ together with a collection of vector bundle isomorphisms $\varphi_g : \E \xrightarrow{\sim} \E^g \otimes L_g$ for each $g \in \Gamma$ such that the diagram
\begin{equation}
\label{eq:twisted_bundle}
\bfig
\Square|aarb|[\E^{gh} \otimes (L_h)^g \otimes L_g ` \E^{gh} \otimes L_{gh} ` \E^g \otimes L_g ` \E; \id \otimes \mu_{g,h} `  ( (\varphi_h)^g)^{-1} \otimes \id ` \varphi_{gh}^{-1} ` \varphi_g^{-1}]
\efig
\end{equation}
commutes. Here $\E^g$ denotes the induced action of $g \in \Gamma$ on the vector bundle $\E \rightarrow M$.
\end{defn}

\begin{rmk}
Note that $\varphi_g$ induces an isomorphism $\varphi_g : \End \E \xrightarrow{\sim} \End \E^g$ as well. More precisely, $\varphi_g : \E \xrightarrow{\sim} \E^g \otimes L_g$ induces an isomorphism
\begin{equation*}
\varphi_g : \End \E \xrightarrow{\sim} \End (\E^g \otimes L_g) \cong \End (\E^g ) \otimes \End (L_g).
\end{equation*}
As $\End (L_g) \cong \C$ we see $\End \E \cong \End \E^g$ via $\varphi_g$.
\end{rmk}

Let $\nabla^\E : C^\infty (M , \E) \to \Omega^1 (M , \E)$ be a connection on $\E$ and $(L,\mu)$ a gerbe on $M \rtimes \Gamma$ with connection $(\nabla, \theta)$ and discrepancy $\alpha$.  Then there is a connection on $\E^g \otimes L_g$ given by ${(\nabla^\E)}^g \otimes 1 + 1 \otimes \nabla_g$ for every $g \in \Gamma$.  By ${(\nabla^\E)}^g$ we mean the connection ${(\nabla^\E)}^g s^g = {(\nabla^\E s)}^g$ where $s^g$ is a section of $\E^g$.  This defines a connection on $\E^g$ which we denote $\nabla^{\E^g} = {(\nabla^\E)}^g$.  The discrepancy between $\nabla^\E$ and $\nabla^{\E^g} \otimes 1 + 1 \otimes \nabla_g$ is an $\End \E$-valued 1-form which we denote by
\begin{equation}
\nabla^\E - \varphi^*_g ( \nabla^{\E^g} \otimes 1 + 1 \otimes \nabla_g)  = -A(g) \in \Omega^1 (M , \End \E)
\end{equation}
for each $g \in \Gamma$. We will call $A \in \Omega^1 (M \times \Gamma , \End \E)$ the \textit{discrepancy} of $\nabla^\E$. Notice that we then have the identity
\begin{align}
\label{eq:twisted_bundle_curvature}
(\nabla^\E + A(g))^2 &= \varphi_g^* ((\nabla^{\E^g}) \otimes 1 + 1 \otimes \nabla_g)^2 \notag \\
 (\nabla^\E)^2 + A(g)^2 + [\nabla^\E , A(g)] &= \varphi_g^* ((\nabla^{\E^g})^2 ) + \theta_g \notag \\
 \theta^\E + A(g)^2 + [\nabla^\E , A(g)] &= \varphi_g^* (\theta^\E)^g + \theta_g 
\end{align}
on $\Omega^2 (M , \End \E)$ where we define $\theta^\E = (\nabla^\E)^2$.


\begin{ex}[Direct Sum Bundle]
\label{ex:direct_sum_bundle}

Any gerbe $(L,\mu)$ on $M \rtimes \Gamma$ determines an $L$-twisted vector bundle called the \textit{direct sum bundle} of $(L, \mu)$. The direct sum bundle is the most relevant example of a twisted vector bundle for our purposes. Given a gerbe $(L, \mu)$, define the infinite dimensional vector bundle $\pi : \E \rightarrow M$ where $\E = \bigoplus_{g \in \Gamma} L_g$.  In addition, define isomorphisms $\varphi_g : \E \to \E^g \otimes L_g$ for each $g \in \Gamma$ by
\begin{equation}
\varphi_g : \bigoplus_{g' \in \Gamma} L_{g'} \to \bigoplus_{g' \in \Gamma} \mu_{g , g'}^{-1} (L_{g'}) = \bigoplus_{g' \in \Gamma} (L_{g'})^g \otimes L_g
\end{equation}
where the factor $L_g'$ is mapped to $\mu_{g , g'}^{-1} (L_{g'})$.  With these isomorphisms, $\E$ is an $L$-twisted vector bundle on $M$.  In addition, a connection $\nabla$ on $L \rightarrow M \rtimes \Gamma$ induces the direct sum connection $\nabla^\E = \bigoplus_{g \in \Gamma} \nabla_g$ on $\E \rightarrow M$.  We will again denote the discrepancy of $\nabla^\E$ by $A$.  The following lemma is a useful identity relating the discrepancy of $\nabla$ with the discrepancy of $\nabla^\E$.

\begin{lem}
\label{lem:discrepancy_identity}
Let $\E$ be the direct sum bundle of $(L,\mu)$, a gerbe on $M \rtimes \Gamma$ with connection $(\nabla , \theta)$, and for each $g \in \Gamma$, let $\varphi_g : \E \xrightarrow{\sim} \E^g \otimes L_g$ be an isomorphism.  Then for any $g, h \in \Gamma$ we have the identity
\begin{equation}
A(g) + (\varphi_g)^* (A(h)^g) - A(gh) = - \alpha (g,h) \in \Omega^1 (M , \End \E),
\end{equation}
where $\alpha$ is the discrepancy of $\nabla$ and $A$ is the discrepancy of $\nabla^\E$.
\begin{proof}
This follows from the commutativity of \eqref{eq:twisted_bundle}.  From the right arrow, we have the condition
\begin{equation}
\nabla^\E - \varphi_{gh}^* (\nabla^{\E^{gh}} \otimes 1 + 1 \otimes \nabla_{gh}) = -A(gh)
\end{equation}
and from the top arrow we have the condition
\begin{multline}
(\nabla^{\E^{gh}} \otimes 1 + 1 \otimes \nabla_{gh}) - (1 \otimes \mu_{gh}^{-1})^* (\nabla^{\E^{gh}} \otimes 1 \otimes 1 + 1 \otimes (\nabla_h)^g \otimes 1 + 1 \otimes 1 \otimes \nabla_g) \\= -(1 \otimes \mu_{gh}^{-1})^* (1 \otimes \alpha(g,h)) = -1 \otimes \alpha (g,h).
\end{multline}
From the composition of these two arrows we have the condition
\begin{multline}
\label{eq:discrepancy_identity_1}
\nabla^\E - \varphi_{gh}^* (1 \otimes \mu_{gh}^{-1})^* (\nabla^{\E^{gh}} \otimes 1 \otimes 1 + 1 \otimes (\nabla_h)^g \otimes 1 + 1 \otimes 1 \otimes \nabla_g) \\= -A(gh) + \varphi_{gh}^* (1 \otimes \alpha (g,h)).
\end{multline}
Next, the bottom arrow implies the condition
\begin{equation}
\nabla^\E - \varphi_g^* (\nabla^{\E^g} \otimes 1 + 1 \otimes \nabla_g) = -A(g)
\end{equation}
and the left arrow implies the condition
\begin{multline}
(\nabla^{\E^g} \otimes 1 + 1 \otimes \nabla_g) - ((\varphi_h)^g)^* (\nabla^{\E^{gh}} \otimes 1 \otimes 1 + 1 \otimes (\nabla_h)^g \otimes 1 + 1 \otimes 1 \otimes \nabla_g) \\=  -(A(h))^g
\end{multline}
which compose to give the condition
\begin{multline}
\label{eq:discrepancy_identity_2}
\nabla^\E - \varphi_g^* ((\varphi_h)^g \otimes 1)^* (\nabla^{\E^{gh}} \otimes 1 \otimes 1 + 1 \otimes (\nabla_h)^g + 1 \otimes 1 \nabla_g) \\= -A(g) - (\varphi_g)^* (A(h)^g).
\end{multline}
Since the diagram \eqref{eq:twisted_bundle} commutes, we can equate \eqref{eq:discrepancy_identity_1} and \eqref{eq:discrepancy_identity_2} to produce
\begin{equation}
-A(gh) + \varphi^*_{gh}  (1 \otimes \alpha (g,h)) = -A(g) - (\varphi_g)^* (A(h)^g)
\end{equation}
which is the above identity.
\end{proof}
\end{lem}

\begin{cor}
\label{cor:discrepancy_identity2}
With the hypothesis of the previous lemma, as well as $1 < i < j$, we have
\begin{multline}
\alpha (g_1 , g_2 \cdots g_i) - \alpha (g_1 , g_2 \cdots g_j) + \alpha (g_1 \cdots g_i , g_{i+1} \cdots g_j) \\=  \varphi_{g_1}^* \alpha (g_2 \cdots g_i , g_{i+1} \cdots g_j)^{g_1} 
\end{multline}
where $\alpha$ is the discrepancy of $\nabla$.
\begin{proof}
This follows from the calculation
\begin{align*}
& \alpha (g_1 , g_2 \cdots g_i) - \alpha (g_1 , g_2 \cdots g_j) + \alpha (g_1 \cdots g_i , g_{i+1} \cdots g_j) = \\
&= A(g_1) + \varphi_{g_1}^* (A(g_2 \cdots g_i)^{g_1}) - A(g_1 \cdots g_i) \\
&\quad - A(g_1) - \varphi_{g_1}^* (A(g_2 \cdots g_j)^{g_1}) + A(g_1 \cdots g_j) \\
&\quad + A(g_1 \cdots g_i) + \varphi_{g_1 \cdots g_i}^* (A(g_{i+1} \cdots g_j))^{g_1 \cdots g_i} - A(g_1 \cdots g_j) \\
&= \varphi_{g_1}^* (A(g_2 \cdots g_i)^{g_1}) - \varphi_{g_1}^* (A(g_2 \cdots g_j)^{g_1}) + \varphi_{g_1 \cdots g_i}^* (A(g_{i+1} \cdots g_j))^{g_1 \cdots g_i} \\
&= \varphi_{g_1}^* (A(g_2 \cdots g_i) - A(g_2 \cdots g_j) + \varphi_{g_2 \cdots g_i}^* (A(g_{i+1} \cdots g_j))^{g_2 \cdots g_i} )^{g_1} \\
&= \varphi_{g_1}^* \alpha (g_2 \cdots g_i , g_{i+1} \cdots g_j)^{g_1}
\end{align*}
where we have used the previous lemma.
\end{proof}
\end{cor}

The bundle $\End \E \rightarrow M$ is vector bundle of finite rank endomorphisms when $\E$ is the direct sum bundle.  Note that any section of $C^\infty (M , \End \E)$ may be decomposed as follows.  Since $\E = \bigoplus_{g \in \Gamma} L_g$, $\End \E$ consists of morphisms $ \phi : \bigoplus_{g \in \Gamma} L_g \rightarrow \bigoplus_{g \in \Gamma} L_g$.  Let $E_{g , h} (\phi) \in \Hom (L_g , L_h)$ denote the restriction of $\phi$ to the component $E_{g,h} (\phi): L_g \rightarrow L_h$, so that $\phi = \bigoplus_{g,h \in \Gamma} E_{g,h} (\phi)$.  So given a section $f \in C^\infty (M , \End \E)$ we may decompose $f$ in the same way, with $E_{g,h} (f)$ a section on $M$ with values in morphisms $L_g \rightarrow L_h$. In other words, $E_{g,h} (f) \in C^\infty (M , \Hom (L_g , L_h))$.

Given such a decomposition we can describe the action of $g \in \Gamma$ on sections $C^\infty (M , \End \E)$ as follows.  Let $f \in C^\infty (M , \Hom(L_{g_1} , L_{g_2}) )$, i.e.\ $f=E_{g_1 , g_2} (f)$.  For any $g \in \Gamma$, the section $f^g$ is a section of $\Hom (L_{g_1}^g , L_{g_2}^g)$ and $\id \otimes f^g$ is a section of $\Hom (L_{g} \otimes (L_{g_1})^{g} , L_{g} \otimes (L_{g_2})^{g})$.  We may denote the corresponding section by $f^{g}$ as well.  As there are isomorphisms $\mu_{g, g_1} : L_{g} \otimes (L_{g_1})^{g} \rightarrow L_{g g_1}$ and $\mu_{g , g_2} : L_{g} \otimes (L_{g_2})^{g} \rightarrow L_{g g_2}$, under these isomorphisms $f^{g}$ corresponds to a section in $C^\infty (M , \Hom (L_{g g_1} , L_{g g_2}) )$ which is denoted by $E_{g g_1 , g g_2} (f^g)$.  In other words
\begin{equation*}
E_{g g_1 , g g_2} (f^g) = \mu_{g , g_2} \circ \id \otimes f^g \circ \mu_{g , g_1}^{-1}
\end{equation*}
Hence the (left) action of $g \in \Gamma$ on $E_{g_1 , g_2} (f)$ is $g \cdot E_{g_1 , g_2} (f) = E_{g g_1 , g g_2} (f^g)$.

For a section $a \in C^\infty (M_g , L_g)$, we may consider $a$ as the section $E_{1, g} ( a) : L_1 \rightarrow L_g$ of $C^\infty (M , \End \E)$ since $a : M_g \rightarrow L_g$ defines a section $\tilde{a} : M \times \C \rightarrow L_g$ by $\tilde{a} ( x, \lambda) = a(x)$ for all $\lambda \in \C$ and $L_1 \cong M \times \C$. We have $L_1 \cong M \times \C$ as $\mu_{1 , 1} : L_1 \otimes L_1 \xrightarrow{\sim} L_1$.  
We will not make a distinction in notation between a section $a : M \rightarrow L_g$ and the corresponding section of $ C^\infty (M , \Hom (L_1 , L_g ) )$.

Let $f_1 , f_2 \in C^\infty (M , \End \E)$ be sections of $\End \E$. The sections $E_{g_1 , g_2} (f_1)$ and $E_{g_2 , g_3} (f_2)$ may be composed to produce a section $E_{g_2 , g_3} (f_2) \circ E_{g_1, g_2} (f_1)$. We will denote this operation of composition by
\begin{equation}
E_{g_1, g_2} (f_1) E_{g_2 , g_3} (f_2) = E_{g_1 , g_3} ( f_1 f_2 )  \in C^\infty (M , \Hom (L_{g_1},  L_{g_3}))
\end{equation}
Thus we may obtain a matrix product on $C^\infty (M , \End \E)$.  Let $f_1 f_2 \in C^\infty (M , \End \E)$ be defined by
\begin{equation}
E_{h , k} (f_1 f_2) = \sum_{g \in \Gamma} E_{h , g} (f_1) E_{g , k} (f_2)
\end{equation}
for $h,k \in \Gamma$ so that $f_1 f_2 = \bigoplus_{h,k \in \Gamma} E_{h,k} (f_1 f_2)$.

\begin{lem}

Given sections $a_g \in C^\infty (M , L_g)$ and $a_h \in C^\infty (M , L_h)$
\begin{equation*}
E_{1,g} (a_g) E_{g,h} (a_h^g) = E_{1, gh} (\mu_{g,h} (a_g \otimes a_h^g) )
\end{equation*}
where we consider  $a_g \in C^\infty (M , \Hom (L_1 , L_g))$, $a_h \in C^\infty (M , \Hom (L_1 , L_h))$, as well as $\mu_{g,h} (a_g \otimes a_h^g) \in C^\infty (M , \Hom (L_1 , L_{gh}) )$.
\begin{proof}
The left hand side is defined by
\begin{equation*}
L_1 \xrightarrow{a_g} L_g \xrightarrow{\mu_{g,1}^{-1} }  L_g \otimes L_1^g \xrightarrow{\id \otimes a_h^g} L_g \otimes L_h^g  \xrightarrow{\mu_{g,h}} L_{gh}
\end{equation*}
while the right hand side is defined by
\begin{equation*}
L_1 \xrightarrow{ a_g \otimes a_h^g } L_g \otimes L_h^g \xrightarrow{ \mu_{g,h} } L_{gh}
\end{equation*}
which are corresponding sections of $\Hom (L_1 , L_{gh})$ under the isomorphism $\mu$.
\end{proof}
\end{lem}

For $f \in C^\infty (M , \End \E)$ there is a trace map $\tr : C^\infty (M , \End \E) \rightarrow C^\infty (M)$ defined by
\begin{equation}
\tr(f) = \sum_{g \in \Gamma} E_{g , g} (f)
\end{equation}
which is an element of $C^\infty (M)$ as $E_{g, g} (f) \in \End (L_{g})$ for all $g \in \Gamma$.  In fact, for $\Omega^* (M, \End \E) = \Omega^* (M) \otimes C^\infty (M , \End \E)$, this trace map induces a trace
\begin{equation}
\tr : \Omega^* (M , \End \E) \rightarrow \Omega^* (M)
\end{equation}
as for $\omega \otimes f \in \Omega^*(M) \otimes C^\infty (M, \End \E)$ we may define $\tr (\omega \otimes f) = \omega \tr(f)$.
\begin{lem}
\label{lem:trace_commutator}
For $\omega_1, \omega_2 \in \Omega^* (M , \End \E)$
\begin{equation}
\tr (\omega_1 \omega_2) = (-1)^{ |\omega_1| |\omega_2|} \tr(\omega_2 \omega_1),
\end{equation}
and furthermore,
\begin{equation}
\tr ( [\omega_1 , \omega_2]) = 0,
\end{equation}
where $[\omega_1 , \omega_2]$ is the graded commutator of $\omega_1$ and $\omega_2$.
\begin{proof}
We will first mention that $\tr (f_1 f_2) = \tr (f_2 f_1)$ for any $f_1 , f_2 \in C^\infty (M, \End \E)$. This follows as usual
\begin{align*}
\tr(f_1 f_2) &= \sum_{g \in \Gamma} E_{g,g} (f_1 f_2) \\
&= \sum_{g \in \Gamma}  \left( \sum_{h \in \Gamma} E_{g , h} (f_1) E_{h,g} (f_2) \right) \\
&= \sum_{h \in \Gamma} \left( \sum_{g \in \Gamma}  E_{h,g} (f_2) E_{g,h} (f_1) \right) \\
&= \sum_{h \in \Gamma} E_{h,h} (f_2 f_1) \\
&= \tr(f_2 f_1)
\end{align*}
by changing the summation order.

Let $\omega_1 = \alpha_1 \otimes f_1 $, $\omega_2 = \alpha_2 \otimes f_2$ where $\alpha_1 , \alpha_2 \in \Omega^* (M)$ and $f_1, f_2 \in C^\infty (M, \End \E)$.  Then $| \alpha_1 | = | \omega_1|$ and $| \alpha_2 | = | \omega_2 |$ and
\begin{align*}
\tr (\omega_1 \omega_2) &= \tr ((\alpha_1 \wedge \alpha_2) \otimes (f_1 f_2) ) \\
&= (\alpha_1 \wedge \alpha_2) \tr ( (f_1 f_2) ) \\
&= (-1)^{|\alpha_1| |\alpha_2|} (\alpha_2 \wedge \alpha_1) \tr(  (f_2 f_1) ) \\
&= (-1)^{|\omega_1| |\omega_2|} \tr(\omega_2 \omega_1)
\end{align*}
From this it follows that
\begin{equation*}
\tr ([\omega_1, \omega_2]) = \tr(\omega_1 \omega_2 - (-1)^{|\omega_1| |\omega_2|} \omega_2 \omega_1) =  \tr(\omega_1 \omega_2) - (-1)^{|\omega_1| |\omega_2|} \tr( \omega_2 \omega_1) = 0.
\end{equation*}
\end{proof}
\end{lem}

\end{ex}


\section{Twisted simplicial cohomology}

\subsection{Nonorientable cohomology}


If $M$ is a non-orientable manifold we must additionally discuss twisting by the orientation bundle of $M$.  Given an atlas $\{(U_\alpha , \phi_\alpha)\}$ of $M$ with transition functions $\{\phi_{\alpha \beta}\}$ the \textit{orientation bundle} of $M$ is the real line bundle $\tau$ on $M$ given by transition functions $\sgn J(\phi_{\alpha \beta})$, where $J(\phi_{\alpha \beta})$ is the Jacobian determinant of the matrix of partial derivatives of $\phi_{\alpha \beta}$ and $\sgn$ is the sign function
\begin{equation}
\sgn(x) = \begin{cases} +1 & \textrm{for $x$ positive} \\ 0 & \textrm{for $x=0$} \\ -1 & \textrm{for $x$ negative.} \end{cases}
\end{equation}
Explicitly, $\tau$ may be constructed from $\bigsqcup_\alpha U_\alpha \times \R^1/\sim$, where $\sim$ is the equivalence relation
\begin{equation}
(x,v) \sim (x, \sgn J(\phi_{\alpha \beta} (x)) v),
\end{equation}
for $(x,v) \in U_\alpha \times \R^1$ and $(x , \sgn J(\phi_{\alpha \beta} (x)) v) \in U_\beta \times \R^1$. The coordinate charts $\phi_\alpha$ define local trivializations $\phi_\alpha' : \tau|_{U_\alpha} \xrightarrow{\sim} U_\alpha \times \R^1$. 

The \textit{$\tau$-twisted de Rham complex},  $\Omega_\tau^* (M) =  \Omega^* (M , \tau)$ is the algebra of $\tau$-valued differential forms and is independent of the choice of atlas.  Also, while the exterior derivative $d_\phi$ depends on the choice of atlas, any another choice of atlas produces an isomorphic algebra, so we will denote the exterior derivative by $d$.  The elements of $\Omega_\tau^* (M)$ are called \textit{densities}. The cohomology of the complex $(\Omega_\tau^* (M) , d)$ is called the \textit{$\tau$-twisted de Rham cohomology} $H_\tau^* (M)$. The \textit{$\tau$-twisted de Rham cohomology} with compact support, $H_{\tau,c}^* (M)$ is defined similarly. For more details see \cite{bott_tu} Chapter I.7.

\subsection{Twisted cohomology}

For a smooth manifold $M$ of dimension $m=\dim M$, let $\Theta \in \Omega^3 (M)$ be a closed $3$-form.  Let $u$ be a formal variable of degree $+2$.  Define 
\begin{equation}
\label{eq:twisted_complex_grading}
\Omega_k (M) = \Omega_\tau^{m - k } (M).
\end{equation}
The de Rham differential $d_{dR} : \Omega_k (M) \rightarrow \Omega_{k-1}$ is of degree $-1$ and exterior multiplication by $\Theta$, $\Theta \wedge \cdot : \Omega_k (M) \rightarrow \Omega_{k-3} (M)$ is of degree $-3$.  Let $\Omega_k (M) [u]$ denote polynomials in $u$ over differential forms in $\Omega_\ell (M)$ such that the form degree $\ell$ plus the degree of powers of $u$ sum to $k$.  In other words elements $\omega \in \Omega_k (M)[u]$ are of the form
\begin{equation}
\label{eq:decomposition_rescaled_de_Rham_complex}
\omega =
\begin{cases}
\omega_k + u \omega_{k-2} + \cdots + u^{k/2} \omega_0 & \textrm{if $k$ is even}\\
\omega_k + u \omega_{k-2} + \cdots + u^{(k-1)/2} \omega_1 & \textrm{if $k$ is odd}
\end{cases}
\end{equation}
where $\omega_j \in \Omega_j (M)$, i.e.\ $\omega_j$ is an $(m-j)$-form on $M$ for $0 \leq j \leq k$.
We can rescale the de Rham differential by $u$ to obtain a degree $+1$ differential
\begin{equation}
u d_{dR} : \Omega_k (M) [u] \rightarrow \Omega_{k+1} (M) [u]
\end{equation}
\begin{lem}
Given $\alpha \in \Omega^k (M) [u]$ and $\beta \in \Omega_\ell (M)[u]$
\begin{equation}
(u d_{dR}) (\alpha \wedge \beta) = (u d_{dR} \alpha) \wedge \beta + (-1)^k \alpha \wedge (u d_{dR} \beta)
\end{equation}
\begin{proof}
This follows by using the decomposition of $\beta$ \eqref{eq:decomposition_rescaled_de_Rham_complex} and a similar decomposition for $\alpha$ as well as using the usual identity $d(\omega \wedge \eta) = d\omega \wedge \eta + (-1)^{|\omega|} \omega \wedge d \eta$ for forms $\omega  \in \Omega_\tau^* (M)$ and $\eta \in \Omega^* (M)$, taking note of the fact that the signs of the second term all agree as $u$ is even degree.
\end{proof}
\end{lem}

The differential 
\begin{equation}
d_\Theta = ud_{dR} + u^2 \Theta \wedge \cdot : \Omega_k (M)[u] \rightarrow \Omega_{k+1} (M)[u]
\end{equation}
is of degree $+1$ due to the grading on $\Omega_k (M) [u]$.  The \textit{$\Theta$-twisted de Rham complex of $M$} $(\Omega_* (M)[u] , d_\Theta)$ is the complex of polynomials in $u$ over differential forms on $M$ with the differential $d_\Theta$.  If $\Theta ' = \Theta + d \eta$ is cohomologous to $\Theta$ then the complexes $(\Omega_* (M)[u] , d_\Theta)$ and $(\Omega_* (M)[u], d_{\Theta'})$ are isomorphic via the isomorphism
\begin{equation}
I_\eta : \xi \mapsto e^{-u \eta} \wedge \xi
\end{equation}

The cohomology of $(\Omega_* (M) [u] , d_\Theta)$ is called the $\Theta$-twisted cohomology of $M$, $H^*_\Theta (M)$.

\subsection{Twisted simplicial complex}

Now we will similarly develop the twisted simplicial cohomology of the simplicial manifold $M_\bullet$ in which the twisting of the complex of compatible forms is given by a closed simplicial $3$-form $\Theta \in \Omega^3 (M_\bullet)$. This definition involves a modification of the bicomplex of compatible forms \eqref{eq:bicompatible_form}. 

First we adjust the grading of $\Omega^{r,s} (M_\bullet)$ to define a bicomplex $\Omega_-^{r , s} (M_\bullet)$.  Let
\begin{equation}
\Omega_-^{r,s} (M_\bullet) [u] = \left( \coprod_{n \geq 0} \Omega_r (M_n) [u] \otimes \Omega^s (\Delta^n) \right) / \sim
\end{equation}
where $\Omega_r (M_n)$ is as in \eqref{eq:twisted_complex_grading}.  Here the equivalence relation $\sim$ is defined in precisely the same way as the compatibility condition in \eqref{eq:bicompatibility_cond}.  For $\alpha^r_{(n)} \in \Omega_r (M_n)[u]$ and $\beta^s_{(n)} \in \Omega^s (\Delta^n)$, we have $\alpha_{(n)}^r \otimes \beta_{(n)}^s \sim \alpha_{(n-1)}^r \otimes \beta_{(n-1)}^s$ if and only if the corresponding statement to \eqref{eq:bicompatibility_cond} holds.

Furthermore, there is a degree $+1$ differential $u\tilde{d}_{dR} : \Omega_-^{r,s} (M_\bullet) [u]\rightarrow \Omega_-^{r+1,s} (M_\bullet) [u]$ induced by the exterior derivative on each $M_n$,
\begin{equation}
(-1)^{r+s}ud \otimes \id : \Omega_r (M_n)[u] \otimes \Omega^s (\Delta^n) \rightarrow \Omega_{r+1} (M_n)[u] \otimes \Omega^s (\Delta^n)
\end{equation}
and there is a degree $+1$ differential $d_\Delta$ induced by the exterior derivative on each $\Delta^n$, i.e.\
\begin{equation}
d_{\Delta} =  \id \otimes d : \Omega_r (M_n)[u] \otimes \Omega^s (\Delta^n) \rightarrow \Omega_r (M_n)[u] \otimes \Omega^{s+1} (\Delta^n).
\end{equation}
As in \eqref{eq:bicompatible_sum}, we define
\begin{equation}
\Omega^k_- (M_\bullet) [u] = \bigoplus_{r+s=k} \Omega_-^{r,s} (M_\bullet)[u].
\end{equation}
Denote $\Omega^{r,s}_- (M_{(n)}) = \Omega^{r,s}_- (M_\bullet) [u] |_{M_n \times \Delta^n} $ and $\Omega^k_- (M_{(n)}) = \Omega^k (M_\bullet) |_{M_n \times \Delta^n}$.


\begin{prop}
The complex $(\Omega^*_- (M_\bullet)[u] , ud + d_\Delta )$ is quasi-isomorphic to the complex $(\Omega^* (M_\bullet) , d_1 + d_2)$.
\end{prop}

\begin{defn}
Given a simplicial differential $3$-form $\Theta \in \Omega^3 (M_\bullet)$ as in Definition \ref{def:compatible_form} we may decompose $\Theta$ with respect to the sum \eqref{eq:bicompatible_sum} as $\Theta = \Theta^{3,0} + \Theta^{2,1} + \Theta^{1,2} + \Theta^{0,3}$ where $\Theta^{j,3-j} \in \Omega^{j,3-j} (M_\bullet)$ for $0 \leq j \leq 3$.  If, further, $\Theta^{0,3} = 0$ then define a rescaling of $\Theta$
\begin{equation}
\label{eq:Theta_rescaled}
\Theta_u = u^2 \Theta^{3,0} + u \Theta^{2,1} + \Theta^{1,2} \in \Omega^* (M_\bullet)[u].
\end{equation}
Hence, exterior multiplication by $\Theta_u$ is an operator of degree $+1$ on
\begin{equation}
\Theta_u \wedge \cdot : \Omega_-^k (M_\bullet)[u] \rightarrow \Omega_-^{k+1} (M_\bullet)[u].
\end{equation}
With these operators we define the \textit{$\Theta$-twisted complex of compatible forms} to be $(\Omega^*_- (M_\bullet)[u] , u\tilde{d}_{dR} + d_\Delta - \Theta_u \wedge \cdot)$.
\end{defn}

\begin{prop}
Let $\Theta , \Theta ' \in \Omega^3 (M_\bullet)$ be compatible forms such that the components in $\Omega^{0,3} (M_\bullet)$ satisfy $(\Theta)^{0,3} = (\Theta')^{0,3}  = 0$.  If $\Theta - \Theta ' = d \eta$ for some $\eta \in \Omega^2 (M_\bullet)$ that satisfies $\eta^{0,2}  = 0$ then $(\Omega^*_- (M_\bullet)[u] , ud + d_\Delta + \Theta_u \wedge \cdot)$ and $(\Omega^*_- (M_\bullet)[u] , ud + d_\Delta + \Theta'_u \wedge \cdot)$ are isomorphic via the isomorphism
\begin{equation}
I_\eta : \xi \mapsto e^{ \eta_u} \wedge \xi
\end{equation}
where $\eta_u = u \eta^{2,0} + \eta^{1,1}$.
\end{prop}

Now let us consider the special case of $M_\bullet = (M \rtimes \Gamma)_\bullet$ in the following lemmas.

\begin{lem}
\label{lem:wedge_product_twisted_simplicial_complex}
For $\omega_1 \in \Omega^k_- (M_\bullet)[u]$ and $\omega_2 \in \Omega^\ell (M_\bullet)[u]$ the induced wedge product satisfies
\begin{equation*}
(-1)^{|\omega_2|} (u \tilde{d}_{dR} + d_\Delta) (\omega_1 \wedge \omega_2) = ( (u \tilde{d}_{dR} + d_\Delta) \omega_1) \wedge \omega_2 +  \omega_1 \wedge ( (u d_{dR} + d_\Delta) \omega_2)
\end{equation*}
\begin{proof}
As the variable $u$ is of even degree and will not alter any signs, let $\tilde{d} = \tilde{d}_{dR} + d_\Delta$ and let $d = d_{dR} + d_\Delta$. Suppose $\omega_1 \in \Omega^{r_1,s_1}_- (M_\bullet)[u]$ and $\omega_2 \in \Omega^{r_2 , s_2} (M_\bullet)[u]$ so that
\begin{align*}
\tilde{d} (\omega_1 \wedge \omega_2) &= (-1)^{\dim M - r_1 + r_2 + s_1 + s_2} d (\omega_1 \wedge \omega_2) \\
&= (-1)^{\dim M - r_1 + r_2 + s_1 + s_2} ( d \omega_1 \wedge \omega_2 + (-1)^{\dim M - r_1 + s_1} \omega_1 \wedge d \omega_2)\\
&= (-1)^{\dim M - r_1 + r_2 + s_1 + s_2} \left( (-1)^{\dim M - r_1 + s_1} \tilde{d} \omega_1 \wedge \omega_2 \right. \\
&\quad \left. + (-1)^{\dim M - r_1 + s_1} \omega_1 \wedge d \omega_2 \right)\\
&= (-1)^{-r_2 + s_2} \tilde{d} \omega_1 \wedge \omega_2 + (-1)^{-r_2 + s_2} \omega_1 \wedge d\omega_2
\end{align*}
Since $|\omega_2| = r_2+s_2$ which is $s_2 - r_2$ modulo 2,
\begin{equation*}
\tilde{d} (\omega_1 \wedge \omega_2) = (-1)^{|\omega_2|} \tilde{d} \omega_1 \wedge \omega_2 + (-1)^{|\omega_2|} \omega_1 \wedge d\omega_2
\end{equation*}
from which the lemma follows.
\end{proof}
\end{lem}

\begin{lem}
\label{lem:Stokes_twisted_simplicial_complex}
For $\alpha_{(n)} \otimes \beta_{(n)} \in \Omega_{0} (M_n) [u] \otimes \Omega^{n-1} (\Delta^n)$ we have the following
\begin{equation*}
\int_M \int_{\Delta^n} (u\tilde{d}_{dR} + d_\Delta) (\alpha_{(n)} \otimes \beta_{(n)}) = (-1)^{n-1} \int_M \alpha_{(n)} \int_{\partial (\Delta^n)} \beta_{(n)}.
\end{equation*}
Also, if $\beta_{(n)}$ is not of degree $n-1$ or $\alpha_{(n)}$ is not of degree $0$ in $\Omega^*_- (M_\bullet)[u]$ then the integral over $M \times \Delta^n$ is zero.
\begin{proof}
Since $\alpha_{(n)}$ is in $\Omega_0 (M_n)$, $\alpha_{(n)} = (\alpha_{(n)})_0$ as $\alpha_{(n)}$ must be an $m$-form on $M_n$.   For any $g_1 , \ldots , g_n \in \Gamma$ denote $\alpha = \alpha_{(n)} (g_1 , \ldots , g_n)$. Then 
\begin{align*}
\int_M \int_{\Delta^n} (u\tilde{d}_{dR} + d_\Delta) (\alpha \otimes \beta_{(n)}) &= \int_M \int_{\Delta^n} (u \tilde{d}_{dR} \alpha) \otimes \beta_{(n)} + \int_M \int_{\Delta^n} \alpha \otimes (d_\Delta \beta_{(n)}) \\
&=  \int_M \alpha  \int_{\Delta^n}  d_\Delta \beta_{(n)} \\
&=  \int_M \alpha  \int_{\Delta^n}  d \beta_{(n)} \\
&=  \int_M \alpha  \int_{\partial (\Delta^n)} \beta_{(n)}
\end{align*}
by Stokes' theorem.
\end{proof}
\end{lem}

Combining the previous two lemmas we have
\begin{lem}
\label{lem:Stokes_wedge_twisted_simplicial_complex}
For any $\omega_1 \in \Omega^*_- (M_\bullet)$ and $ \omega_2 \in \Omega^* (M_\bullet)$ 
\begin{multline*}
\int_M \int_{\Delta^n} ((u d_{dR} + d_\Delta) \omega_1) \wedge \omega_2  \\ = - \int_M \int_{\Delta^n} \omega_1 \wedge ( (ud_{dR} + d_\Delta) \omega_2) + (-1)^{|\omega_2|} \int_M \int_{\partial (\Delta^n)} \omega_1 \wedge \omega_2
\end{multline*}
where $\omega_2$ is of fixed degree.
\begin{proof}
We have
\begin{align*}
&\int_M \int_{\Delta^n} ((u d_{dR} + d_\Delta) \omega_1) \wedge \omega_2 \\
&= -  \int_M \int_{\Delta^n} \omega_1 \wedge ( (u d_{dR} + d_\Delta) \omega_2) + (-1)^{|\omega_2|}  \int_M \int_{\Delta^n} (u d_{dR} + d_\Delta) (\omega_1 \wedge \omega_2) \\
&= -  \int_M \int_{\Delta^n} \omega_1 \wedge ( (u d_{dR} + d_\Delta) \omega_2) + (-1)^{|\omega_2|} \int_M \int_{\partial(\Delta^n)} \omega_1 \wedge \omega_2
\end{align*}
where the first inequality follows from Lemma \ref{lem:wedge_product_twisted_simplicial_complex} and the second from Lemma \ref{lem:Stokes_twisted_simplicial_complex}.
\end{proof}
\end{lem}

\chapter{Simplicial Dixmier-Douady Form}
\label{chap:3}

In this chapter we will construct a 3-form $\Theta \in \Omega^3 (M_\bullet)$ that is a representative of the Dixmier-Douady class for $M_\bullet = (M \rtimes \Gamma)_\bullet$.

\section{Derivations $\nabla^k$}

Let $(L , \mu)$ be a gerbe on $M \rtimes \Gamma$ with connection $\nabla$ on $L \rightarrow M \rtimes \Gamma$.  Our goal now will be to construct a representative of the Dixmier-Douady class of the gerbe $(L, \mu)$ with connection $\nabla$ in the complex of compatible forms on $M_\bullet$.  We will again denote the discrepancy of $\nabla$ by $\alpha$.  For the entirety of the chapter, define the vector bundle $\pi : \E \rightarrow M$ where $\E = \bigoplus_{g \in \Gamma} L_g$ as in Example \ref{ex:direct_sum_bundle}.

Now let $p_k : M_{(k)} = M \times \Gamma^k \rightarrow M$ be defined by $p_k : (x, g_1 , \ldots , g_k) \mapsto x \in M$.  The pullback bundle $\E_k = p_k^* \E$ is a vector bundle over $M_{(k)}$.  There is also a connection on $\E_k$ defined by $\nabla^{\E_k} = p_k^* \nabla^{\E}$.  We will denote the connection $\nabla^{\E_k} |_{M \times \{g_1 \} \times \cdots \times \{g_k\}}$ by $\nabla^{\E_k} (g_1 , \ldots , g_k)$.  Note with this definition that $\nabla^{\E_k} (g_1 , \ldots , g_k) = \nabla^{\E}$ for all $g_1 , \ldots , g_k$.

\begin{rmk}
In general the collection $\{ \pi_k : \E_k \rightarrow M_{(k)} \}$ cannot be made into simplicial vector bundle using the maps \eqref{eq:simp_mfld_face_maps} and \eqref{eq:simp_mfld_deg_maps}.  In order for the diagram \eqref{eq:simp_bund_diag} of face maps to commute, we must have that $\delta_i^* \E_{k-1} \cong \E_{k}$ for $0 \leq i \leq k$ which does not hold in general.  Specifically, consider $\delta_0^* \E_{k-1}$.  For any $g_1 , \ldots , g_k \in \Gamma$ and $m \in M$,
\begin{align*}
\delta_0^* \E_{k-1} (x , g_1 , \ldots ,g_k) &= \delta_0^* (p_{k-1}^* \E) (x , g_1 ,\ldots , g_k) \\
&= (p_{k-1}^* \E) (xg_1 , g_2 ,\ldots , g_k) \\
&= \E^{g_1} \\
&\cong \E \otimes L_{g_1^{-1}},
\end{align*}
where the last isomorphism is given by $\varphi_{g_1^{-1}}$.
\end{rmk}

On the other hand, the vector bundles $\End \E_k \rightarrow M_{(k)}$ do form a simplicial vector bundle, hence we will denote $\End \E_k$ by ${\End \E}_{(k)}$.

\begin{prop}
The vector bundles $\pi_{(k)} : \End \E_{(k)} \rightarrow M_{(k)}$ form a simplicial vector bundle $\pi_\bullet : \End \E_\bullet \rightarrow M_\bullet$.
\begin{proof}
This follows as $\End (\E \otimes L_g) \cong \End \E$.  So $\delta_i^* \End \E_{k-1} \cong \End \E_k$ and the diagram \eqref{eq:simp_bund_diag} of face maps commutes for $0 \leq i \leq k$.
\end{proof}
\end{prop}

\begin{defn}
Denote $M_{(k)}^\Delta = M_{(k)} \times \Delta^k$ and $\E_k^\Delta = \E_k \times \Delta^k$.  Consider the bundle $ \pi_k \times \id :   \E_k^\Delta \rightarrow  M_{(k)}^\Delta $. Then define a connection $\nabla^{k} : C^\infty (M_{(k)}^\Delta , \E_k^\Delta) \rightarrow \Omega^1 (M_{(k)}^\Delta , \E_k^\Delta)$ by 
\begin{multline}
((\nabla^{k} (g_1 , \ldots , g_k) ) s) ( x_1 , \ldots , x_m, t_1 , \ldots , t_k ) \\
= ((\nabla^{\E_k} (g_1 , \ldots , g_k) + t_1 A(g_1) + \cdots + t_k A(g_1 \cdots g_k)) s ) ( x_1 , \ldots , x_m, t_1 , \ldots , t_k )
\end{multline}
where $s \in C^\infty (M_{(k)}^\Delta , \E_k^\Delta)$, $x_1 , \ldots , x_m$ are coordinates on $M$ (with $m = \dim M$), and $t_1 , \ldots , t_k$ are Cartesian coordinates on $\Delta^k$.

Here we use $\nabla^{\E_k}$ to denote the derivation on $\E_k^\Delta \rightarrow  M_{(k)}^\Delta$ that is the derivation $\nabla^{\E_k}$ on $\E_k \rightarrow M_{(k)}$ at every point of $\Delta^k$.  In other words, if $\bar{p}_k : M_{(k)}^\Delta \rightarrow M_{(k)}$ is the projection $(x,t) \mapsto x$ then we denote $\bar{p}_k^* \nabla^{\E_k}$ merely by $\nabla^{\E_k}$.  Note that, while $\nabla^k$ is a connection on $\E_k^\Delta \rightarrow  M_{(k)}^\Delta$, $\nabla^k$ is not a simplicial connection as $\E_k \rightarrow M_{(k)}$ does not form a simplicial bundle.

We may use the formal variable $u$ to rescale the connection $\nabla^k$ by considering $(\nabla^k)^{1,0}$ the component of the operator $\nabla^k$ 
\begin{equation*}
(\nabla^k)^{1,0} : C^\infty (M_{(k)}^\Delta , \E_k^\Delta) \rightarrow \Omega^1 (M_{(k)} , \E_k) \otimes C^\infty (\Delta^k , \Delta^k)
\end{equation*}
and $(\nabla^k)^{0,1} = d_\Delta$ the component that maps
\begin{equation*}
(\nabla^k)^{0,1} : C^\infty (M_{(k)}^{\Delta} , \E_k^\Delta) \rightarrow C^\infty (M_{(k)} , \E_k) \otimes \Omega^1 (\Delta^k, \Delta^k).
\end{equation*}
Let $\nabla^k_u = (\nabla^k_u)^{1,0} + u^{-1} d _\Delta$.

We may extend $\nabla^k$ to an operator
\begin{equation}
\nabla^k : \Omega^\ell (M_{(k)}^\Delta , \E_k^\Delta) \rightarrow \Omega^{\ell + 1} (M_{(k)}^\Delta , \E_k^\Delta)
\end{equation}
in the usual way.  Given a section $s \in C^\infty (M_{(k)}^\Delta , \E_k^\Delta)$ and a form $\omega \in C^\infty (M_{(k)}^\Delta)$ we define
\begin{equation}
\nabla^k (\omega \otimes s) = d \omega \otimes s + (-1)^\ell \omega \wedge \nabla^k s.
\end{equation}
\end{defn}

\begin{rmk}
The operator $\nabla^k$ may be defined to act on the algebra of $\End \E$-valued forms $\Omega^* (M , \End \E)$ and in particular on sections $C^\infty (M , \End \E)$ as follows.  Given a section $\eta \in \Omega^* (M , \End \E)$, define
\begin{equation}
\nabla^k \eta = [\nabla^k , \bar{p}_k^* p_k^* \eta]
\end{equation}
where $p_k : M_{(k)} \rightarrow M$ is the projection $(x, g_1 , \ldots , g_k) \mapsto x$, $\bar{p}_k : M_{(k)} \times \Delta^k \rightarrow M_{(k)}$ is the natural projection $(x,t) \mapsto x$ and we consider $\nabla^k$ to be degree $+1$ with respect to the graded commutator. In other words, $\nabla^k \eta$ operates on sections $s \in C^\infty (M_{(k)}^\Delta , \E_k^\Delta)$ by
\begin{equation}
(\nabla^k \eta) s = [\nabla^k , \bar{p}_k^* p_k^* \eta] s = \nabla^k ((\bar{p}_k^* p_k^* \eta) s) - (-1)^{|\eta| + 1} (\bar{p}_k^* p_k^* \eta) (\nabla^k s).
\end{equation}
\end{rmk}

With this extension the rescaling $\nabla^k_u =  (\nabla^k)^{1,0} + u^{-1} d_\Delta$ still makes sense.

\begin{lem}
\label{lem:trace_connection}
For any $\eta \in \Omega^* (M, \End \E)$,
\begin{equation}
(u d_{dR} + d_\Delta) \tr (\eta) = u\tr (\nabla^k_u (\eta))
\end{equation}
where by $\tr (\eta)$ we mean $\tr(\bar{p}_k^* p_k^* \eta)$.
\begin{proof}
We will omit the pullback maps $\bar{p}_k^* p_k^*$ as they should be clear from the context.  As $\nabla^k_u =  (\nabla^k)^{1,0} + u^{-1}d_\Delta$ and $(\nabla^k)^{1,0} = d_{dR} + \omega$ locally for some $\omega \in \Omega^1 (M, \End \E)$.  For a local section $s \in C^\infty (M_{(k)}^\Delta , \E_k^\Delta)$
\begin{align*}
((\nabla^k)^{1,0} \eta) (s) &= [(\nabla^k)^{1,0} , \eta] s \\
&= [d_{dR}  + \omega , \eta] s \\
&= d_{dR} (\eta s) - (-1)^{|\eta|} \eta d_{dR} s + [\omega , \eta] s \\
&= (d_{dR} \eta + [\omega , \eta]) s
\end{align*}
So by Lemma \ref{lem:trace_commutator}
\begin{align*}
\tr ((\nabla^k)^{1,0} \eta) &= \tr ( (d_{dR} \eta +  [\omega , \eta])) \\
&= \tr (d_{dR} \eta) + \tr([\omega , \eta]) \\
&= \tr (d_{dR} \eta) \\
&= d_{dR} \tr (\eta).
\end{align*}
Hence
\begin{equation*}
u\tr ( ( (\nabla^k)^{1,0} + u^{-1}d_\Delta) \eta) = u \tr ( (\nabla^k)^{1,0} \eta) + \tr (d_\Delta \eta) = u d_{dR} \tr(\eta) + d_{\Delta} \tr (\eta),
\end{equation*}
as required. 
\end{proof}
\end{lem}

\section{Simplicial 2-form}

The collection of $\End \E_{(k)}$-valued $2$-forms on $M_{(k)} \times \Delta^k$ given by $(\nabla^k)^2$ for each $k \geq 0$ does not define a simplicial differential form on $M_\bullet$.  In the theorem below we adjust $(\nabla^k)^2$ by a scalar-valued form to obtain a simplicial differential $2$-form on $M_\bullet$.

\begin{thm}
\label{thm:simplicial_2_form}
Let $(L, \mu)$ be a gerbe on $M \rtimes \Gamma$ and let $\nabla$ be a connection on $L \rightarrow M \rtimes \Gamma$.  For each $k \geq 0$, let $\vartheta_{(k)} \in \Omega^2 ( M_{(k)}^\Delta, \End \E_{(k)}^\Delta)$ be defined by the formula
\begin{multline}
\label{eq:theta_def}
\vartheta_{(k)} (g_1 , \ldots , g_k) = (\nabla^k (g_1, \ldots , g_k))^2  \\ - \sum_{i=1}^k t_i \theta_{g_1 \cdots g_i}  + \sum_{1 \leq i < j \leq k} \alpha(g_1 \cdots g_i , g_{i+1} \cdots g_j) (t_i dt_j - t_j dt_i).
\end{multline}
Then $\vartheta = \{ \vartheta_{(k)} \}$ is a compatible $2$-form on $M_\bullet$ with values in $\End \E_\bullet$, which means that $\vartheta \in \Omega^2 (M_\bullet , \End \E_\bullet)$. We call $\vartheta$ the \textit{simplicial $2$-form associated to $\nabla$}.
\begin{proof}
Our goal is to check the compatibility conditions of \eqref{eq:compatibility_cond}. First it will be helpful to rewrite $\vartheta_{(k)} (g_1 , \ldots , g_k)$. With a few algebraic calculations we obtain
\begin{align*}
 \vartheta_{(k)}(g_1 , \ldots , g_k) &= (\nabla^{\E_k} (g_1 , \ldots , g_k) + t_1 A(g_1) + t_2 A(g_1 g_2) + \cdots + t_k A(g_1 \cdots g_k))^2 \displaybreak[0]\\
&\quad - \sum_{i=1}^k t_i \theta_{g_1 \cdots g_i} + \sum_{1 \leq i < j \leq k} \alpha(g_1 \cdots g_i , g_{i+1} \cdots g_j) (t_i dt_j - t_j dt_i) \displaybreak[2]\\
&= (\nabla^{\E})^2 + \sum_{i=1}^k t_i^2 A(g_1 \cdots g_i)^2 +  [\nabla^{\E}  , t_i A(g_1 \cdots g_i) ]  - t_i \theta_{g_1 \cdots g_i} \displaybreak[0]\\
& \quad + \sum_{1 \leq i < j \leq k} [t_i A(g_1 \cdots g_i) , t_j A(g_1 \cdots g_j)] \\& \qquad \qquad+ \alpha(g_1 \cdots g_i , g_{i+1} \cdots g_j) (t_i dt_j - t_j dt_i) \displaybreak[2]\\
\intertext{and by expanding $[ \nabla^{\E}  , t_i A(g_1 \cdots g_i) ] = dt_i A(g_1 \cdots g_i) +  t_i [ \nabla^{\E}  , A(g_1 \cdots g_i) ]$}
&= \theta^\E  + \sum_{i=1}^k t_i^2 A(g_1 \cdots g_i)^2 + dt_i A(g_1 \cdots g_i) \\
&\qquad \qquad+ t_i [ \nabla^{\E}  , A(g_1 \cdots g_i) ]  - t_i \theta_{g_1 \cdots g_i} \displaybreak[0]\\
& \quad + \sum_{1 \leq i < j \leq k} [t_i A(g_1 \cdots g_i) , t_j A(g_1 \cdots g_j)] \\
&\qquad \qquad+ \alpha(g_1 \cdots g_i , g_{i+1} \cdots g_j) (t_i dt_j - t_j dt_i) \displaybreak[2]\\
&= \left(1 - \sum_{i = 1}^k t_i  \right)  \theta^\E + \sum_{i=1}^k  t_i \theta^\E  + t_i^2  A(g_1 \cdots g_i)^2 + dt_i A(g_1 \cdots g_i) \\
&\qquad \qquad+ t_i [ \nabla^{\E}  , A(g_1 \cdots g_i) ]  - t_i \theta_{g_1 \cdots g_i} \displaybreak[0]\\
& \quad + \sum_{1 \leq i < j \leq k} [t_i A(g_1 \cdots g_i) , t_j A(g_1 \cdots g_j)] \\
&\qquad \qquad + \alpha(g_1 \cdots g_i , g_{i+1} \cdots g_j) (t_i dt_j - t_j dt_i) \displaybreak[2]\\
&= \left(1 - \sum_{i = 1}^k t_i  \right) \theta^\E  \displaybreak[0]\\
&\quad + \sum_{i=1}^k  t_i \left( \theta^\E +  A(g_1 \cdots g_i)^2 + [ \nabla^{\E}  , A(g_1 \cdots g_i) ]  - \theta_{g_1 \cdots g_i} \right) \displaybreak[0]\\
&\qquad \qquad + dt_i A(g_1 \cdots g_i) - (1-t_i)t_i A(g_1 \cdots g_i)^2 \displaybreak[0]\\
& \quad + \sum_{1 \leq i < j \leq k} [t_i A(g_1 \cdots g_i) , t_j A(g_1 \cdots g_j)] \\
&\qquad \qquad+ \alpha(g_1 \cdots g_i , g_{i+1} \cdots g_j) (t_i dt_j - t_j dt_i) \displaybreak[2] \\
\intertext{and by Equation \eqref{eq:twisted_bundle_curvature} }
&= \left(1 - \sum_{i = 1}^k t_i  \right) \theta^\E  + \sum_{i=1}^k  t_i \varphi_{g_1\cdots g_i}^* \theta^{\E^{g_1 \cdots g_i}} \\
&\qquad \qquad+ dt_i A(g_1 \cdots g_i) - (1-t_i)t_i A(g_1 \cdots g_i)^2 \displaybreak[0]\\
& \quad + \sum_{1 \leq i < j \leq k} [t_i A(g_1 \cdots g_i) , t_j A(g_1 \cdots g_j)] \\
&\qquad \qquad+ \alpha(g_1 \cdots g_i , g_{i+1} \cdots g_j) (t_i dt_j - t_j dt_i). 
\end{align*}
So we have obtained
\begin{multline}
\label{eq:theta_rewritten}
\vartheta_{(k)}(g_1 , \ldots , g_k) \\ =  \left(1 - \sum_{i = 1}^k t_i  \right) \theta^\E  +  \sum_{i=1}^k  t_i \varphi_{g_1\cdots g_i}^* \theta^{\E^{g_1 \cdots g_i}} + dt_i A(g_1 \cdots g_i) - (1-t_i)t_i A(g_1 \cdots g_i)^2 \\+ \sum_{1 \leq i < j \leq k} [t_i A(g_1 \cdots g_i) , t_j A(g_1 \cdots g_j)] + \alpha(g_1 \cdots g_i , g_{i+1} \cdots g_j) (t_i dt_j - t_j dt_i) 
\end{multline}

As we are using Cartesian coordinates, when $1 \leq \ell \leq k$ we have that 
\begin{equation}
( \id \times \partial_\ell )^* \vartheta_{(k)} (g_1 , \ldots , g_k) = \vartheta_{(k)} (g_1 , \ldots , g_k) |_{t_\ell = 0}
\end{equation}
and when $\ell = 0$ we have that 
\begin{equation}
( \id \times \partial_0 )^* \vartheta_{(k)} (g_1 , \ldots , g_k) = \vartheta_{(k)} (g_1 ,\ldots , g_k) |_{t_1+ \cdots + t_k = 1}
\end{equation}
To show that $( \id \times \partial_\ell )^* \vartheta_{(k)} (g_1 , \ldots , g_k) = (\delta_\ell \times \id)^* \vartheta_{(k-1)}$ for $1 \leq \ell \leq k-1$ note that when $t_\ell = 0$, $dt_\ell = 0$.  From the above
\begin{align*}
 \vartheta_{(k)} (g_1 , \ldots , g_k) |_{t_\ell = 0} &= \theta^\E + \sum_{\substack{ i=1 \\ i \neq \ell} }^k t_i^2 A(g_1 \cdots g_i)^2 + [ \nabla^{\E}  , t_i A(g_1 \cdots g_i) ] - t_i \theta_{g_1 \cdots g_i} \displaybreak[0]\\
& \quad + \sum_{\substack{ 1 \leq i < j \leq k \\ i,j \neq \ell} } [t_i A(g_1 \cdots g_i) , t_j A(g_1 \cdots g_j)] \\&\qquad \qquad+ \alpha(g_1 \cdots g_i , g_{i+1} \cdots g_j) (t_i dt_j - t_j dt_i) \displaybreak[2]
\end{align*}
which may be written as $\vartheta_{(k-1)} (g_1 , \ldots , g_\ell g_{\ell +1} , \ldots , g_k)$ upon reindexing $t_{i+1}$ to $t_i$ when $\ell \leq i+1 \leq k$.

Now for the $\ell = 0$ case note that when $t_1 + \cdots + t_k = 1$ we have $dt_1 + \cdots + dt_k = 0$. In the above expression for $\vartheta_{(k)} (g_1 , \ldots , g_k)$ we set $t_1 + \cdots + t_k =1$ to obtain
\begin{align*}
& \vartheta(g_1 , \ldots , g_k)|_{\sum_{i=1}^k t_i = 1}= \displaybreak[0]\\
&= \left(1 - \sum_{i = 1}^k t_i  \right)  \theta^\E + \sum_{i=1}^k  t_i \varphi_{g_1\cdots g_i}^* \theta^{\E^{g_1 \cdots g_i}}  + dt_i A(g_1 \cdots g_i) - (1-t_i)t_i A(g_1 \cdots g_i)^2 \displaybreak[0]\\
& \quad + \sum_{1 \leq i < j \leq k} t_i t_j [ A(g_1 \cdots g_i) , A(g_1 \cdots g_j)] +  (t_i dt_j - t_j dt_i)  \alpha(g_1 \cdots g_i , g_{i+1} \cdots g_j)\displaybreak[2]\\
&=  \sum_{i=1}^k  t_i \varphi_{g_1\cdots g_i}^* \theta^{\E^{g_1 \cdots g_i}}  + dt_i A(g_1 \cdots g_i) - \left( \sum_{ j=1,  j\neq i}^k t_j \right)t_i A(g_1 \cdots g_i)^2 \displaybreak[0]\\
& \quad + \sum_{1 \leq i < j \leq k} t_i t_j [ A(g_1 \cdots g_i) , A(g_1 \cdots g_j)] +  (t_i dt_j - t_j dt_i)  \alpha(g_1 \cdots g_i , g_{i+1} \cdots g_j)\displaybreak[2]\\
&=  \sum_{i=1}^k  t_i \varphi_{g_1\cdots g_i}^* \theta^{\E^{g_1 \cdots g_i}}  + dt_i A(g_1 \cdots g_i) \displaybreak[0] \\
& \quad + \sum_{1 \leq i < j \leq k} t_i t_j \left( -A(g_1 \cdots g_i)^2 - A(g_1 \cdots g_j)^2 + [ A(g_1 \cdots g_i) , A(g_1 \cdots g_j)] \right)\\
&\qquad \qquad  +  (t_i dt_j - t_j dt_i)  \alpha(g_1 \cdots g_i , g_{i+1} \cdots g_j)\displaybreak[2]\\
&=  \sum_{i=1}^k  t_i \varphi_{g_1\cdots g_i}^* \theta^{\E^{g_1 \cdots g_i}}  + dt_i A(g_1 \cdots g_i) \displaybreak[0] \\
& \quad + \sum_{1 \leq i < j \leq k} - t_i t_j \left( A(g_1 \cdots g_i) - A(g_1 \cdots g_j) \right)^2 \displaybreak[0]\\
&\qquad \qquad +  (t_i dt_j - t_j dt_i)  \alpha(g_1 \cdots g_i , g_{i+1} \cdots g_j)\displaybreak[2]\\
\intertext{where by Lemma \ref{lem:discrepancy_identity}}
&=  \sum_{i=1}^k  t_i \varphi_{g_1\cdots g_i}^* \theta^{\E^{g_1 \cdots g_i}}  + dt_i A(g_1 \cdots g_i)  \displaybreak[0]\\
& \quad + \sum_{1 \leq i < j \leq k} -t_i t_j \left(-\alpha (g_1 \cdots g_i , g_{i+1} \cdots g_j) - \varphi_{g_1 \cdots g_i}^* A(g_{i+1} \cdots g_j)^{g_1 \cdots g_i} \right)^2 \displaybreak[0]\\
&\qquad \qquad +  (t_i dt_j - t_j dt_i)  \alpha(g_1 \cdots g_i , g_{i+1} \cdots g_j) \displaybreak[2]\\
\intertext{and because $\alpha$ is a scalar-valued 1-form we have}
&=  \sum_{i=1}^k  t_i \varphi_{g_1\cdots g_i}^* \theta^{\E^{g_1 \cdots g_i}}  + dt_i A(g_1 \cdots g_i)  \displaybreak[0]\\
& \quad + \sum_{1 \leq i < j \leq k} - t_i t_j \left( \varphi_{g_1 \cdots g_i}^* A(g_{i+1} \cdots g_j)^{g_1 \cdots g_i} \right)^2 \displaybreak[0] \\
& \qquad \qquad +  (t_i dt_j - t_j dt_i)  \alpha(g_1 \cdots g_i , g_{i+1} \cdots g_j) \displaybreak[2]\\
&=  \sum_{i=1}^k  t_i \varphi_{g_1\cdots g_i}^* \theta^{\E^{g_1 \cdots g_i}}  + dt_i A(g_1 \cdots g_i) \displaybreak[0] \\
& \quad + \sum_{1 \leq i < j \leq k} - t_i t_j \left( \varphi_{g_1 \cdots g_i}^* A(g_{i+1} \cdots g_j)^{g_1 \cdots g_i} \right)^2 \displaybreak[0] \\
& \quad +  \sum_{2 \leq i < j \leq k} (t_i dt_j - t_j dt_i)  \alpha(g_1 \cdots g_i , g_{i+1} \cdots g_j) + \sum_{j=2}^k (t_1 dt_j - t_j dt_1) \alpha (g_1 , g_2 \cdots g_j) \displaybreak[2]\\
&=  \sum_{i=1}^k  t_i \varphi_{g_1\cdots g_i}^* \theta^{\E^{g_1 \cdots g_i}} + \sum_{i=2}^k dt_i (A(g_1 \cdots g_i) - A(g_1))  \displaybreak[0] \\
& \quad + \sum_{1 \leq i < j \leq k} - t_i t_j \left( \varphi_{g_1 \cdots g_i}^* A(g_{i+1} \cdots g_j)^{g_1 \cdots g_i} \right)^2 \displaybreak[0] \\
& \quad +  \sum_{2 \leq i < j \leq k} (t_i dt_j - t_j dt_i)  \alpha(g_1 \cdots g_i , g_{i+1} \cdots g_j)\displaybreak[0] \\
& \quad + \sum_{j=2}^k \left( \left(1-\sum_{\ell=2}^k t_\ell \right) dt_j - t_j \left( \sum_{\ell=2}^k (-dt_\ell) \right) \right) \alpha (g_1 , g_2 \cdots g_j) \displaybreak[2]\\
&=  \sum_{i=1}^k  t_i \varphi_{g_1\cdots g_i}^* \theta^{\E^{g_1 \cdots g_i}}  + \sum_{i=2}^k dt_i (A(g_1 \cdots g_i) - A(g_1))  \displaybreak[0] \\
& \quad + \sum_{1 \leq i < j \leq k} - t_i t_j \left( \varphi_{g_1 \cdots g_i}^* A(g_{i+1} \cdots g_j)^{g_1 \cdots g_i} \right)^2  \displaybreak[0]\\
& \quad +  \sum_{2 \leq i < j \leq k} (t_i dt_j - t_j dt_i)  \alpha(g_1 \cdots g_i , g_{i+1} \cdots g_j) \displaybreak[0]\\
& \quad + \sum_{j,\ell=2, j \neq \ell}^k (t_j dt_\ell - t_\ell dt_j) \alpha (g_1 , g_2 \cdots g_j) + \sum_{j=2}^k dt_j \alpha(g_1 , g_2 \cdots g_j) \displaybreak[2]\\
&=  \sum_{i=1}^k  t_i \varphi_{g_1\cdots g_i}^* \theta^{\E^{g_1 \cdots g_i}}  + \sum_{i=2}^k dt_i (A(g_1 \cdots g_i) - A(g_1) - \alpha (g_1 , g_2 \cdots g_i)) \displaybreak[0] \\
& \quad + \sum_{1 \leq i < j \leq k} - t_i t_j \left( \varphi_{g_1 \cdots g_i}^* A(g_{i+1} \cdots g_j)^{g_1 \cdots g_i} \right)^2 \displaybreak[0] \\
& \quad +  \sum_{2 \leq i < j \leq k} (t_i dt_j - t_j dt_i)  (\alpha(g_1 , g_2 \cdots g_i) - \alpha( g_1 , g_2 \cdots g_j) + \alpha(g_1 \cdots g_i , g_{i+1} \cdots g_j)) \displaybreak[2]\\
\intertext{and by Lemma \ref{lem:discrepancy_identity} and Corollary \ref{cor:discrepancy_identity2} we have}
&=  \sum_{i=1}^k  t_i \varphi_{g_1\cdots g_i}^* \theta^{\E^{g_1 \cdots g_i}}  - \sum_{i=2}^k dt_i \varphi_{g_1}^* A(g_2 \cdots g_i)^{g_1} \displaybreak[0] \\
& \quad + \sum_{1 \leq i < j \leq k} - t_i t_j \left( \varphi_{g_1 \cdots g_i}^* A(g_{i+1} \cdots g_j)^{g_1 \cdots g_i} \right)^2 \displaybreak[0] \\
& \quad +  \sum_{2 \leq i < j \leq k} (t_i dt_j - t_j dt_i)  \varphi_{g_1}^* \alpha(g_2 \cdots g_i , g_{i+1} \cdots g_j)^{g_1} \displaybreak[2]\\
&=  \left( 1- t_2 - \cdots - t_k \right) \theta^\E + \sum_{i=2}^k  t_i \varphi_{g_1\cdots g_i}^* \theta^{\E^{g_1 \cdots g_i}}  - \sum_{i=2}^k dt_i \varphi_{g_1}^* A(g_2 \cdots g_i)^{g_1}  \displaybreak[0]\\
& \quad + \sum_{j=2}^k - \left( 1 - t_2 - \cdots - t_k \right) t_j (\varphi_{g_1}^* A(g_2 \cdots g_j)^{g_1})^2 \displaybreak[0]\\
& \quad + \sum_{2 \leq i < j \leq k} - t_i t_j \left( \varphi_{g_1 \cdots g_i}^* A(g_{i+1} \cdots g_j)^{g_1 \cdots g_i} \right)^2  \displaybreak[0]\\
& \quad +  \sum_{2 \leq i < j \leq k} (t_i dt_j - t_j dt_i)  \varphi_{g_1}^* \alpha(g_2 \cdots g_i , g_{i+1} \cdots g_j)^{g_1} \displaybreak[2]\\
&=  \left( 1- t_2 - \cdots - t_k \right) \theta^\E + \sum_{i=2}^k  t_i \varphi_{g_1\cdots g_i}^* \theta^{\E^{g_1 \cdots g_i}}  - \sum_{i=2}^k dt_i \varphi_{g_1}^* A(g_2 \cdots g_i)^{g_1} \displaybreak[0] \\
& \quad + \sum_{j=2}^k - \left( 1 -  t_j \right) t_j (\varphi_{g_1}^* A(g_2 \cdots g_j)^{g_1})^2 \displaybreak[0]\\
& \quad + \sum_{2 \leq i < j \leq k} - t_i t_j \left( \varphi_{g_1 \cdots g_i}^* (A(g_{i+1} \cdots g_j)^{g_1 \cdots g_i})^2 \right. \displaybreak[0] \\
&\qquad \qquad \left. - (\varphi_{g_1}^* A(g_2 \cdots g_i)^{g_1})^2 - (\varphi_{g_1}^* A(g_2 \cdots g_j)^{g_1})^2 \right)  \displaybreak[0]\\
& \quad +  \sum_{2 \leq i < j \leq k} (t_i dt_j - t_j dt_i)  \varphi_{g_1}^* \alpha(g_2 \cdots g_i , g_{i+1} \cdots g_j)^{g_1}\displaybreak[2] \\
&=  \left( 1- t_2 - \cdots - t_k \right) \theta^\E + \sum_{i=2}^k  t_i \varphi_{g_1\cdots g_i}^* \theta^{\E^{g_1 \cdots g_i}}  - \sum_{i=2}^k dt_i \varphi_{g_1}^* A(g_2 \cdots g_i)^{g_1} \displaybreak[0] \\
& \quad + \sum_{j=2}^k - \left( 1 -  t_j \right) t_j (\varphi_{g_1}^* A(g_2 \cdots g_j)^{g_1})^2 \displaybreak[0]\\
& \quad + \sum_{2 \leq i < j \leq k} - t_i t_j \left( \varphi_{g_1 \cdots g_i}^* (A(g_{i+1} \cdots g_j)^{g_1 \cdots g_i})^2 \right. \\
&\qquad \qquad \left. - (\varphi_{g_1}^* A(g_2 \cdots g_i)^{g_1})^2 - (\varphi_{g_1}^* A(g_2 \cdots g_j)^{g_1})^2 \right)  \displaybreak[0]\\
& \quad + \sum_{2 \leq i < j \leq k} t_i t_j \left([\varphi_{g_1}^* A(g_2 \cdots g_i)^{g_1} , \varphi_{g_1}^* A(g_2 \cdots g_j)^{g_1} ] \right. \\
&\qquad \qquad \left. - [\varphi_{g_1}^* A(g_2 \cdots g_i)^{g_1} , \varphi_{g_1}^* A(g_2 \cdots g_j)^{g_1} ] \right) \displaybreak[0]\\
& \quad +  \sum_{2 \leq i < j \leq k} (t_i dt_j - t_j dt_i)  \varphi_{g_1}^* \alpha(g_2 \cdots g_i , g_{i+1} \cdots g_j)^{g_1} \displaybreak[2]\\
&=  \left( 1- t_2 - \cdots - t_k \right) \theta^\E \\
&\quad + \sum_{i=2}^k  t_i \varphi_{g_1\cdots g_i}^* \theta^{\E^{g_1 \cdots g_i}}  -  dt_i \varphi_{g_1}^* A(g_2 \cdots g_i)^{g_1} - \left( 1 -  t_i \right) t_i (\varphi_{g_1}^* A(g_2 \cdots g_i)^{g_1})^2 \displaybreak[0] \\
& \quad + \sum_{2 \leq i < j \leq k} - t_i t_j \left( \varphi_{g_1 \cdots g_i}^* (A(g_{i+1} \cdots g_j)^{g_1 \cdots g_i})^2 \right. \displaybreak[0] \\
& \qquad \qquad \left. - (\varphi_{g_1}^* A(g_2 \cdots g_i)^{g_1} - \varphi_{g_1}^* A(g_2 \cdots g_j)^{g_1})^2 \right) \displaybreak[0] \\
& \quad + \sum_{2 \leq i < j \leq k} - t_i t_j  [\varphi_{g_1}^* A(g_2 \cdots g_i)^{g_1} , \varphi_{g_1}^* A(g_2 \cdots g_j)^{g_1} ]  \displaybreak[0] \\
&\qquad \qquad  +(t_i dt_j - t_j dt_i)  \varphi_{g_1}^* \alpha(g_2 \cdots g_i , g_{i+1} \cdots g_j)^{g_1} \displaybreak[2]\\
\intertext{and using Lemma \ref{lem:discrepancy_identity} again we have}
&=  \left( 1- t_2 - \cdots - t_k \right) \theta^\E \displaybreak[0]\\
&\quad + \sum_{i=2}^k  t_i \varphi_{g_1\cdots g_i}^* \theta^{\E^{g_1 \cdots g_i}}  -  dt_i \varphi_{g_1}^* A(g_2 \cdots g_i)^{g_1} - \left( 1 -  t_i \right) t_i (\varphi_{g_1}^* A(g_2 \cdots g_i)^{g_1})^2 \displaybreak[0]\\
& \quad + \sum_{2 \leq i < j \leq k} - t_i t_j \left( \varphi_{g_1 \cdots g_i}^* (A(g_{i+1} \cdots g_j)^{g_1 \cdots g_i})^2 \right. \\
& \qquad \qquad \left. - (-\alpha (g_2 \cdots g_i , g_{i+1} \cdots g_j) -\varphi_{g_1 \cdots g_i}^* A(g_{i+1} \cdots g_j)^{g_1 \cdots g_i})^2 \right)  \displaybreak[0]\\
& \quad + \sum_{2 \leq i < j \leq k} - t_i t_j  [\varphi_{g_1}^* A(g_2 \cdots g_i)^{g_1} , \varphi_{g_1}^* A(g_2 \cdots g_j)^{g_1} ]  \displaybreak[0]\\
&\qquad \qquad + (t_i dt_j - t_j dt_i)  \varphi_{g_1}^* \alpha(g_2 \cdots g_i , g_{i+1} \cdots g_j)^{g_1} \displaybreak[2]\\
&=  \left( 1- t_2 - \cdots - t_k \right) \theta^\E \displaybreak[0] \\
&\quad + \sum_{i=2}^k  t_i \varphi_{g_1\cdots g_i}^* \theta^{\E^{g_1 \cdots g_i}}  -  dt_i \varphi_{g_1}^* A(g_2 \cdots g_i)^{g_1} - \left( 1 -  t_i \right) t_i (\varphi_{g_1}^* A(g_2 \cdots g_i)^{g_1})^2 \displaybreak[0]\\
& \quad + \sum_{2 \leq i < j \leq k} - t_i t_j  [\varphi_{g_1}^* A(g_2 \cdots g_i)^{g_1} , \varphi_{g_1}^* A(g_2 \cdots g_j)^{g_1} ] \displaybreak[0] \\
&\qquad \qquad + (t_i dt_j - t_j dt_i)  \varphi_{g_1}^* \alpha(g_2 \cdots g_i , g_{i+1} \cdots g_j)^{g_1} \displaybreak[0] \\
&= \varphi_{g_1}^* \theta_{(k-1)} (g_2,  \ldots , g_k)^{g_1}.
\end{align*}
which is the form of $\vartheta_{(k-1)}$ we obtained in \eqref{eq:theta_rewritten}.
\end{proof}
\end{thm}

\begin{rmk}
If we define $\vartheta_{(k)}'$ by $(\nabla^k)^2 = [\vartheta_{(k)} ' , \cdot]$ then the discrepancy between $\vartheta_{(k)}'$ and the simplicial 2-form associated to $\nabla$, $\vartheta_{(k)}$, i.e.\ $(\vartheta_{(k)}' - \vartheta_{(k)}) (g_1 , \ldots , g_k) $, is given by
\begin{equation}
\sum_{i=1}^k t_i \theta_{g_1 \cdots g_i}  - \sum_{1 \leq i < j \leq k} \alpha(g_1 \cdots g_i , g_{i+1} \cdots g_j) (t_i dt_j - t_j dt_i)
\end{equation}
which is a scalar-valued form.  In other words, $\vartheta_{(k)} ' - \vartheta_{(k)} \in \Omega^2 (M_{(k)}^\Delta )$.  Since a scalar valued form in $\Omega^2 (M_{(k)}^\Delta)$ is central in $\Omega^2 (M_{(k)}^\Delta , \End \E_{(k)}^\Delta)$, $[\vartheta_{(k)}' , \cdot] = [\vartheta_{(k)}, \cdot]$.  Hence
\begin{equation}
\label{eq:nabla_squared_equals_theta}
(\nabla^k)^2 (a) = [\vartheta_{(k)} , a]
\end{equation}
for any $a \in C^\infty (M_{(k)}^\Delta, \End \E_{(k)}^\Delta)$.
\end{rmk}

\begin{defn}
Notice from the previous formula for $\vartheta$ that $\vartheta^{0,2} = 0$. We will then define $\vartheta_u = u \vartheta^{2,0} + \vartheta^{1,1}$.  With this rescaling we can view $\vartheta_u$ as an element of $\Omega^0_- (M_\bullet)[u]$.  
\end{defn}

\section{Simplicial Dixmier-Douady form}

\begin{thm}
Let $(L, \mu)$ be a gerbe on $M \rtimes \Gamma$ and let $\nabla$ be a connection on $L \rightarrow M \rtimes \Gamma$ with simplicial $2$-form $\vartheta$ as in Theorem \ref{thm:simplicial_2_form}.  Let $\Theta_{(k)} = \nabla^k \vartheta_{(k)}$. The collection of $3$-forms $\Theta_{(k)}$ define a simplicial $3$-form $\{ \Theta_{(k)} \} \in \Omega^3 (M_\bullet)$.  Furthermore, $\Theta_{(k)}^{0,3} = 0$ as an element of $\Omega^{r,s} (M_\bullet)$. 
\begin{proof}
First let us provide a formula for $\Theta_{(k)}$. Using the Bianchi identity
\begin{align*}
\Theta_{(k)} (g_1, \ldots , g_k) &= (\nabla^k \vartheta_{(k)} ) (g_1 , \ldots , g_k) \displaybreak[2]\\
&= (\nabla^k)(g_1 , \ldots , g_k) \left( (\nabla^k)^2 (g_1 , \ldots , g_k) - \sum_{i=1}^k t_i \theta_{g_1 \cdots g_i} \right. \displaybreak[0]\\
&\qquad \qquad  \left. + \sum_{1 \leq i < j \leq k} \alpha ( g_1 \cdots g_i , g_{i+1} \cdots g_j ) (t_i dt_j - t_j dt_i ) \right) \displaybreak[2]\\
&=  (\nabla^k)^3 (g_1 , \ldots , g_k) -  (\nabla^k (g_1 , \ldots , g_k) ) \left(  \sum_{i=1}^k t_i \theta_{g_1 \cdots g_i}  \right. \displaybreak[0]\\ 
&\qquad \qquad \left. + \sum_{1 \leq i < j \leq k} \alpha ( g_1 \cdots g_i , g_{i+1} \cdots g_j ) (t_i dt_j - t_j dt_i ) \right) \displaybreak[0]\\
&= -\sum_{i=1}^k dt_i \theta_{g_1 \cdots g_i} \\
&\quad + \sum_{1 \leq i < j \leq k} d\alpha (g_1 \cdots g_i, g_{i+1} \cdots g_j) (t_i dt_j - t_j dt_i) \\
&\qquad \qquad + 2 \alpha (g_1 \cdots g_i , g_{i+1} \cdots g_j) dt_i dt_j,
\end{align*}
since $\theta_g$ and $\alpha(g , g')$ are scalar-valued forms.  Since there are no terms containing $dt_{i_1} dt_{i_2} dt_{i_3}$ where $1 \leq i_1 , i_2 , i_3 \leq k$, we see that $\Theta_{(k)}^{0,3} = 0$.

Using this formula, we can see $\Theta_{(k)}$ is simplicial by a direct calculation. First we check that $\Theta_{(k)} (g_1 , \ldots , g_k) |_{t_\ell = 0} = \Theta_{(k-1)} (g_1 , \ldots , g_\ell g_{\ell + 1} , \ldots , g_k)$ for $1 \leq \ell \leq k$.  When $t_\ell = 0$, we have
\begin{align*}
\Theta_{(k)} (g_1 , \ldots  , g_k) |_{t_\ell = 0} &= -\sum_{\substack{ i=1 \\ i \neq \ell } }^k dt_i \theta_{g_1 \cdots g_i} \displaybreak[0]\\
&\quad + \sum_{\substack{1 \leq i < j \leq k \\ i,j \neq \ell} } d\alpha (g_1 \cdots g_i, g_{i+1} \cdots g_j) (t_i dt_j - t_j dt_i) \displaybreak[0]\\
&\qquad \qquad + 2 \alpha (g_1 \cdots g_i , g_{i+1} \cdots g_j) dt_i dt_j \displaybreak[1]\\
&= \Theta_{(k-1)} (g_1 , \ldots , g_\ell g_{\ell + 1} , \ldots , g_k)
\end{align*}
by substituting $t_i$ for $t_{i+1}$ whenever $i > \ell$.  Now on the face $t_1 + \cdots + t_k = 1$ we have $dt_1 + \cdots + dt_k = 0$ so that by replacing $t_1$
\begin{align*}
&\Theta_{(k)} (g_1 , \ldots  , g_k) |_{t_1 + \cdots + t_k = 1} = \displaybreak[1]\\
&= -\sum_{i=2}^k dt_i \theta_{g_1 \cdots g_i} + \sum_{i=2}^k dt_i \theta_{g_1} \displaybreak[0]\\
&\quad + \sum_{2 \leq i < j \leq k} d\alpha (g_1 \cdots g_i, g_{i+1} \cdots g_j) (t_i dt_j - t_j dt_i) \displaybreak[0]\\
&\quad + \sum_{i=2}^k d \alpha (g_1 , g_2 \cdots g_i) ( (1- t_2 - \cdots - t_k) dt_i - t_i (- dt_2 - \cdots - dt_k))\displaybreak[0] \\
&\quad + \sum_{2 \leq i < j \leq k} 2 \alpha (g_1 \cdots g_i , g_{i+1} \cdots g_j) dt_i dt_j \displaybreak[0]\\
&\quad + \sum_{i=2}^k 2 \alpha (g_1 , g_2 \cdots g_i) (-dt_2 -\cdots -dt_k) dt_i \displaybreak[2]
\intertext{By Proposition \ref{prop:gerbe_connection2_group} applied in the first line and Corollary \ref{cor:discrepancy_identity2} as well as $t_idt_i -t_i dt_i = 0$ in the third line and $dt_i dt_j = - dt_j dt_i$ applied in the other lines}
&= \sum_{i=2}^k dt_i ( d \alpha (g_1 , g_2 \cdots g_i) - \theta_{g_2 \cdots g_i}^{g_1}) \displaybreak[0]\\
&\quad + \sum_{2 \leq i < j \leq k} (d \alpha (g_1 , g_2 \cdots g_i) - d\alpha (g_1 , g_2 \cdots g_j) + d \alpha (g_1 \cdots g_i , g_{i+1} \cdots g_j) ) (t_i dt_j - t_j dt_i) \displaybreak[0]\\
&\quad + \sum_{i=2}^k d\alpha (g_1 , g_2 \cdots g_i) dt_i\displaybreak[0] \\
&\quad + \sum_{2 \leq i < j \leq k} 2 (\alpha (g_1 , g_2 \cdots g_i ) - \alpha (g_1 , g_2 \cdots g_j) + \alpha ( g_1 \cdots g_i , g_{i+1} \cdots g_j)) dt_i dt_j \\
&= - \sum_{i=2}^k dt_i \theta_{g_2 \cdots g_i}^{g_1} \displaybreak[2]\\
&\quad + \sum_{2 \leq i < j \leq k} d\alpha (g_2 \cdots g_i , g_{i+1} \cdots g_j)^{g_1} (t_i dt_j -t_j dt_i) \displaybreak[0]\\
&\quad + \sum_{2 \leq i < j \leq k} 2\alpha (g_2 \cdots g_i , g_{i+1} \cdots g_j)^{g_1} dt_i dt_j \displaybreak[1]\\
&= \Theta_{(k-1)} (g_2 ,\ldots , g_k)^{g_1}
\end{align*}
by a cancellation and application of Corollary \ref{cor:discrepancy_identity2} again.
\end{proof}
\end{thm}

\begin{defn}
Given a gerbe $(L , \mu)$ on $M \rtimes \Gamma$ with connection $\nabla$ on $L \rightarrow M \rtimes \Gamma$ and simplicial $2$-form $\{\vartheta_{(k)} \}$ associated to $\nabla$, the simplicial $3$-form $-\Theta = -\{ \nabla^k \vartheta_{(k)} \}$ is the \textit{simplicial Dixmier-Douady form} associated to $\nabla$.

\begin{thm}
Let $(L ,\mu)$ be a gerbe on $M \rtimes \Gamma$ and $\nabla$ a connection on $L \rightarrow M \rtimes \Gamma$ with discrepancy $\alpha$ and simplicial 2-form $\vartheta$. Then the simplicial Dixmier-Douady form $-\Theta \in \Omega^3 (M_\bullet)$ associated to $\nabla$ is a representative of the Dixmier-Douady class of $\nabla$ on $(L, \mu)$.
\begin{proof}

We will apply the quasi-isomorphism of integration along simplices $\mathcal{I}_\Delta$ described in Theorem \ref{thm:simplex_integration} to see $\mathcal{I}_\Delta (\Theta) = (-\alpha , \theta) \in \Omega^1 ( M_{(2)}) \oplus \Omega^2 ( M_{(1)} )$.  For any $g \in \Gamma$,
\begin{equation*}
\Theta_{(1)} (g) = - dt_1 \theta_g
\end{equation*}
so that restricted to $M \times \{g\}$
\begin{equation*}
\mathcal{I}_\Delta : \Theta_{(1)} (g) \mapsto \int_{\Delta^1} \Theta_{(1)} (g) = \int_0^1 - \theta_g dt_1 = -\theta_g.
\end{equation*}
For any $g_1 , g_2 \in \Gamma$,
\begin{equation*}
\Theta_{(2)} (g_1 , g_2) = -dt_1 \theta_{g_1} - dt_2 \theta_{g_2} + d\alpha (g_1 , g_2) (t_1 dt_2 - t_2 dt_1) + 2 \alpha (g_1 , g_2) dt_1 dt_2
\end{equation*}
so that restricted to $M \times \{g_1\} \times \{g_2\}$ we have
\begin{align*}
\mathcal{I}_\Delta : \Theta_{(2)} (g_1 , g_2) &\mapsto \int_{\Delta^2} \Theta_{(2)} (g_1 , g_2)\\
&= \int_{\Delta^2} -dt_1 \theta_{g_1} - dt_2 \theta_{g_2} + d\alpha (g_1 , g_2) (t_1 dt_2 - t_2 dt_1) + 2 \alpha (g_1 , g_2) dt_1 dt_2 \\
&= \int_{\Delta^2} 2 \alpha (g_1 , g_2) dt_1 dt_2 \\
&= \int_0^1 \int_0^{1-t_1} 2 \alpha (g_1 , g_2) dt_1 dt_2 \\
&= \alpha (g_1 , g_2).
\end{align*}
Furthermore, $\Theta_{(k)}^{r,s} = 0$ for $s \geq 3$, i.e.\ $\Theta_{(k)}$ has at most two $dt$'s in any term, so that 
\begin{equation*}
\mathcal{I}_\Delta (\Theta_{(k)}) = \int_{\Delta^k} \Theta_{(k)}= 0
\end{equation*}
whenever $k \geq 3$. Hence $\mathcal{I}_\Delta : -\Theta \mapsto -\alpha + \theta$ so that $-\Theta$ is a representative of the Dixmier-Douady class defined by $(-\alpha, \theta)$ in the bicomplex of compatible forms.
\end{proof}
\end{thm}

We may also introduce a rescaled version of the simplicial Dixmier-Douady form denoted $\Theta_u$ in the complex $\Omega^*_- (M_\bullet)$ as in \eqref{eq:Theta_rescaled}.  With the rescaling $\nabla^k_u =  (\nabla^k)^{1,0} + u^{-1}d_\Delta$ we have $\nabla^k_u (\vartheta_u)_{(k)} = u^{-1} (\Theta_u)_{(k)}$ as
\begin{align*}
( (\nabla^k)^{1,0} + u^{-1} d_\Delta) ( u (\vartheta_{(k)})^{2,0} + (\theta_{(k)})^{1,1}) &= u ((\nabla^k)^{1,0} (\vartheta_{(k)})^{2,0})^{3,0} \\
&\quad +  (( \nabla^k)^{1,0} (\vartheta_{(k)})^{1,1} + d_\Delta (\vartheta_{(k)})^{2,0} )^{2,1} \\
&\quad + u^{-1} (d_\Delta (\vartheta_{(k)}))^{1,2} \\
&=u (\Theta_{(k)})^{3,0} +  (\Theta_{(k)})^{2,1} + u^{-1} (\Theta_{(k)})^{1,2}
\end{align*}
Also note that under this rescaling, $u (\nabla^k_u)^2 = [ (\vartheta_u)_{(k)} , \cdot]$.
\end{defn}

%

\chapter{Cocycles on the twisted convolution algebra}
\label{chap:4}

Now we will construct a morphism from the twisted simplicial complex of a gerbe $(L,\mu)$ on $M \rtimes \Gamma$ to the periodic cyclic complex of the twisted convolution algebra $C^\infty_c (M \rtimes \Gamma, L)$.  We will continue to denote $M_\bullet = (M \rtimes \Gamma)_\bullet$.

\section{Twisted Convolution Algebra}

\begin{defn}
For a gerbe $(L, \mu)$ on $\G$, the \textit{twisted convolution algebra} $C^\infty_c (\G, L)$ is defined as follows. The convolution product on $C^\infty_c (\G , L)$ is defined by the (finite) sum
\begin{equation}
 (f_1 * f_2 ) (g) = \sum_{\{(g_1, g_2) \in \G^{(2)} | g_1 g_2 = g \} } f_1 (g_1) \cdot f_2 (g_2)
 \end{equation}
for sections $f_1 , f_2 \in C^\infty_c (\G , L)$, where $f_1 (g_1) \cdot f_2 (g_2)$ is computed using the product $\mu_{(g_1,g_2)} : L|_{g_1} \otimes L|_{g_2} \xrightarrow{\sim} L_{g}$.  This gives $C^\infty_c (\G, L)$ the structure of an algebra.

In particular, for $\G = M \rtimes \Gamma$, the convolution product on $C^\infty_c (M \rtimes \Gamma, L)$ may be written
\begin{align}
(f_1 * f_2) (x,g) &= \sum_{g_1 g_2 = g}  \mu_{g_1 , g_2} (f_1|_{M_{g_1}} \otimes (f_2|_{M_{g_2}})^{g_1} ) (x , g) \\
&= \sum_{g_1 g_2 = g} f_1 (x,g_1) \cdot f_2 ( xg_1 , g_2)
\end{align}
for sections $f_1 , f_2 \in C^\infty_c (M \times \Gamma , L)$, where $f_1 (x,g_1) \cdot f_2 (xg_1 , g_2)$ is computed using the product $\mu_{g_1 , g_2} : L_{g_1} \otimes (L_{g_2})^{g_1} \xrightarrow{\sim} L_g$.
\end{defn}


\section{JLO morphism}

Following \cite{jlo88}, \cite{gorokhovsky99}, \cite{mathai_stevenson06}, and \cite{tu_xu06} we now construct a JLO-type morphism from the twisted simplicial complex of $M \rtimes \Gamma$.

For the discrete group $\Gamma$ and a $\Gamma$-module $K$, let $C^n (\Gamma , K)$ denote the space of degree $n$ $\Gamma$-cochains, i.e.\ the space of maps $\Gamma^n \rightarrow K$.  For any $f \in C^n (\Gamma, K)$ define $\delta_\Gamma : C^n (\Gamma , K) \rightarrow C^{n+1} (\Gamma , K)$ to be the group coboundary $\delta_\Gamma = \sum_{i=0}^n (-1)^i {\delta_\Gamma}_i$ where
\begin{equation*}
({\delta_\Gamma}_i f) (g_1 , \ldots , g_{n+1}) = 
\begin{cases}
f(g_2 , \ldots , g_{n+1})^{g_1} & \text{for $i=0$}\\
 (-1)^i f( g_1 , \ldots , g_{i-1} , g_i g_{i+1} , g_{i+2} , \ldots , g_{n+1}) & \text{for $1 \leq i \leq n$}\\
 (-1)^{n+1} f(g_1 , \ldots , g_n) & \text{for $i=n+1$},
\end{cases}
\end{equation*}
where the superscript ${}^{g_1}$ denotes the action of $g_1 \in \Gamma$ on $K$.  The cohomology of the complex $(C^\bullet (\Gamma, K) , {\delta_\Gamma})$ is the group cohomology of $\Gamma$ with coefficients in $K$.  In particular we will consider the complex 
\begin{equation}
(C^k (\Gamma, \overline{C}^n( C^\infty_c (M , \End \E) ) [u^{-1} , u] ) , (b+uB) + (-1)^n {\delta_\Gamma})
\end{equation}
of $\Gamma$-cochains with values in the periodic cyclic complex of $C^\infty_c (M , \End \E)$.  By the operators $b$ and $uB$ on such $\Gamma$-cochains we mean the operators defined by 
\begin{equation}
( b f) (g_1 , \ldots , g_n) = b (f(g_1 , \ldots , g_n))
\end{equation}
 and 
\begin{equation}
(uB f) (g_1 , \ldots , g_n) = uB (f (g_1 , \ldots , g_n))
\end{equation}
for any $f \in C^n (\Gamma, \overline{C}^\bullet ( C^\infty_c (M , \End \E) ) [u^{-1} , u] )$. Let us also denote ${\delta_\Gamma}' = (-1)^n {\delta_\Gamma}$.

\begin{thm}
Given a gerbe $(L , \mu)$ on $M \rtimes \Gamma$ with connection $\nabla$, associated simplicial $2$-form $\vartheta$ and associated simplicial Dixmier-Douady form $\Theta$, and direct sum bundle $\E$, there is a morphism 
\begin{multline}
\tau_\nabla : (\Omega^*_- (M_\bullet) [u] , u\tilde{d}_{dR} + d_\Delta - \Theta_u \wedge \cdot) \\ \to (C^\bullet (\Gamma , \overline{C}^\bullet ( C^\infty_c (M , \End \E) )  [u^{-1} , u] ) , b+uB + {\delta_\Gamma}')
\end{multline}
from the $\Theta$-twisted simplicial complex of $M_\bullet$ to the complex of $\Gamma$-cochains taking values in the periodic cyclic complex of $C^\infty_c (M , \End \E)$.
The morphism $\tau_\nabla$ is defined by the JLO-type formula
\begin{multline}
\tau_\nabla (  \omega ) (\tilde{a_0} , a_1 , \ldots , a_n) = \sum_k \int_M \int_{\Delta^k}  \\ \omega_{(k)} \wedge \left( \int_{\Delta^n} \tr (\tilde{a_0} e^{-\sigma_0 \left( \vartheta_u \right)_{(k)}} \nabla_u^k (a_1) e^{-\sigma_1  \left( \vartheta_u \right)_{(k)}} \cdots \nabla_u^k (a_n) e^{-\sigma_n  \left( \vartheta_u \right)_{(k)}} ) d \sigma_1 \cdots d \sigma_n \right) 
\end{multline}
where $\omega = \{ \omega_{(k)} \} \in \Omega^*_- (M_\bullet)$, $a_1 , \ldots , a_n \in C^\infty_c (M , \End \E)$, $\tilde{a_0} \in \widetilde{C^\infty_c (M , \End \E)}$, and $\sigma_0 , \ldots , \sigma_n$ are barycentric coordinates on $\Delta^n$.
\begin{proof}
Let us first discuss the reason the sum is finite.  We must check in particular that the degrees of the forms being integrated on the simplices $\Delta^k$ being summed over are bounded.  Given $\omega \in \Omega^\ell_- (M_\bullet)$, and $a_0 , \ldots , a_n \in \overline{C}^\bullet ( C^\infty_c (M , \End \E) )$, consider
\begin{equation}
\omega_{(k)} \wedge \left( \int_{\Delta^n} \tr (\tilde{a_0} e^{-\sigma_0 \left( \vartheta_u \right)_{(k)}} \nabla_u^k (a_1) e^{-\sigma_1  \left( \vartheta_u \right)_{(k)}} \cdots \nabla_u^k (a_n) e^{-\sigma_n  \left( \vartheta_u \right)_{(k)}} ) d \sigma_1 \cdots d \sigma_n \right) 
\end{equation}
The maximum number of $dt_i$'s due to $\omega$ is $\ell$, where $t_i$ are coordinates on the simplex $\Delta^k$, as the number of $dt_i$'s cannot exceed the degree of $\omega$.  Within the trace, the terms $\nabla^k_u (a_i)$ do not contribute any $dt_i$'s because $a_j$ is extended trivially to $\Delta^n$ and $(\nabla^k)^{0,1} = d_\Delta$.  Since $\vartheta_{(k)} = \vartheta_{(k)}^{2,0} + \vartheta_{(k)}^{1,1}$ and in particular $\vartheta_{(k)}^{0,2} = 0$, the contribution to the degree of the form from all the $(\vartheta_u)_{(k)}$ appearing must be at most $\dim M$, i.e.\ the number of $dx_i$'s (where $x_i$'s are coordinates on the manifold $M$) is greater than or equal to the number of $dt_j$'s. Hence the simplex degree of the form above is at most $\ell  + \dim M$ so the sum is finite.

Without loss of generality we will compute $\tau_\nabla ( (u\tilde{d}_{dR} + d_{\Delta}) \omega) (a_0 , a_1 , \ldots , a_n)$ as $\tilde{a_0} e^{-\sigma_0 \left( \vartheta_u \right)_{(k)}} \nabla_u^k (a_1) \cdots  e^{-\sigma_n  \left( \vartheta_u \right)_{(k)}} $ in the integrand above means
\begin{equation*}
a_0 e^{-\sigma_0 \left( \vartheta_u \right)_{(k)}} \nabla_u^k (a_1)  \cdots  e^{-\sigma_n  \left( \vartheta_u \right)_{(k)}} + \lambda e^{-\sigma_0 \left( \vartheta_u \right)_{(k)}} \nabla_u^k (a_1)  \cdots e^{-\sigma_n  \left( \vartheta_u \right)_{(k)}} 
\end{equation*}
for $\tilde{a}_0 = (a_0 , \lambda)$.

Following \cite{mathai_stevenson06}, we will compute the effect of $(u\tilde{d}_{dR} + d_{\Delta})$.  This gives
\begin{align}
&\tau_\nabla ( (u\tilde{d}_{dR} + d_{\Delta}) \omega) (a_0 , a_1 , \ldots , a_n)  \notag \\
&= \sum_k  \int_M \int_{\Delta^k}  (u\tilde{d}_{dR} + d_{\Delta}) \omega_{(k)} \wedge \bigg( \int_{\Delta^n} \tr (a_0 e^{-\sigma_0 \left( \vartheta_u \right)_{(k)}} \nabla_u^k (a_1)  \cdots \displaybreak[0]\notag\\
&\qquad \qquad \cdots e^{-\sigma_n  \left( \vartheta_u \right)_{(k)}} ) d \sigma_1 \cdots d \sigma_n \bigg)  \displaybreak[3] \notag\\
&=- \sum_k  \int_M \int_{\Delta^k}  \omega_{(k)} \wedge \bigg( \int_{\Delta^n} (ud_{dR} + d_{\Delta}) \tr (a_0 e^{-\sigma_0 \left( \vartheta_u \right)_{(k)}} \nabla_u^k (a_1) \cdots \displaybreak[0] \notag\\
&\qquad \qquad \cdots e^{-\sigma_n  \left( \vartheta_u \right)_{(k)}} ) d \sigma_1 \cdots d \sigma_n \bigg)  \displaybreak[2]\notag\\
&\quad + (-1)^{n} \int_M \int_{\partial(\Delta^k)}  \omega_{(k)} \wedge \bigg( \int_{\Delta^n} \tr (a_0 e^{-\sigma_0 \left( \vartheta_u \right)_{(k)}} \nabla_u^k (a_1)  \cdots \displaybreak[0] \notag \\
&\qquad \qquad \cdots  e^{-\sigma_n  \left( \vartheta_u \right)_{(k)}} ) d \sigma_1 \cdots d \sigma_n \bigg)  \notag\\
&= -\sum_k  \int_M \int_{\Delta^k}  \omega_{(k)} \wedge \left( \int_{\Delta^n} u \tr (\nabla^k_u (a_0 e^{-\sigma_0 \left( \vartheta_u \right)_{(k)}}   \cdots e^{-\sigma_n  \left( \vartheta_u \right)_{(k)}} ) ) d \sigma_1 \cdots d \sigma_n \right)  \displaybreak[2] \label{eq:JLO_formula_1}\\
&\quad + (-1)^{n } \int_M \int_{\partial(\Delta^k)}  \omega_{(k)} \wedge \left( \int_{\Delta^n} \tr (a_0 e^{-\sigma_0 \left( \vartheta_u \right)_{(k)}}   \cdots  e^{-\sigma_n  \left( \vartheta_u \right)_{(k)}} ) d \sigma_1 \cdots d \sigma_n \right) , \label{eq:JLO_formula_2}
\end{align}
where the second equality follows from Lemma \ref{lem:Stokes_wedge_twisted_simplicial_complex} and the last from Lemma \ref{lem:trace_connection}. By commuting forms inside the trace within the first term \eqref{eq:JLO_formula_1} we have
\begin{align*}
&u \tr \left( \nabla^k_u (a_0 e^{-\sigma_0 \left( \vartheta_u \right)_{(k)}} \nabla_u^k (a_1) \cdots e^{-\sigma_n  \left( \vartheta_u \right)_{(k)}} ) \right) \\
&=  u \tr \left( \nabla^k_u  (a_0) e^{-\sigma_0 \left( \vartheta_u \right)_{(k)}} \nabla_u^k (a_1) \cdots e^{-\sigma_n  \left( \vartheta_u \right)_{(k)}} \right) \\
&\quad + \sum_{i=0}^{n-1} (-1)^i \left( u \tr \left(a_0 e^{-\sigma_0 \left( \vartheta_u \right)_{(k)}} \cdots \nabla^k_u (a_i) \nabla^k_u \left(e^{-\sigma_i \left( \vartheta_u \right)_{(k)}} \right) \nabla^k_u (a_{i+1}) \cdots e^{-\sigma_n  \left( \vartheta_u \right)_{(k)}}  \right) \right. \\
&\qquad \qquad + \left. u \tr \left(  a_0 e^{-\sigma_0 \left( \vartheta_u \right)_{(k)}} \cdots \nabla^k_u (a_i) e^{-\sigma_i \left( \vartheta_u \right)_{(k)}} (\nabla^k_u)^2 (a_{i+1}) \cdots e^{-\sigma_n  \left( \vartheta_u \right)_{(k)}} \right) \right) \\
&\quad + (-1)^n u \tr \left( a_0 e^{-\sigma_0 \left( \vartheta_u \right)_{(k)}} \nabla_u^k (a_1) \cdots \nabla^k_u (a_n) \nabla^k_u (e^{-\sigma_n  \left( \vartheta_u \right)_{(k)}} ) \right) \displaybreak[2]\\
&=  u\tr \left( \nabla^k_u  (a_0) e^{-\sigma_0 \left( \vartheta_u \right)_{(k)}} \nabla_u^k (a_1) \cdots e^{-\sigma_n  \left( \vartheta_u \right)_{(k)}} \right) \\
&\quad + \sum_{i=0}^{n}  (-1)^{i+i(n+1-i)} u \tr \bigg( \nabla^k_u \left(e^{-\sigma_i \left( \vartheta_u \right)_{(k)}} \right) \nabla^k_u (a_{i+1}) \cdots \\
&\qquad \qquad \cdots e^{-\sigma_n  \left( \vartheta_u \right)_{(k)}} a_0 e^{-\sigma_0 \left( \vartheta_u \right)_{(k)}} \cdots \nabla^k_u (a_i)  \bigg) \\
&\quad + \sum_{i=0}^{n-1}  (-1)^i  \tr \left(  a_0 e^{-\sigma_0 \left( \vartheta_u \right)_{(k)}} \cdots \nabla^k_u (a_i) e^{-\sigma_i \left( \vartheta_u \right)_{(k)}} [(\vartheta_u)_{(k)} , a_{i+1}] \cdots e^{-\sigma_n  \left( \vartheta_u \right)_{(k)}} \right) \displaybreak[2]\\
&=  u\tr \left( \nabla^k_u  (a_0) e^{-\sigma_0 \left( \vartheta_u \right)_{(k)}} \nabla_u^k (a_1) \cdots e^{-\sigma_n  \left( \vartheta_u \right)_{(k)}} \right) \\
&\quad + \sum_{i=0}^{n}  (-1)^{i+i(n+1-i)+i(n-i)}  (-\sigma_i (\Theta_u)_{(k)}) \wedge \tr \left(  a_0 e^{-\sigma_0 \left( \vartheta_u \right)_{(k)}} \nabla_u^k (a_1) \cdots e^{-\sigma_n  \left( \vartheta_u \right)_{(k)}}\right) \\
&\quad + \sum_{i=0}^{n-1}  (-1)^i  \tr \left(  a_0 e^{-\sigma_0 \left( \vartheta_u \right)_{(k)}} \cdots \nabla^k_u (a_i) e^{-\sigma_i \left( \vartheta_u \right)_{(k)}} [(\vartheta_u)_{(k)} , a_{i+1}] \cdots e^{-\sigma_n  \left( \vartheta_u \right)_{(k)}} \right)\displaybreak[2] \\
&=  u\tr \left( \nabla^k_u  (a_0) e^{-\sigma_0 \left( \vartheta_u \right)_{(k)}} \nabla_u^k (a_1) \cdots e^{-\sigma_n  \left( \vartheta_u \right)_{(k)}} \right) \\
&\quad - (\Theta_u)_{(k)} \wedge \tr \left(  a_0 e^{-\sigma_0 \left( \vartheta_u \right)_{(k)}} \nabla_u^k (a_1) \cdots e^{-\sigma_n  \left( \vartheta_u \right)_{(k)}}\right) \displaybreak[2]\\
&\quad + \sum_{i=0}^{n-1}  (-1)^i  \tr \left(  a_0 e^{-\sigma_0 \left( \vartheta_u \right)_{(k)}} \cdots \nabla^k_u (a_i) e^{-\sigma_i \left( \vartheta_u \right)_{(k)}} [(\vartheta_u)_{(k)} , a_{i+1}] \cdots e^{-\sigma_n  \left( \vartheta_u \right)_{(k)}} \right) 
\end{align*}
where we have used $u (\nabla_u^k)^2 (a) = [(\vartheta_u)_{(k)} , a]$, $u \nabla^k_u \left( \left( \vartheta_u \right)_{(k)} \right) =\left( \Theta_u \right)_{(k)}$ and
\begin{equation}
u \nabla^k_u \left( e^{-\sigma_i \left( \vartheta_u \right)_{(k)}} \right) = - \left(d \sigma_i \left( \vartheta_u \right)_{(k)} + \sigma_i  \left( \Theta_u \right)_{(k)} \right) e^{-\sigma_i \left( \vartheta_u \right)_{(k)}}
\end{equation}
as well as $\sum_{i=0}^n d \sigma_i = 0$ and $\sum_{i=0}^n \sigma_i = 1$ on $\Delta^n$.  Therefore by adding $\tau_\nabla ( - \Theta_u \wedge \omega ) $ to $\tau_\nabla ( (u \tilde{d}_{dR} + d_\Delta))$ we obtain
\begin{align}
& - \sum_k  \int_M \int_{\Delta^k}  \omega_{(k)}  \wedge \left( \int_{\Delta^n} u \tr (\nabla^k_u (a_0) e^{-\sigma_0 \left( \vartheta_u \right)_{(k)}}  \cdots e^{-\sigma_n  \left( \vartheta_u \right)_{(k)}} )  d \sigma_1 \cdots d \sigma_n \right) \label{eq:JLO_formula_B1} \\
&\quad - \int_M \int_{\Delta^k}  \omega_{(k)}  \wedge \bigg( \int_{\Delta^n} \sum_{i=0}^{n-1}  (-1)^i  \tr \Big(  a_0 e^{-\sigma_0 \left( \vartheta_u \right)_{(k)}} \cdots \notag\\
&\qquad \qquad \cdots \nabla^k_u (a_i) e^{-\sigma_i \left( \vartheta_u \right)_{(k)}} [(\vartheta_u)_{(k)} , a_{i+1}] \cdots e^{-\sigma_n  \left( \vartheta_u \right)_{(k)}} \Big) \bigg) \label{eq:JLO_formula_b1} \\
&\quad + (-1)^n \int_M \int_{\partial(\Delta^k)} \omega_{(k)} \wedge \bigg( \int_{\Delta^n} \tr \big(a_0 e^{-\sigma_0 \left( \vartheta_u \right)_{(k)}} \nabla_u^k (a_1)  \cdots \notag \\
&\qquad \qquad e^{-\sigma_n  \left( \vartheta_u \right)_{(k)}} \big) d \sigma_1 \cdots d \sigma_n \bigg) \label{eq:JLO_formula_G1} 
\end{align}
Next we compute the effect of the algebraic operators.  The first term above, \eqref{eq:JLO_formula_B1}, is equivalent to $-(uB) \tau_\nabla (\omega) (a_0 , \ldots , a_n)$ as
\begin{align*}
&(uB) \tau_\nabla (\omega) (a_0 , \ldots , a_n) \displaybreak[2]\\
&= \sum_{i=0}^n (-1)^{ni} u \tau_\nabla (\omega) (1, a_i , \ldots , a_n , a_0 , \ldots a_{i-1}) \displaybreak[2] \\
&= \sum_{i=0}^n (-1)^{ni} u \sum_k \int_M \int_{\Delta^k} \omega_{(k)} \wedge \bigg( \int_{\Delta^{n+1}} \tr \Big( e^{-\sigma_0 \left( \vartheta_u \right)_{(k)}} \nabla_u^k (a_i) e^{-\sigma_1  \left( \vartheta_u \right)_{(k)}} \cdots  \displaybreak[0]\\ 
&\quad \cdots \nabla_u^k (a_n) e^{-\sigma_{n-i+1}  \left( \vartheta_u \right)_{(k)}} \nabla_u^k (a_0) e^{-\sigma_{n-i+2}  \left( \vartheta_u \right)_{(k)}} \cdots \displaybreak[0]\\
&\quad \cdots \nabla_u^k (a_{i-1}) e^{-\sigma_{n+1}  \left( \vartheta_u \right)_{(k)}} \Big) d \sigma_1 \cdots d \sigma_{n+1} \bigg) \displaybreak[2]\\
&= \sum_{i=0}^n (-1)^{ni} u \sum_k \int_M \int_{\Delta^k} \omega_{(k)} \wedge \bigg( \int_{\Delta^{n+1}} \tr \Big(  \nabla_u^k (a_i) e^{-\sigma_1  \left( \vartheta_u \right)_{(k)}} \cdots  \\ 
&\quad \cdots  \nabla_u^k (a_n) e^{-\sigma_{n-i+1}  \left( \vartheta_u \right)_{(k)}} \nabla_u^k (a_0) e^{-\sigma_{n-i+2}  \left( \vartheta_u \right)_{(k)}} \cdots \\
&\quad \cdots \nabla_u^k (a_{i-1}) e^{-(\sigma_0+\sigma_{n+1})  \left( \vartheta_u \right)_{(k)}} \Big) d \sigma_1 \cdots d \sigma_{n+1} \bigg) \displaybreak[2]\\
&= \sum_{i=0}^n (-1)^{ni+(n-i+1)i} u \sum_k \int_M \int_{\Delta^k} \omega_{(k)} \wedge \bigg( \int_{\Delta^{n+1}} \tr \Big( \nabla_u^k (a_0) e^{-\sigma_{n-i+2}  \left( \vartheta_u \right)_{(k)}} \cdots  \\ 
&\quad  \cdots \nabla_u^k (a_{i-1}) e^{-(\sigma_0+\sigma_{n+1})  \left( \vartheta_u \right)_{(k)}} \nabla_u^k (a_i) e^{-\sigma_1  \left( \vartheta_u \right)_{(k)}} \cdots \\
&\quad \cdots \nabla_u^k (a_n) e^{-\sigma_{n-i+1}  \left( \vartheta_u \right)_{(k)}}  \Big) d \sigma_1 \cdots d \sigma_{n+1} \bigg) \displaybreak[2]\\
&= \sum_{i=0}^n  u \sum_k \int_M \int_{\Delta^k} \omega_{(k)} \wedge \bigg( \int_{\Delta^{n}} \sigma_i \tr \Big( \nabla_u^k (a_0) e^{-\sigma_{0}  \left( \vartheta_u \right)_{(k)}} \nabla^k_u (a_1) \cdots \\
&\quad \cdots \nabla_u^k (a_n) e^{-\sigma_{n}  \left( \vartheta_u \right)_{(k)}}  \Big) d \sigma_1 \cdots d \sigma_{n} \bigg) \displaybreak[2]\\
&= u \sum_k \int_M \int_{\Delta^k} \omega_{(k)} \wedge \bigg( \int_{\Delta^{n}} \tr \big( \nabla_u^k (a_0) e^{-\sigma_{0}  \left( \vartheta_u \right)_{(k)}} \nabla^k_u (a_1) \cdots \displaybreak[1]\\
&\quad \cdots \nabla_u^k (a_n) e^{-\sigma_{n}  \left( \vartheta_u \right)_{(k)}}  \Big) d \sigma_1 \cdots d \sigma_{n} \bigg)
\end{align*}
where we have changed variables, performed a partial integration, and then used that $\sum \sigma_i = 1$.

The second term in $\tau_\nabla ( (u\tilde{d}_{dR} + d_\Delta - \Theta_u \wedge \cdot) \omega)$, \eqref{eq:JLO_formula_b1}, is $-b \tau_\nabla (\omega) (a_0 , \ldots , a_n)$ as
\begin{align*}
& b \tau_\nabla (\omega) (a_0 , \ldots , a_n) \\
&= \sum_{i=0}^{n-1} (-1)^i \tau_\nabla (a_0 , \ldots , a_i a_{i+1} , \ldots , a_n) + (-1)^{n} \tau_\nabla (a_n a_0 , a_1 , \ldots , a_{n-1}) \displaybreak[2]\\
&=  \sum_k \int_M \int_{\Delta^k} \omega_{(k)} \wedge  \bigg( \int_{\Delta^{n-1}} \tr \Big(\nabla^k_u (a_0 a_1) e^{-\sigma_0 \left( \vartheta_u \right)_{(k)}} \nabla_u^k (a_2)  \cdots \\
&\quad \cdots \nabla_u^k (a_n) e^{-\sigma_{n-1}  \left( \vartheta_u \right)_{(k)}} \Big) d \sigma_1 \cdots d \sigma_{n-1} \bigg) \displaybreak[0]\\
&\quad + \sum_{i=1}^{n-1} (-1)^i \sum_k \int_M \int_{\Delta^k} \omega_{(k)} \wedge  \Big( \int_{\Delta^{n-1}} \big( \tr (a_0 e^{-\sigma_0 \left( \vartheta_u \right)_{(k)}} \nabla_u^k (a_1) \cdots \\
&\qquad \qquad  \cdots e^{-\sigma_{i-1} \left( \vartheta_u \right)_{(k)}} \nabla^k_u (a_i a_{i+1}) e^{-\sigma_{i} \left( \vartheta_u \right)_{(k)}} \cdots \nabla^k_u (a_n) e^{-\sigma_{n-1}  \left( \vartheta_u \right)_{(k)}} \big) d \sigma_1 \cdots d \sigma_{n-1} \Big) \displaybreak[0]\\
& \quad + (-1)^n  \sum_k \int_M \int_{\Delta^k} \omega_{(k)} \wedge  \Big( \int_{\Delta^{n-1}} \tr \big(\nabla^k_u (a_n a_0) e^{-\sigma_0 \left( \vartheta_u \right)_{(k)}} \nabla_u^k (a_1)  \cdots   \\
& \qquad \qquad \cdots  \nabla_u^k (a_{n-1}) e^{-\sigma_{n-1}  \left( \vartheta_u \right)_{(k)}} \big) d \sigma_1 \cdots d \sigma_{n-1} \Big) \displaybreak[2]\\
\intertext{and because $\nabla^k_u (a_i a_{i+1}) = \nabla^k_u (a_i) a_{i+1} + a_i \nabla^k (a_{i+1})$ we have}
&= \sum_k \int_M \int_{\Delta^k} \omega_{(k)} \wedge  \Big( \int_{\Delta^{n-1}} \tr \big(\nabla^k_u (a_0) a_1 e^{-\sigma_0 \left( \vartheta_u \right)_{(k)}} \nabla_u^k (a_2)  \cdots \\
&\cdots \nabla_u^k (a_n) e^{-\sigma_{n-1}  \left( \vartheta_u \right)_{(k)}} \big) d \sigma_1 \cdots d \sigma_{n-1} \Big) \displaybreak[1]\\
&\quad + \sum_k \int_M \int_{\Delta^k} \omega_{(k)} \wedge  \Big( \int_{\Delta^{n-1}} \tr \big(a_0 \nabla^k_u ( a_1) e^{-\sigma_0 \left( \vartheta_u \right)_{(k)}} \nabla_u^k (a_2)  \cdots \\
&\cdots \nabla_u^k (a_n) e^{-\sigma_{n-1}  \left( \vartheta_u \right)_{(k)}} \big) d \sigma_1 \cdots d \sigma_{n-1} \Big) \displaybreak[1]\\
&\quad + \sum_{i=1}^{n-1} (-1)^i \sum_k \int_M \int_{\Delta^k} \omega_{(k)} \wedge  \Big( \int_{\Delta^{n-1}} \tr \big(a_0 e^{-\sigma_0 \left( \vartheta_u \right)_{(k)}} \nabla_u^k (a_1) \cdots  \\
&\qquad \qquad \cdots e^{-\sigma_{i-1} \left( \vartheta_u \right)_{(k)}} \nabla^k_u (a_i) a_{i+1} e^{-\sigma_{i} \left( \vartheta_u \right)_{(k)}} \cdots \nabla^k_u (a_n) e^{-\sigma_{n-1}  \left( \vartheta_u \right)_{(k)}} \big) d \sigma_1 \cdots d \sigma_{n-1} \Big) \displaybreak[1]\\
&\quad + \sum_{i=1}^{n-1} (-1)^i \sum_k \int_M \int_{\Delta^k} \omega_{(k)} \wedge  \Big( \int_{\Delta^{n-1}} \tr \big(a_0 e^{-\sigma_0 \left( \vartheta_u \right)_{(k)}} \nabla_u^k (a_1) \cdots  \\
&\qquad \qquad \cdots e^{-\sigma_{i-1} \left( \vartheta_u \right)_{(k)}} a_i \nabla^k_u (a_{i+1})  e^{-\sigma_{i} \left( \vartheta_u \right)_{(k)}} \cdots \nabla^k_u (a_n) e^{-\sigma_{n-1}  \left( \vartheta_u \right)_{(k)}} \big) d \sigma_1 \cdots d \sigma_{n-1} \Big) \displaybreak[1]\\
& \quad + (-1)^n  \sum_k \int_M \int_{\Delta^k} \omega_{(k)} \wedge  \Big( \int_{\Delta^{n-1}} \tr \big(\nabla^k_u (a_n) a_0 e^{-\sigma_0 \left( \vartheta_u \right)_{(k)}} \nabla_u^k (a_1)  \cdots  \\
& \qquad \qquad \cdots \nabla_u^k (a_{n-1}) e^{-\sigma_{n-1}  \left( \vartheta_u \right)_{(k)}} \big) d \sigma_1 \cdots d \sigma_{n-1} \Big)\displaybreak[1] \\
& \quad + (-1)^n  \sum_k \int_M \int_{\Delta^k} \omega_{(k)} \wedge  \Big( \int_{\Delta^{n-1}} \tr \big(a_n\nabla^k_u ( a_0) e^{-\sigma_0 \left( \vartheta_u \right)_{(k)}} \nabla_u^k (a_1)  \cdots \\
& \qquad \qquad \cdots  \nabla_u^k (a_{n-1}) e^{-\sigma_{n-1}  \left( \vartheta_u \right)_{(k)}} \big) d \sigma_1 \cdots d \sigma_{n-1} \Big) \displaybreak[2] \\
\intertext{and by combining terms as well as commuting terms in the trace of the last two terms}
&= \sum_{i=1}^{n} (-1)^{i-1} \sum_k \int_M \int_{\Delta^k} \omega_{(k)} \wedge  \Big( \int_{\Delta^{n-1}} \tr \big(a_0 e^{-\sigma_0 \left( \vartheta_u \right)_{(k)}} \nabla_u^k (a_1) \cdots \\
&\qquad \qquad \cdots e^{-\sigma_{i-2} \left( \vartheta_u \right)_{(k)}} \nabla^k_u (a_{i-1}) [a_{i}, e^{-\sigma_{i-1} \left( \vartheta_u \right)_{(k)}} ] \nabla^k_u (a_{i+1}) \cdots \\
&\qquad \qquad \cdots \nabla^k_u (a_n) e^{-\sigma_{n-1}  \left( \vartheta_u \right)_{(k)}} \big) d \sigma_1 \cdots d \sigma_{n-1} \Big) .
\end{align*}
As $[a_i , \cdot]$ is a derivation, we can use the integration formula of \cite{quillen}, Equation 7.2, to obtain
\begin{align*}
[a_i , e^{- \sigma_i  \left( \vartheta_u \right)_{(k)}} ] &= \int_0^1 e^{- (1-\sigma) \sigma_i  \left( \vartheta_u \right)_{(k)}} [a_i , - \sigma_i  \left( \vartheta_u \right)_{(k)} ] e^{- \sigma \sigma_i  \left( \vartheta_u \right)_{(k)}} d \sigma \\
&= \sigma_i \int_0^1 e^{- (1-\sigma) \sigma_i  \left( \vartheta_u \right)_{(k)}} [   \left( \vartheta_u \right)_{(k)} , a_i ] e^{- \sigma \sigma_i  \left( \vartheta_u \right)_{(k)}} d \sigma 
\end{align*}
so that we have
\begin{align*}
&= \sum_{i=1}^{n} (-1)^{i-1} \sum_k \int_M \int_{\Delta^k} \omega_{(k)} \wedge  \left( \int_{\Delta^{n-1}} \int_0^1 \sigma_i \tr (a_0 e^{-\sigma_0 \left( \vartheta_u \right)_{(k)}} \nabla_u^k (a_1) \cdots \right. \\
&\qquad \qquad \left. \cdots e^{-\sigma_{i-2} \left( \vartheta_u \right)_{(k)}} \nabla^k_u (a_{i-1}) e^{- (1-\sigma) \sigma_i  \left( \vartheta_u \right)_{(k)}} [   \left( \vartheta_u \right)_{(k)} , a_i ] e^{- \sigma \sigma_i  \left( \vartheta_u \right)_{(k)}} \nabla^k_u (a_{i+1}) \cdots \right. \\
&\qquad \qquad \left. \cdots \nabla^k_u (a_n) e^{-\sigma_{n-1}  \left( \vartheta_u \right)_{(k)}} ) d \sigma d \sigma_1 \cdots d \sigma_{n-1} \right) .
\end{align*}
Now by the change of variables in the $i$-th term, $(1-\sigma)\sigma_i = \sigma_{i-1}'$ and $\sigma_i' = \sigma_{i-1} \sigma$ (as well as $\sigma_0 = \sigma_0', \ldots , \sigma_{i-2} = \sigma_{i-2}'$ and $\sigma_i = \sigma_{i+1}' , \ldots , \sigma_{n-1} = \sigma_n '$) we have $\sigma_{i-1} = \sigma_{i-1}' + \sigma_i'$, and in fact $\sigma_0 ' + \cdots + \sigma_n ' =1$ and $d\sigma_0' + \cdots + d\sigma_n' = 0$, i.e.\ barycentric coordinates on $\Delta^n$.  So with this change of variables the above expression is equivalent to \eqref{eq:JLO_formula_b1}.

The third term in $\tau_\nabla ( (u\tilde{d}_{dR} + d_\Delta - \Theta_u \wedge \cdot) \omega)$, \eqref{eq:JLO_formula_G1}, equals ${\delta_\Gamma}' \tau_\nabla (\omega) (a_0 , \ldots , a_n)$ because $\{ \vartheta_{(k)} \}$ and $\{ \omega_{(k)} \}$ are simplicial forms.  In particular, if we define Cartesian coordinates on $\Delta^k$ by $t_1 , \ldots , t_k$ then the restriction of $\vartheta_{(k)}$ to the face $t_i = 0$ is a form on $\Delta^{k-1}$ given by
\begin{equation*}
\vartheta_{(k)} (g_1 , \ldots , g_k) |_{t_i=0} = \vartheta_{(k-1)} (g_1 , \ldots , g_i g_{i+1} , \ldots , g_k)
\end{equation*}
for $1 \leq i \leq k-1$ and 
\begin{equation*}
\vartheta_{(k)} (g_1 , \ldots , g_k) |_{t_k=0} = \vartheta_{(k-1)} (g_1 , \ldots , g_{k-1})
\end{equation*}
The restriction to the face $\sum_{i=1}^k t_i = 1$ is a form on $\Delta^{k-1}$ given by
\begin{equation*}
\vartheta_{(k)} (g_1 , \ldots , g_k) |_{\sum t_i=1} = \vartheta_{(k-1)} (g_2 , \ldots , g_k)^{g_1}
\end{equation*}
as shown in Theorem \ref{thm:simplicial_2_form}.  Hence the corresponding identities hold for $\vartheta_u$. We also have the same sort of identities for $\{ \omega_{(k)} \}$ because $\{ \omega_{(k)} \}$ is a simplicial form. In addition $\nabla^k (g_1 , \ldots , g_k)$ satisfies the following identities for $1 \leq i \leq k-1$
\begin{align*}
\nabla^k (g_1 , \ldots , g_k) |_{t_i =0} &= \nabla^\E + t_1 A(g_1) + \cdots + t_{i-1} A(g_1 \cdots g_{i-1}) + t_{i+1} A(g_1 \cdots g_{i+1}) + \cdots \\
&\qquad \qquad  \cdots + t_k A(g_1 \cdots g_k) \\
&= \nabla^{k-1} (g_1 , \ldots , g_i g_{i+1} , \ldots , g_k)
\end{align*}
with a relabeling of variables.  Similarly $\nabla^k (g_1 , \ldots , g_k)|_{t_k=0} = \nabla^{k-1} (g_1 , \ldots , g_{k-1})$. Furthermore,
\begin{align*}
\nabla^k (g_1 , \ldots , g_k) |_{\sum t_i = 1} &= \nabla^\E + A ( g_1) + t_2 (A (g_1 g_2) - A(g_1)) + \cdots \\
&\qquad \cdots + t_k ( A(g_1 \cdots g_k) - A(g_1)) \\
&= \varphi_{g_1}^* (\nabla^{\E^g} \otimes 1 + 1 \otimes \nabla_g) + t_2 (\alpha( g_1 , g_2) + \varphi_{g_1}^* A(g_2)^{g_1}) + \cdots \\
&\qquad \qquad \cdots + t_k (\alpha (g_1 , g_2 \cdots g_k) + \varphi_{g_1}^* A(g_2 \cdots g_k) ) \\
&= \varphi_{g_1}^* (\nabla^{\E^g} + t_2 A(g_2) + \cdots + t_k A(g_2 \cdots g_k))^{g_1}  \\
&\quad + \varphi_{g_1}^* \nabla_g + t_2 \alpha (g_1 , g_2) + \cdots t_k \alpha(g_1 , g_2 \cdots g_k)\\
&= \varphi_{g_1}^* (\nabla^{\E^g} + t_2 A(g_2) + \cdots + t_k A(g_2 \cdots g_k))^{g_1}  \\
&\cong \nabla^{k-1} (g_2 , \ldots , g_k)^{g_1}
\end{align*}
by Lemma \ref{lem:discrepancy_identity} and the definition of $A(g_1)$. Since $\nabla_g + t_2 \alpha (g_1 , g_2) + \cdots t_k \alpha(g_1 , g_2 \cdots g_k)$ is scalar-valued,
\begin{multline*}
[\nabla^{k-1} (g_2 , \ldots , g_k)^{g_1} + \varphi_{g_1}^* \nabla_g + t_2 \alpha (g_1 , g_2) + \cdots t_k \alpha(g_1 , g_2 \cdots g_k) , \eta ] \\ = [\nabla^{k-1} (g_2 , \ldots , g_k)^{g_1} , \eta]
\end{multline*}
for any $\eta \in C^\infty_c (M , \End \E)$.  The corresponding identities hold for $\nabla^k_u$ as well.

Taking orientation into account we may write 
\begin{equation*}
\partial (\Delta^k) = \Delta^k |_{\sum t_i = 1} - \Delta^k |_{t_1 = 0} + \Delta^k |_{t_2 = 0}  \cdots + (-1)^{k+1} \Delta^k |_{t_k = 0}
\end{equation*}
where a negative sign denotes the orientation opposite the orientation induced by $\Delta^k$. With this the third term of $\tau_\nabla ( (u\tilde{d}_{dR} + d_\Delta - \Theta_u \wedge \cdot) \omega) (a_0 , \ldots , a_n)$, \eqref{eq:JLO_formula_G1}, may be written
\begin{align*}
&(-1)^n \sum_k  \int_M \int_{\Delta^k |_{\sum t_i = 1}} \omega_{(k)} \wedge \left( \int_{\Delta^n} \tr (a_0 e^{-\sigma_0 \left( \vartheta_u \right)_{(k)}} \nabla_u^k (a_1)  \cdots  e^{-\sigma_n  \left( \vartheta_u \right)_{(k)}} ) d \sigma_1 \cdots d \sigma_n \right) \\
&+ (-1)^n \sum_{i=1}^k (-1)^i \int_M \int_{\Delta^k |_{t_i = 0}} \omega_{(k)} \wedge \left( \int_{\Delta^n} \tr (a_0 e^{-\sigma_0 \left( \vartheta_u \right)_{(k)}}   \cdots  e^{-\sigma_n  \left( \vartheta_u \right)_{(k)}} ) d \sigma_1 \cdots d \sigma_n \right).
\end{align*}
Upon evaluation on $g_1 , \ldots , g_k$
\begin{align*}
&(-1)^n \sum_k  \int_M \int_{\Delta^{k-1}} \omega_{(k-1)} (g_2 , \ldots , g_k)^{g_1} \wedge \left( \int_{\Delta^n} \tr (a_0 e^{-\sigma_0 \left( \left( \vartheta_u \right)_{(k-1)} (g_2 , \ldots , g_k) \right)^{g_1} } \right. \\
&\qquad \qquad \left. (\nabla_u^{k-1} (g_2 , \ldots , g_k) (a_1))^{g_1} \cdots  e^{-\sigma_n \left( \left( \vartheta_u \right)_{(k-1)} (g_2 , \ldots , g_k) \right)^{g_1}}  ) d \sigma_1 \cdots d \sigma_n \right) \\
&+ (-1)^n \sum_{i=1}^{k-1} (-1)^i \int_M \int_{\Delta^{k-1}} \omega_{(k-1)} (g_1 , \ldots , g_i g_{i+1} , \ldots , g_k) \wedge  \\
&\qquad \qquad  \wedge \left( \int_{\Delta^n} \tr (a_0 e^{-\sigma_0 \left( \vartheta_u \right)_{(k-1)} (g_1 , \ldots , g_i g_{i+1} , \ldots , g_k) } (\nabla_u^{k-1} (g_1 , \ldots , g_i g_{i+1} , \ldots , g_k)) (a_1)  \cdots \right. \\
&\qquad \qquad \cdots \left. e^{-\sigma_n  \left( \vartheta_u \right)_{(k-1)} (g_1 , \ldots , g_i g_{i+1} , \ldots , g_k) } ) d \sigma_1 \cdots d \sigma_n \right) \\
&\quad + (-1)^n (-1)^{k} \int_M \int_{\Delta^{k-1}} \omega_{(k-1)} (g_1 , \ldots , g_{k-1}) \wedge  \\
&\qquad \qquad  \wedge \left( \int_{\Delta^n} \tr (a_0 e^{-\sigma_0 \left( \vartheta_u \right)_{(k-1)} (g_1 , \ldots , g_{k-1}) } (\nabla_u^{k-1} (g_1 , \ldots , g_{k-1})) (a_1)  \cdots \right. \\
&\qquad \qquad \cdots \left. e^{-\sigma_n  \left( \vartheta_u \right)_{(k-1)} (g_1 , \ldots , g_{k-1}) } ) d \sigma_1 \cdots d \sigma_n \right) \\
&= (-1)^n \tau_\nabla (\omega) (g_2 , \ldots , g_k)^{g_1} (a_0 , \ldots , a_n) \\
&\quad + (-1)^n \sum_{i=1}^{k-1} \tau_\nabla (\omega) (g_1 , \ldots , g_i g_{i+1}, \ldots , g_k) (a_0 , \ldots , a_n)\\
&\quad + (-1)^n (-1)^k \tau_\nabla (\omega) (g_1 , \ldots , g_{k-1}) (a_0 , \ldots , a_n) \\
&={\delta_\Gamma}' \tau_\nabla (\omega) (g_1 , \ldots , g_k) (a_0 , \ldots , a_n)
\end{align*}
This proves that $\tau_\nabla \circ (u \tilde{d}_{dR} + d_\Delta - \Theta_u \wedge \cdot) = (-(b+ uB) + {\delta_\Gamma}') \circ \tau_\nabla$.
\end{proof}
\end{thm}

\section{Algebraic morphisms}

Now we will need a different cochain complex computing group cohomology.  For a discrete group $\Gamma$ and a $\Gamma$-module $K$, let $\widetilde{C}^n (\Gamma , K)$ denote the space of \textit{homogenous $n$-cochains of $\Gamma$ with values in $K$}.  This is the space of degree $n+1$ $\Gamma$-cochains, i.e.\ the space of maps $\Gamma^{n+1} \rightarrow K$, that satisfy the condition
\begin{equation*}
g f (g_0 , \ldots , g_n) = f (g g_0 , \ldots , g g_n).
\end{equation*}
For any $f \in \widetilde{C}^n (\Gamma , K)$ define $\widetilde{\delta}_\Gamma : \widetilde{C}^n (\Gamma , K) \to \widetilde{C}^{n+1} (\Gamma , K)$ by 
\begin{equation}
({\widetilde{\delta}_\Gamma} f) (g_0 , \ldots , g_n) = \sum_{i=0}^n (-1)^i f(g_0 , \ldots , \hat{g_i} , \ldots , g_n) .
\end{equation}
There is a bijection from $C^n (\Gamma , K) \to \widetilde{C}^n (\Gamma ,K)$ given by $f \mapsto \widetilde{f}$ where $\widetilde{f}$ is defined by
\begin{equation}
\widetilde{f} (g_0 , \ldots , g_n) = f(g_0^{-1} g_1 , (g_0 g_1)^{-1} g_2 , \ldots , (g_0 \cdots g_{n-1})^{-1} g_n )
\end{equation}
for any $f \in C^n (\Gamma , K)$ and $g_0 , \ldots , g_n \in \Gamma$.  Therefore complex $(\widetilde{C} (\Gamma ,K) , {\widetilde{\delta}_\Gamma})$ also computes group cohomology. 

\begin{prop}

The above bijection induces a map
\begin{multline}
\Psi_0 : (C^\bullet (\Gamma , \overline{C}^\bullet ( C^\infty_c (M , \End \E) )  [u^{-1} , u]) , b+uB + {\delta_\Gamma}') \\ \to (\widetilde{C}^\bullet (\Gamma , \overline{C}^\bullet ( C^\infty_c (M , \End \E) ) [u^{-1} , u] ) , b+uB + {\widetilde{\delta}_\Gamma}') 
\end{multline}
by the formula
\begin{equation}
\Psi_0 (c) (g_0 , \ldots , g_n) = c (g_0^{-1} g_1 , (g_0 g_1)^{-1} g_2 , \ldots , (g_0 \cdots g_{n-1})^{-1} g_n)
\end{equation}
for any $c \in C^\bullet (\Gamma , \overline{C}^\bullet ( C^\infty_c (M , \End \E) ) [u^{-1} , u] )$

\end{prop}

\begin{lem}
There is a contracting homotopy
\begin{multline*}
h : (\widetilde{C}^k (\Gamma , \overline{C}^n ( C^\infty_c (M , \End \E) )  [u^{-1} , u] ) , {\widetilde{\delta}_\Gamma} , (-1)^k (b+uB) ) \\ \to (\widetilde{C}^{k-1} (\Gamma , \overline{C}^n ( C^\infty_c (M , \End \E) )  [u^{-1} , u] ) , {\widetilde{\delta}_\Gamma} , (-1)^k (b+uB) )
\end{multline*}

\begin{proof}

First we will define a map $\gamma: C^\infty_c (M , \End \E)^{\otimes n} \to \C \Gamma$ by
\begin{equation}
\gamma : E_{g_0 , g_1} (a_0) \otimes E_{g_1 , g_2} (a_1) \otimes \ldots \otimes E_{g_n , g_0} (a_n) \mapsto g_0
\end{equation}
where $E_{g_i , g_{i+1}} (a_i), E_{g_n , g_0} (a_n) \in C^\infty_c (M , \End \E)$ for $0 \leq i \leq n-1$.  Note that by multilinearity,
\begin{equation}
\gamma ( E_{g_0 , g_1} (a_0) + E_{g_0' , g_1} (a_0') , E_{g_1 , g_2} (a_1) , \ldots , E_{g_n , g_0} (a_n) + E_{g_n , g_0'} (a_n') ) = g_0 + g_0' .
\end{equation}
Furthermore, the map $\gamma$ is equivariant:
\begin{align*}
\gamma ( g \cdot (  E_{g_0 , g_1} (a_0)  \otimes \ldots \otimes E_{g_n , g_0} (a_n)  ) ) &= \gamma ( E_{gg_0 , gg_1} (a_0^g)  \otimes \ldots \otimes E_{gg_n , gg_0} (a_n^g) ) \\
&= g g_0 \\
&= g \gamma ( E_{g_0 , g_1} (a_0)  \otimes \ldots \otimes E_{g_n , g_0} (a_n) )
\end{align*}

Then we can define a map
\begin{equation}
h : \widetilde{C}^{k+1} (\Gamma , \overline{C}^n ( C^\infty_c (M , \End \E) ) ) \to \widetilde{C}^{k} (\Gamma , \overline{C}^n ( C^\infty_c (M , \End \E) ) )
\end{equation}
by
\begin{equation}
(h \varphi) (g_0 , \ldots , g_k) ( s ) = (-1)^k \varphi (g_0 , \ldots , g_k , \gamma (s)) (s)
\end{equation}
for any $s \in C^\infty_c (M , \End \E)^{\otimes n}$.

The map $h$ is a contracting homotopy for group cochains as we have
\begin{align*}
({\widetilde{\delta}_\Gamma} h) f (g_0 , \ldots , g_k) (s) &= \sum_{i=0}^k (-1)^i (hf) (g_0 , \ldots , \hat{g_i} , \ldots , g_k) (s) \\
&= (-1)^{k-1} \sum_{i=0}^k (-1)^i f (g_0 , \ldots , \hat{g_i} , \ldots , g_k , \gamma (s) ) (s) \displaybreak[2]\\
(h {\widetilde{\delta}_\Gamma} ) f (g_0 , \ldots , g_k) (s) &= (-1)^k ({\widetilde{\delta}_\Gamma} f) (g_0 , \ldots , g_k , \gamma (s) ) (s) \\
&= (-1)^k \sum_{i=0}^k (-1)^i f( g_0 , \ldots , \hat{g_i} , \ldots , g_k , \gamma(s) ) (s)  \\
&\quad+ (-1)^{k+(k+1)} f( g_0 , \ldots , g_k) (s) 
\end{align*}
so that
\begin{equation}
({\widetilde{\delta}_\Gamma} h + h {\widetilde{\delta}_\Gamma})  = 1
\end{equation}
is the identity operator.
\end{proof}
\end{lem}

\begin{thm}
Given a gerbe $(L , \mu)$ on $M \rtimes \Gamma$ with an $L$-twisted vector bundle $\E$ there is a morphism
\begin{multline}
\Psi_1 : (\widetilde{C}^\bullet (\Gamma , \overline{C}^\bullet ( C^\infty_c (M , \End \E) )  [u^{-1} , u] ) , b+uB + {\widetilde{\delta}_\Gamma}') \\ \to  ( \overline{C}^\bullet (C^\infty_c (M , \End \E)) [u^{-1} , u], b+uB)^\Gamma
\end{multline}
from the complex of cochains on $\Gamma$ taking values in the periodic cyclic complex of $C^\infty_c (M , \End \E)$ to the periodic cyclic complex of $\Gamma$-invariant cochains on $C^\infty_c (M , \End \E)$.

\begin{proof}

Now we will view $(\widetilde{C}^\bullet (\Gamma , \overline{C}^\bullet ( C^\infty_c (M , \End \E) )  [u^{-1} , u]) , b+uB + {\widetilde{\delta}_\Gamma}')$ as a bicomplex $(\widetilde{C}^k (\Gamma , \overline{C}^n ( C^\infty_c (M , \End \E) )  [u^{-1} , u] ) , {\widetilde{\delta}_\Gamma} , (-1)^k (b+uB) )$ with horizontal differential $\widetilde{\delta}_\Gamma$ and vertical differential $b+uB$.  More explicitly this is done via a morphism
\begin{multline*}
\Psi_{\text{\textinterrobang}} : (\widetilde{C}^\bullet (\Gamma , \overline{C}^\bullet ( C^\infty_c (M , \End \E) ) [u^{-1} , u] ) , b+uB + {\widetilde{\delta}_\Gamma}') \\ \to (\widetilde{C}^k (\Gamma , \overline{C}^n ( C^\infty_c (M , \End \E) )  [u^{-1} , u] ) , {\widetilde{\delta}_\Gamma} , (-1)^k (b+uB) )
\end{multline*}
which is defined by exchanging the vertical and horizontal differentials and making the appropriate sign adjustments ${\widetilde{\delta}_\Gamma}' \mapsto {\widetilde{\delta}_\Gamma}$ and $(b+ uB) \mapsto (-1)^k (b+uB)$.

Since the map $h$ is a contracting homotopy for the rows of this bicomplex, i.e.\ the rows of the augmented double complex are exact, the cohomology of the bicomplex is the cohomology of the initial column.  As described in \cite{bott_tu} in the case of the \v{C}ech-de Rham complex, an explicit morphism
\begin{multline*}
\Psi_1^k : (\widetilde{C}^k (\Gamma , \overline{C}^n ( C^\infty_c (M , \End \E) )  [u^{-1} , u]) , {\widetilde{\delta}_\Gamma} , (-1)^k (b+uB) ) \\ \to  (\overline{C}^{n+k} (C^\infty_c (M , \End \E))  [u^{-1} , u] , b+uB)^\Gamma
\end{multline*}
is given by the formula
\begin{equation*}
\Psi_1^k (c) = ( -(b+uB)h)^k c - h (-(b+uB) h)^{k-1}  ((b+uB) c) - h (-(b+uB) h)^k ( {\widetilde{\delta}_\Gamma} c ).
\end{equation*}
For cochains in column $k$ denote $D= (-1)^k (b+uB)$ for the sake of brevity.  Then this formula is
\begin{equation*}
\Psi_1^k (c) = ( -Dh)^k c - h (-D h)^{k-1}  (D c) - h (-D h)^k ( {\widetilde{\delta}_\Gamma} c ).
\end{equation*}
Then let $\Psi_1' = \sum_k \Psi_1^k$ and let $\Psi_1 = \Psi_1' \circ \Psi_{\text{\textinterrobang}}$.  We can see immediately that if $c$ is a group $k$-cochain then 
\begin{equation*}
\Psi_1 ( c) \in  C^0 ( \Gamma , \overline{C}^{n+k} (C^\infty_c (M , \End \E))  [u^{-1} , u])
\end{equation*}
because $(Dh)$ maps an element of $ C^i ( \Gamma , \overline{C}^{j} (C^\infty_c (M , \End \E))  [u^{-1} , u] )$ to an element of $ C^{i-1} ( \Gamma , \overline{C}^{j+1} (C^\infty_c (M , \End \E))  [u^{-1} , u])$.

At this point let us recount a lemma from \cite{bott_tu}.
\begin{lem}
For $i \geq 1$
\begin{equation*}
{\widetilde{\delta}_\Gamma} (D h)^i = (D h)^i {\widetilde{\delta}_\Gamma} - (D h)^{i-1} D
\end{equation*}

\begin{proof}[Proof of lemma]
Since ${\widetilde{\delta}_\Gamma}$ anticommutes with $D$ and ${\widetilde{\delta}_\Gamma} h + h {\widetilde{\delta}_\Gamma} = 1$,
\begin{align*}
{\widetilde{\delta}_\Gamma} (D h)(D h)^{i-1} &= -D {\widetilde{\delta}_\Gamma} h (D h)^{i-1} \\
&= - D (1- h{\widetilde{\delta}_\Gamma}) (D h)^{i-1} \\
&= (D h) {\widetilde{\delta}_\Gamma} (D h)^{i-1}
\end{align*}
Therefore we can commute ${\widetilde{\delta}_\Gamma}$ and $(D h)$ until we reach $(D h)^{i-1} {\widetilde{\delta}_\Gamma} (Dh)$.  We obtain the formula by
\begin{align*}
{\widetilde{\delta}_\Gamma} (D h)(D h)^{i-1} &= (D h)^{i-1} {\widetilde{\delta}_\Gamma} (Dh) \\
&= - (Dh)^{i-1} D (1 - h {\widetilde{\delta}_\Gamma}) \\
&= - (Dh)^{i-1} D + (Dh)^i {\widetilde{\delta}_\Gamma}
\end{align*}
\end{proof}
\end{lem}

The first step in proving $\Psi_1$ is a morphism is to check that $\Psi_1 (c)$ does in fact land in the complex of $\Gamma$-invariant cyclic cochains $\overline{C}^{n+k} (C^\infty_c (M , \End \E)  [u^{-1} , u] , b+uB)^\Gamma$.  To check this we compute ${\widetilde{\delta}_\Gamma} \Psi_1  (c) $. By the above lemma, for $k \geq 2$,
\begin{align*}
{\widetilde{\delta}_\Gamma} \Psi_1^k (c) &= {\widetilde{\delta}_\Gamma} (-Dh)^k c - {\widetilde{\delta}_\Gamma} h (-D h)^{k-1} D c - {\widetilde{\delta}_\Gamma} h (-D h)^k {\widetilde{\delta}_\Gamma} c \\
&= (-Dh)^k {\widetilde{\delta}_\Gamma} c + (-Dh)^{k-1} D c - {\widetilde{\delta}_\Gamma} h (-Dh)^{k-1} D c - {\widetilde{\delta}_\Gamma} h (-Dh)^k {\widetilde{\delta}_\Gamma} c \\
&= (1- {\widetilde{\delta}_\Gamma} h) (-Dh)^k {\widetilde{\delta}_\Gamma} c + (1- {\widetilde{\delta}_\Gamma} h) (-Dh)^{k-1} D c \\
&= (1- {\widetilde{\delta}_\Gamma} h) (-Dh)^{k-1} (-Dh {\widetilde{\delta}_\Gamma} + D) c \\
&= (1- {\widetilde{\delta}_\Gamma} h) (-Dh)^{k-1} (D (1- h {\widetilde{\delta}_\Gamma} )) c \\
&= ( h {\widetilde{\delta}_\Gamma} ) (-Dh)^{k-1} (D  {\widetilde{\delta}_\Gamma} h ) c \\
\intertext{and by applying the lemma again}
&= h ( (-Dh)^{k-1} {\widetilde{\delta}_\Gamma} + (-Dh)^{k-2} D) (D {\widetilde{\delta}_\Gamma} h ) c \\
&= h ( (-Dh)^{k-1} {\widetilde{\delta}_\Gamma} (D {\widetilde{\delta}_\Gamma} h ) c  + h (-Dh)^{k-2} D) (D {\widetilde{\delta}_\Gamma} h ) c \\
&=0
\end{align*}
by $D^2 = 0$ and $({\widetilde{\delta}_\Gamma})^2 = 0$ as well as commuting $D$ and ${\widetilde{\delta}_\Gamma}$. The cases ${\widetilde{\delta}_\Gamma} \Psi_1^1 (c) = 0$ and ${\widetilde{\delta}_\Gamma} \Psi_1^0 (c) = 0$ are straightforward. Since ${\widetilde{\delta}_\Gamma} \Psi_1 (c) = 0$ we see that $\Psi_1 (c)$ is $\Gamma$-invariant.

Next we will check that the boundary operators commute with $\Psi_1$. Given any $c \in C^k ( \Gamma , \overline{C}^{n} (C^\infty_c (M , \End \E))  [u^{-1} , u] )$
\begin{align*}
\Psi_1 ( ({\widetilde{\delta}_\Gamma} + D) c) &= \Psi_1^k (D c) + \Psi_1^{k+1} ({\widetilde{\delta}_\Gamma} c) \\
&= (-Dh)^k (D c) - h (-Dh )^{k-1} D (D c) - h (-Dh)^k {\widetilde{\delta}_\Gamma} (D c) \\
&\quad + (-Dh)^{k+1} ({\widetilde{\delta}_\Gamma} c) - h (-Dh )^{k} D ({\widetilde{\delta}_\Gamma} c) - h (-Dh)^{k+1} {\widetilde{\delta}_\Gamma} ({\widetilde{\delta}_\Gamma} c) \\
\intertext{which by the commutation relations for the boundary operators}
&= (-Dh)^k D c + (-Dh)^{k+1} {\widetilde{\delta}_\Gamma} c \\
&= -Dh (-Dh)^{k-1} D c - Dh (-Dh)^k {\widetilde{\delta}_\Gamma} c \\
&= D ( (-Dh)^k c - h (-Dh )^{k-1} D c - h(-Dh)^k {\widetilde{\delta}_\Gamma} c) \\
&= D \Psi_1^k (c) = D \Psi_1 (c)
\end{align*}
Hence $\Psi_1$ is a morphism.

\end{proof}

\end{thm}

\begin{thm}
Given a gerbe $(L, \mu)$ on $M \rtimes \Gamma$ with an $L$-twisted vector bundle $\E$ there is a morphism 
\begin{equation}
\Psi_2 : \overline{C}^\bullet (C^\infty_c (M , \End \E)  [u^{-1} , u], b+uB)^\Gamma \to \overline{C}^\bullet (C^\infty_c (M \rtimes \Gamma , L)  [u^{-1} , u], b+uB)
\end{equation}
from the periodic cyclic complex of $\Gamma$-invariant cochains on $C^\infty_c (M , \End \E)$ to the periodic cyclic complex of the twisted convolution algebra $C^\infty_c (M \rtimes \Gamma, L)$.  The morphism is defined by the formula
\begin{multline}
\label{eq:alg_morphism2}
\Psi_2 (c) (\widetilde{a_{g_0}} , \ldots , a_{g_n} ) \\= c(E_{1,g_0} (\widetilde{a_{g_0}}) , \ldots , E_{g_0 \cdots g_{i-1} , g_0 \cdots g_i} (a_{g_i}^{g_0 \cdots g_{i-1}}) , \ldots , E_{g_0 \cdots g_{n-1}, 1} (a_{g_n}^{g_0 \cdots g_{n-1}}))
\end{multline}
for any $g_0 , \ldots g_n \in \Gamma$. Here $c$ is a cochain in $\overline{C}^n (C_c^\infty (M , \End \E))^\Gamma$ and $a_{g_i}$ is a section in $C^\infty_c (M_{g_i} , L_{g_i} )$ for $0 \leq i \leq n$, and $\widetilde{a_{g_0}} = (a_{g_0} , \lambda)$ for $\lambda \in \C$.
\begin{proof}
Note that $E_{1, g_0} (\widetilde{a_{g_0}}) = \widetilde{E_{1,g_0} (a_{g_0})}$ in $\widetilde{C^\infty_c (M , \End \E)}$.  Since $a \in C^\infty_c (M \rtimes \Gamma, L)$ may be written $a = \sum_{g \in \Gamma} a_g$ where $a_g \in C^\infty_c (M_g , L_g)$ for $g \in \Gamma$, it is sufficient to define $\Psi_2$ on such sections.  Also note that $(a_{g} * a_{g'} ) \in C^\infty_c (M_{gg'} , L_{gg'})$ since 
\begin{align*}
(a_g * a_{g'}) (x, h) &= \sum_{h_1 h_2 = h} a_g ( x , h_1) \cdot a_{g'} (xh_1 , h_2) \\
&=a_g ( x , g) \cdot a_{g'} (x g , g') \\
&= \mu_{g , g'} (a_g \otimes a_{g'}^g) (x,g g').
\end{align*}
As the coboundary $b^n = \sum_{i=1}^n d^i_n$, we will prove that $d^i \Psi_2 = \Psi_2 d^i$ for $1\leq i \leq n$.  When $1 \leq i \leq n$
\begin{align*}
&(d^i \Psi_2) ( c ) (\widetilde{a_{g_0}} , \ldots , a_{g_{n+1}} ) \displaybreak[2] \\
&= \Psi_2 (c) (\widetilde{a_{g_0}} , \ldots , a_{g_i} * a_{g_{i+1}} , \ldots , a_{g_{n+1}}) \displaybreak[2] \\
&=\Psi_2 (c) \left(\widetilde{a_{g_0}} , \ldots , \mu_{g_i , g_{i+1}} (a_{g_i} \otimes a_{g_{i+1}}^{g_i}) , \ldots , a_{g_{n+1}}\right) \displaybreak[2] \\
&=c \left( E_{1, g_0} (\widetilde{a_{g_0}})  , \ldots , E_{g_0 \cdots g_{i-1} , g_0 \cdots g_{i+1}} ((\mu_{g_i , g_{i+1}} (a_{g_i} \otimes a_{g_{i+1}}^{g_i} ))^{g_0 \cdots g_{i-1} } ) , \ldots \right. \\
&\quad \left.\ldots, E_{g_0 \cdots g_n , 1} (a_{g_{n+1}}^{g_0 \cdots g_n}) \right)\displaybreak[2] \\
&=c \left( E_{1, g_0} (\widetilde{a_{g_0}})  , \ldots , E_{g_0 \cdots g_{i-1} , g_0 \cdots g_{i+1}} (\mu_{g_i , g_{i+1}}^{g_0 \cdots g_{i-1}} (a_{g_i}^{g_0 \cdots g_{i-1} } \otimes a_{g_{i+1}}^{g_0 \cdots g_i } ) ) , \ldots \right. \\
&\quad \left. \ldots, E_{g_0 \cdots g_n , 1} (a_{g_{n+1}}^{g_0 \cdots g_n}) \right) \displaybreak[2] \\
&= c \left( E_{1, g_0} (\widetilde{a_{g_0}})  , \ldots , E_{g_0 \cdots g_{i-1} , g_0 \cdots g_i} (a_{g_i}^{g_0 \cdots g_{i-1}} ) E_{g_0 \cdots g_i , g_0 \cdots g_{i+1}} (a_{g_{i+1}}^{g_0 \cdots g_i} ) , \ldots \right. \\
&\quad \left. \ldots, E_{g_0 \cdots g_n , 1} (a_{g_{n+1}}^{g_0 \cdots g_n}) \right) \displaybreak[2] \\
&= (d^i c) \left(E_{1, g_0} (\widetilde{a_{g_0}})  , \ldots , E_{g_0 \cdots g_{i-1} , g_0 \cdots g_i} (a_{g_i}^{g_0 \cdots g_{i-1}} ) , E_{g_0 \cdots g_i , g_0 \cdots g_{i+1}} (a_{g_{i+1}}^{g_0 \cdots g_i} ) , \ldots \right. \\
&\quad \left. \ldots, E_{g_0 \cdots g_n , 1} (a_{g_{n+1}}^{g_0 \cdots g_n}) \right) \displaybreak[2] \\
&= (\Psi_2 d^i) (c) (\widetilde{a_{g_0}} , \ldots , a_{g_{n+1}} ).
\end{align*}
For the case of $d^0 \Psi_2 = \Psi_2 d^0$ note that $\widetilde{a_{g_0}} * a_{g_1} = a_{g_0} * a_{g_1} + \lambda a_{g_1}$ so that 
\begin{align*}
&(d^0 \Psi_2) ( c ) (\widetilde{a_{g_0}} , \ldots , a_{g_{n+1}} ) \displaybreak[2] \\
&= \Psi_2 (c) (\widetilde{a_{g_0}} * a_{g_1}, a_{g_2} , \ldots  , a_{g_{n+1}}) \displaybreak[2] \\
&= \Psi_2 (c) (a_{g_0} * a_{g_1} + \lambda a_{g_1}, a_{g_2} , \ldots  , a_{g_{n+1}}) \displaybreak[2] \\
&= \Psi_2 (c) (a_{g_0} * a_{g_1} , a_{g_2} , \ldots  , a_{g_{n+1}}) + \Psi_2 (c) ( \lambda a_{g_1}, a_{g_2} , \ldots  , a_{g_{n+1}}) \displaybreak[2] \\
&= \Psi_2 (c) (\mu_{g_0 , g_1} (a_{g_0} \otimes a_{g_1}^{g_0} ) , a_{g_2} , \ldots  , a_{g_{n+1}}) + \Psi_2 (c) ( \lambda a_{g_1},  a_{g_2} , \ldots  , a_{g_{n+1}}) \displaybreak[2] \\
&= c ( E_{1, g_0 g_1} (\mu_{g_0 , g_1} (a_{g_0} \otimes a_{g_1}^{g_0} )) , E_{g_0 g_1 , g_0 g_1 g_2} (a_{g_2}^{g_0 g_1}) , \ldots , E_{g_0 \dots g_n , 1} (a_{g_{n+1}}^{g_0 \cdots g_n})) \\
&\quad + c (\lambda E_{1, g_1} (a_{g_1}) , E_{g_1 , g_1 g_2} (a_{g_2}^{g_1}) , \ldots , E_{g_1 \cdots g_n ,  g_0^{-1} } (a_{g_{n+1}}^{g_1 \cdots g_n})) \displaybreak[2] \\
\intertext{as the second term is nonzero only if $g_0 \cdots g_{n+1} =1$ and $g_1 \cdots g_{n+1} = 1$, i.e.\ $g_0=1$, }
&= c \left( E_{1,g_0} (a_{g_0}) E_{g_0 , g_0 g_1} (a_{g_1}^{g_0}) , E_{g_0 g_1 , g_0 g_1 g_2} (a_{g_2}^{g_0 g_1}),  \ldots \right. \\
&\quad \left. \ldots , E_{g_0 \cdots g_n , 1} (a_{g_{n+1}}^{g_0 \cdots g_n})  \right)\\
&\quad + c \left( \lambda E_{g_0 , g_0 g_1} (a_{g_1}^{g_0}) , E_{g_0 g_1 , g_0 g_1 g_2} (a_{g_2}^{g_0 g_1}),  \ldots , E_{g_0 \cdots g_n , 1} (a_{g_{n+1}}^{g_0 \cdots g_n})  \right) \displaybreak[2] \\
&= c \left( E_{1,g_0} (a_{g_0}) E_{g_0 , g_0 g_1} (a_{g_1}^{g_0}) + \lambda E_{g_0 , g_0 g_1} (a_{g_1}^{g_0}), E_{g_0 g_1 , g_0 g_1 g_2} (a_{g_2}^{g_0 g_1}),  \ldots \right. \\
&\quad \left. \ldots, E_{g_0 \cdots g_n , 1} (a_{g_{n+1}}^{g_0 \cdots g_n})  \right) \displaybreak[2] \\
&= c \left( E_{1,g_0} (\widetilde{a_{g_0}}) E_{g_0 , g_0 g_1} (a_{g_1}^{g_0}), E_{g_0 g_1 , g_0 g_1 g_2} (a_{g_2}^{g_0 g_1}),  \ldots , E_{g_0 \cdots g_n , 1} (a_{g_{n+1}}^{g_0 \cdots g_n}) \right) \displaybreak[2] \\
&= (d^0 c) \left( E_{1,g_0} (\widetilde{a_{g_0}}) , E_{g_0 , g_0 g_1} (a_{g_0}^{g_1}),  \ldots , E_{g_0 \cdots g_n , 1} (a_{g_{n+1}}^{g_0 \cdots g_n}) \right) \displaybreak[2] \\
&= (\Psi_2 d^0) (c) (\widetilde{a_{g_0}} , \ldots , a_{g_{n+1}} ) .
\end{align*}

We also have $d^{n+1} \Psi_2 = \Psi_2 d^{n+1}$ as
\begin{align*}
&(d^{n+1} \Psi_2) ( c ) (a_{g_0} , \ldots , a_{g_{n+1}} ) \\
&= \Psi_2 (c) (a_{g_{n+1}} * a_{g_0}  , \ldots , a_{g_{n}}) \\
&= c (E_{1 , g_{n+1} g_0} (\mu (a_{g_{n+1}} \otimes a_{g_0}^{g_{n+1}}) ) , E_{g_{n+1} g_0 , g_{n+1} g_0 g_1} (a_{g_1}^{g_{n+1} g_0}) , \ldots \\
&\quad \ldots, E_{g_{n+1} g_0 \cdots g_{n-1} , 1} (a_{g_n}^{g_{n+1} g_0 \cdots g_{n-1}})) \\
&=  c(E_{1 , g_{n+1} g_0} (\mu (a_{g_{n+1}} \otimes a_{g_0}^{g_{n+1}}) ) , E_{g_{n+1} g_0 , g_{n+1} g_0 g_1} (a_{g_1}^{g_{n+1} g_0}) , \ldots \\
&\quad \ldots, E_{g_{n+1} g_0 \cdots g_{n-1} , 1} (a_{g_n}^{g_{n+1} g_0 \cdots g_{n-1}})) \cdot g_{n+1}^{-1}  \\
&= c(E_{g_0 \cdots g_n , g_0} (\mu (a_{g_{n+1}} \otimes a_{g_0}^{g_{n+1}}) ) , E_{g_0 , g_0 g_1} (a_{g_1}^{g_0}) , \ldots , E_{g_0 \cdots g_{n-1} , g_0 \cdots g_n} (a_{g_n}^{g_0 \cdots g_{n-1}}) ) \\
&=  c(E_{g_0 \cdots g_n , 1} (a_{g_{n+1}}) E_{1, g_0} ( a_{g_0}^{g_{n+1}} ) , E_{g_0 , g_0 g_1} (a_{g_1}^{g_0}) , \ldots , E_{g_0 \cdots g_{n-1} , g_0 \cdots g_n} (a_{g_n}^{g_0 \cdots g_{n-1}}) ) \\
&= (d^{n+1} c) ( E_{1, g_0} ( a_{g_0}^{g_{n+1}} )) , E_{g_0 , g_0 g_1} (a_{g_1}^{g_0}) , \ldots , E_{g_0 \cdots g_n , 1} (a_{g_{n+1}}^{g_0 \cdots g_n}) ) \\
&= (\Psi_2  d^{n+1}) (c) (a_{g_0} , \ldots , a_{g_{n+1}} )
\end{align*}
and the extension to $\widetilde{a_{g_0}} = (a_{g_0} , \lambda)$ is not much different from the previous calculation.  This shows that $b \Psi_2 = \Psi_2 b$.

Next we show that $B \Psi_2 = \Psi_2 B$.  As $B = \sum_i \bar{s}^i$ we will show that $\bar{s}^i \Psi_2 = \Psi_2 \bar{s}^i$ for $0 \leq i \leq n$.  The identity element in $C^\infty_c (M \rtimes \Gamma , L)$ is $(0,1)$ where $0$ may be considered as the sum of zero sections $0 : M \rightarrow L_g$ for any $g$, so we consider $0$ as $0 : M \rightarrow L_1$ below.  
\begin{align*}
& (\Psi_2 \bar{s}^i) (c) (a_{g_0} , \ldots , a_{g_n}) \\
&= (\bar{s}^i c) \left(E_{1, g_0} (a_{g_0}) , \ldots , E_{g_0 \cdots g_{i-1} , g_0 \cdots g_i} (a_{g_i}^{g_0 \cdots g_{i-1}}) , \ldots , E_{g_0 \cdots g_{n-1} , 1} (a_{g_n}^{g_0 \cdots g_{n-1}}) \right) \\
&= (\bar{s}^i c) \left(E_{1, g_0} (a_{g_0}) , \ldots , E_{g_0 \cdots g_{i-1} , g_0 \cdots g_i} (a_{g_i}^{g_0 \cdots g_{i-1}}) , \ldots , E_{g_0 \cdots g_{n-1} , 1} (a_{g_n}^{g_0 \cdots g_{n-1}}) \right) \cdot (g_i \cdots g_n) \\
&= (\bar{s}^i c) \left(E_{g_i \cdots g_n, g_i \cdots g_n g_0} (a_{g_0}^{g_i \cdots g_n}) , \ldots , E_{g_i \cdots g_n g_0 \cdots g_{i-2} , 1} (a_{g_{i-1}}^{g_i \cdots g_n g_0 \cdots g_{i-2}}), \right. \\
&\qquad \qquad \left. E_{1 , g_i } (a_{g_{i}}) , \ldots , E_{g_i \cdots g_{n-1} , g_i \cdots g_n} (a_{g_n}^{g_i \cdots g_{n-1}}) \right)  \\
&= c \left(1, E_{1,g_i} (a_{g_i}) , E_{g_i , g_i g_{i+1}} (a_{g_{i+1}}^{g_i}) , \ldots , E_{g_i \cdots g_{n-1} , g_i \cdots g_n} (a_{g_n}^{g_i \cdots g_n}), \right. \\ 
&\qquad \qquad \left. E_{g_i \cdots g_n , g_i \cdots g_n g_0} (a_{g_0}^{g_i \cdots g_n}) , \ldots , E_{g_i \cdots g_n g_0 \cdots g_{i-2} , 1} (a_{g_{i-1}}^{g_i \cdots g_n g_0 \cdots g_{i-2}})    \right) \\
&= (\Psi_2 c) (1, a_{g_i} , \ldots , a_{g_n} , a_{g_0} , \ldots , a_{g_i}) \\
&= (\bar{s}^i \Psi_2) (c) (a_{g_0} , \ldots , a_{g_n}). 
\end{align*}
Hence $\Psi_2$ is a morphism.
\end{proof}
\end{thm}

We have now proved the main theorem.

\begin{thm}
Given a gerbe $(L , \mu)$ on $M \rtimes \Gamma$ with connection $\nabla$ there is a morphism
\begin{multline}
\Psi_2 \circ \Psi_1 \circ \Psi_0 \circ \tau_\nabla : (\Omega^*_- ((M\rtimes \Gamma) _\bullet) [u] , u\tilde{d}_{dR} + d_\Delta - \Theta_u \wedge \cdot) \\ \to \overline{C}^\bullet (C^\infty_c (M \rtimes \Gamma , L)  [u^{-1} , u], b+uB)
\end{multline}
from the twisted simplicial complex of $M \rtimes \Gamma$ to the $(b,B)$ complex of the twisted convolution algebra $C^\infty_c ( M \rtimes \Gamma , L)$.
\end{thm}

\bibliographystyle{plain}	


\bibliography{refs}		

 \end{document}